\documentclass{article}

\usepackage{microtype}
\usepackage{graphicx}
\usepackage{subfigure}
\usepackage{booktabs} 
\usepackage[usenames]{color}
\usepackage{enumitem}

\usepackage{hyperref}

\usepackage[accepted]{icml2023}

\usepackage{amsmath}
\usepackage{amssymb}
\usepackage{mathtools}
\usepackage{amsthm}
\usepackage[capitalize,noabbrev]{cleveref}
\usepackage{natbib}

\theoremstyle{plain}
\newtheorem{theorem}{Theorem}[section]

\newtheorem{lemma}[theorem]{Lemma}
\newtheorem{corollary}[theorem]{Corollary}
\theoremstyle{definition}

\newtheorem{assumption}[theorem]{Assumption}
\theoremstyle{remark}

\usepackage[textsize=tiny]{todonotes}

\usepackage{caption}
\icmltitlerunning{Constrained Optimization via Exact Augmented Lagrangian and Randomized Iterative Sketching}
\usepackage{math_command}

\begin{document}

\twocolumn[
\icmltitle{Constrained Optimization via Exact Augmented Lagrangian \\ and Randomized Iterative Sketching}

\icmlsetsymbol{equal}{*}
\begin{icmlauthorlist}
\icmlauthor{Ilgee Hong}{1,equal}
\icmlauthor{Sen Na}{2,3,equal}
\icmlauthor{Michael W. Mahoney}{2,3,4}
\icmlauthor{Mladen Kolar}{5}
\icmlaffiliation{1}{Department of Statistics, University of Chicago}
\icmlaffiliation{2}{Department of Statistics, University of California, Berkeley}
\icmlaffiliation{3}{International Computer Science Institute}
\icmlaffiliation{4}{Lawrence Berkeley National Laboratory}
\icmlaffiliation{5}{Booth School of Business, University of Chicago}
\icmlcorrespondingauthor{Sen Na}{senna@berkeley.edu}
\end{icmlauthorlist}
\icmlkeywords{Inexact Newton method, Constrained optimization, Randomized numerical linear algebra, Sketching, Augmented Lagrangian}
\vskip 0.1in
]

\printAffiliationsAndNotice{\icmlEqualContribution}

\begin{abstract}

We consider solving equality-constrained nonlinear, nonconvex optimization problems. This class of problems appears widely in a variety of applications in machine learning and engineering,~ranging from constrained deep neural networks, to optimal control, to PDE-constrained optimization.~We develop an adaptive inexact Newton method for this problem class. In each iteration, we solve the Lagrangian Newton system inexactly via a~\textit{randomized iterative sketching} solver, and select a suitable stepsize by performing line search on an \textit{exact augmented Lagrangian} merit function. The randomized solvers have advantages over deterministic linear system solvers by significantly reducing per-iteration flops complexity and storage cost, when equipped with suitable sketching matrices. Our method adaptively controls the accuracy of the randomized solver and the penalty parameters of the exact augmented Lagrangian, to ensure that the inexact Newton direction is a descent direction of the exact augmented Lagrangian. This allows us to establish a \textit{global almost sure convergence}. We also show that a unit stepsize is admissible locally, so that our method exhibits a \textit{local linear convergence}. Furthermore, we prove that the linear convergence can be strengthened to \textit{superlinear convergence} if we gradually sharpen the adaptive accuracy condition on the randomized solver. We demonstrate the superior performance of our method on benchmark nonlinear problems in CUTEst test set, constrained logistic regression with data from LIBSVM, and~a~PDE-constrained problem.

\end{abstract}

\section{Introduction}\label{sec:1}

We study constrained optimization problems of the form\vskip-13pt
\begin{equation}\label{eq:opt_main}
\min_{\vx\in \mR^n}\;f(\vx)\quad\quad \text{s.t.}\;\; c(\vx) = \vzero,
\end{equation}\vskip-7pt
where $f:\mR^n\rightarrow \mR$ is the objective function,  $c:\mR^n\rightarrow\mR^m$ is the equality constraint function, and both of them~are~nonlinear, possibly nonconvex, and twice continuously differentiable. Problem \eqref{eq:opt_main} appears widely (as a subproblem) in a variety of applications, including constrained deep neural networks \citep{Chen2018Constraint}, physical informed neural~networks \citep{Cuomo2022Scientific}, PDE-constrained optimization \citep{Kouri2014Inexact}, optimal control \citep{Lewis2012Optimal}, and constrained model estimations \citep{Scott2022Solving}. There exist numerous methods for solving Problem \eqref{eq:opt_main}, such as projected first- and second-order methods, penalty methods, augmented Lagrangian methods, and sequential quadratic programming (SQP) methods \citep{Nocedal2006Numerical}. Among these methods, SQP is one of the most effective methods for both small and large problems. Compared to projected methods, SQP does not involve a projection step that is generally expensive for nonlinear equality constraints. Compared to penalty methods, SQP does not regularize the objective that destroys the structure of the problem and suffer from ill-conditioning issue \citep{Krishnapriyan2021Characterizing}. Compared to augmented Lagrangian methods, SQP is more robust to the (dual) initialization \citep{Curtis2014adaptive}.

With only equality constraints, SQP is equivalent to an application of Newton method to the Karush–Kuhn–Tucker (KKT) conditions (i.e., first-order optimality conditions).~As the second-order method, SQP is (not surprisingly) very efficient, taking only a few iterations to find a local solution to the problem; while, at the same time, it can struggle~to~solve the Newton system (of the Lagrangian) efficiently at each step. Solving Newton systems dominates the computational cost of Newton methods, and deriving exact~Newton directions is often prohibitive for large-scale problems. 

\subsection{Main Results}\label{sec:1.1}

Motivated by the above bottleneck, we propose an adaptive inexact Newton (SQP) method for Problem \eqref{eq:opt_main} that applies randomized iterative sketching solvers to solve Newton systems inexactly.~Our method is called \textbf{$\myalg$}. At~each~step,~$\myalg$~solves~the~\mbox{Lagrangian}~Newton system by leveraging the \textit{sketch-and-project} framework, which was originally developed in \citet{Gower2015Randomized} and further investigated for various (unconstrained) optimization problems \citep{Luo2016Efficient, Doikov2018Randomized, Gower2019Rsn, Mutny2020Convergence, Derezinski2020Precise, Na2022Asymptotic, Derezinski2022Sharp}. The framework~unifies~popular~randomized iterative linear~system~solvers,~such~as~the~randomized (block) coordinate descent and randomized (block)~Kaczmarz method. Compared to deterministic solvers,~randomized solvers require less per-iteration flops complexity and storage cost when using suitable~sketching matrices (e.g., a sparse sketching matrix). See \citet{Murray2023Randomized} for recent software development. However, using inexact 
randomized solver leads to a much more involved algorithm
design as well as a more challenging analysis, since the approximation error of the solution is random and does not monotonically decrease as more iterations of the solver are performed. As we review in Section~\ref{sec:1.2}, existing inexact Newton methods with randomized sketching solvers are designed for unconstrained problems; we provide the~first~inexact Newton method with a randomized~solver~for constrained~problems.

In the proposed $\myalg$ algorithm, we use an~\textit{exact augmented Lagrangian merit function}, taking the form \vskip-18pt
\begin{equation}\label{eq:opt_merit}
\resizebox{0.89\columnwidth}{!}{$\mathcal{L}_{\veta}(\vx,\vlambda)=\mathcal{L}(\vx,\vlambda)+\dfrac{\eta_1}{2}\Vert c(\vx)\Vert^2+\dfrac{\eta_2}{2}\Vert\nabla_{\vx}\mathcal{L}(\vx,\vlambda)\Vert^2$},
\end{equation}\vskip-10pt
where $\mathcal{L}(\vx,\vlambda)=f(\vx)+\vlambda^Tc(\vx)$ is the Lagrangian function of \eqref{eq:opt_main}, with $\vlambda\in\mR^m$ being the Lagrangian multipliers and $\veta=(\eta_1,\eta_2)$ being the penalty parameters. By ``exact,'' we mean that, compared with standard augmented Lagrangian, it has an additional penalty term $\eta_2\|\nabla_{\vx}\mL\|^2/2$ that biases the optimality error. With this additional term, one can show that the solution of the unconstrained problem $\min_{\vx, \vlambda}\mL_{\veta}$~is also the solution of \eqref{eq:opt_main} and vice versa, provided that $\veta$ are suitably specified. By ``merit function," we mean \eqref{eq:opt_merit} is (only) used to determine whether or not a new iterate approaches toward a local solution to Problem \eqref{eq:opt_main}. In other words, \eqref{eq:opt_merit} does not affect the computation of inexact Newton direction (thus, our method does not suffer ill-conditioning issue caused by $\veta$); but given a direction, we select a suitable stepsize by performing line search on \eqref{eq:opt_merit}. The merit function plays a crucial role in constrained optimization, since the objective function $f$ alone is not suitable to justify the direction (the step that decreases $f$ may severely violate the constraint $c$). 
$\myalg$ adaptively controls the accuracy of the randomized sketching solver and selects suitable penalty parameters $\veta$ to ensure~that the inexact Newton direction is a descent direction of \eqref{eq:opt_merit}. Our adaptive design~balances the number of global outer loop iterations and  inner loop iterations of the solver, while at the same time securing a fast local convergence. By virtue of the adaptivity, all the input parameters of the algorithm are problem independent.

Under mild assumptions, we show that $\myalg$ enjoys a global almost sure convergence --- starting from~any initial point, the KKT residual converges to zero. We also show that our adaptive design ensures a unit stepsize to be admissible locally, and establish a local linear~convergence rate. Furthermore, we prove that the linear rate~can be strengthened to a superlinear rate, as long as we gradually sharpen the adaptive accuracy condition on the randomized solver, which is simply controlled by an input tuning~parameter. As we review in Section \ref{sec:1.2}, the local rate of inexact Newton methods is mostly investigated for unconstrained problems \citep{Gower2019Rsn, Gower2021Adaptive, Hanzely2020Stochastic, Derezinski2022Sharp, Yuan2022Sketched}, while is largely missing in constrained cases \citep{Byrd2008Inexact,Byrd2010inexact, Curtis2014Inexact, Gu2017new, Burke2020Inexact}.

We implement $\myalg$~and~\mbox{benchmark}~\mbox{it against} \mbox{constrained} nonlinear problems in CUTEst test set \citep{Gould2014CUTEst}, constrained logistic regression with datasets from LIBSVM \citep{Chang2011LIBSVM}, and a PDE-constrained~problem~\citep{Hintermueller2002Primal}. We show the~superior performance of our method, in terms of both accuracy and efficiency, by comparing it with a~prominent inexact Newton method~that employs a deterministic solver and $\ell_1$ merit function \citep{Byrd2008Inexact} and standard augmented Lagrangian method \citep{Nocedal2006Numerical}. We also show our method is robust to tuning parameters.

\subsection{Related Literature}\label{sec:1.2}

Our method relates to several inexact Newton-type methods in the literature that we briefly review below. We divide the review into unconstrained and constrained optimization.

\textbf{Unconstrained optimization}. The majority of inexact Newton methods are designed for unconstrained problems. For example, some inexact Newton methods via randomized~iterative sketching have been proposed,~including Stochastic Dual Newton Ascent (SDNA) \citep{Qu2016SDNA}, Sketched Online Newton (SON) \citep{Luo2016Efficient}, Stochastic Subspace Cubic Newton (SSCN) \citep{Hanzely2020Stochastic}, Randomized Block Cubic Newton (RBCN) \citep{Doikov2018Randomized}, Randomized Subspace Newton (RSN) \citep{Gower2019Rsn}, and Randomized Subspace Regularized Newton (RSRN) \citep{Fuji2022Randomized}. These methods incorporate sketching techniques into classical Newton methods or cubic regularized Newton methods. Instead of solving the original large-scale subproblem in each iteration, the methods solve a sketched small-scale subproblem, and it is proved that the sketched direction decreases the objective, due to the convexity assumption and/or unconstrained~nature of the problem. With a good understanding on the~convergence of sketching solvers \citep{Gower2021Adaptive, Derezinski2022Sharp}, the (local) convergence rates of these methods are also established. Empirical results illustrate that~randomized Newton methods integrate the benefits~of~randomization and the efficiency of second-order methods. The proposed $\myalg$ is the first generalization of the aforementioned methods to constrained nonlinear problems \eqref{eq:opt_main}. Our method adaptively controls the accuracy of the sketching solver and the penalty parameters of the merit function, and takes the constraint violation into account when showing the descent property of the sketched inexact Newton direction.

 There are numerous Newton methods where subproblems can be solved inexactly~by~deterministic solvers such as MINRES and CG \citep{Fong2012CG}. For example, different variants of Newton-CG methods have been developed with line search \citep{Royer2018Complexity, Royer2019Newton, Bollapragada2018Exact, Yao2022Inexact}, trust region \citep{Curtis2021Trust}, and cubic regularization \citep{Curtis2018inexact, Yao2021Inexact}. Newton-MR has also recently been reported and analyzed \citep{Liu2022Newton, Roosta2022Newton}. Another line of work focuses on \mbox{constructing} cheap stochastic Hessian approximations by either sketching or subsampling \citep{Pilanci2017Newton, RoostaKhorasani2018Sub, Derezinski2020Debiasing, Lacotte2021Adaptive, Derezinski2021Newton, Na2022Hessian}. We refer to \citet{Berahas2020investigation} for a comparison of these methods. For some of the methods (e.g., \citet{RoostaKhorasani2018Sub}), the resulting Newton systems can also be solved inexactly via deterministic solvers. However, the convergence guarantees of those methods hold in high probability or expectation and (more or less) rely on the convexity of problems. This differs from our almost sure convergence guarantee for nonlinear problems.

\textbf{Constrained optimization}. To our knowledge, there are only a handful of inexact Newton methods for constrained problems, which all employ deterministic solvers for solving Newton systems and are often called inexact SQP methods. \citet{Byrd2008Inexact, Byrd2010inexact} proposed two inexact methods for Problem \eqref{eq:opt_main}, which adopt nonsmooth merit functions in line search step. The methods bound the residuals of the solver by a few \textit{fixed} tuning parameters, and are then refined and adapted to incorporate inequality constraints \citep{Curtis2014Inexact, Burke2020Inexact}. Recently, stochastic SQP methods have been reported \citep{Berahas2021stochastic, Berahas2021Sequential, Curtis2021Inexact, Na2023Inequality, Na2022adaptive, Fang2022Fully, Na2022Asymptotic}. The \mbox{stepsize} in~most~of these~\mbox{methods} is controlled by prespecified decaying sequences (i.e., the~methods are designed under the stochastic approximation regime), with the only exceptions being \citet{Na2023Inequality} and \citet{Na2022adaptive}, which adopt line search to make the methods more adaptive. The aforementioned methods established global convergence, while local rates of the methods are largely missing. Our method~is~the~first design that incorporates~the sketching technique into SQP, adaptively selects suitable stepsizes by line search, and exhibits global and local linear/superlinear convergence as exact SQP/Newton methods do.

\vskip-4pt
\textbf{Notation.} We use $\|\cdot\|$ to denote the $\ell_2$ norm for vectors and the operator norm for matrices. For any matrix $A\in\mR^{p\times q}$, $row(A)$ denotes the row space of $A$ and $col(A)$ denotes the column space of $A$. We let $H(\vx,\vlambda)=\nabla^2_{\vx}\mathcal{L}(\vx,\vlambda)\in\mR^{(n+m)\times(n+m)}$ be the Lagrangian Hessian with respect to $\vx$ and $G(\vx)=\nabla^Tc(\vx)\in\mR^{m\times n}$ be the constraint Jacobian. At the $k$-th iteration, we let $f_k = f(\vx_k)$, $H_k = H(\vx_k,\vlambda_k)$, etc., to ease the notation. For any $a,b\in\mR$, we let $a\vee b=\max\{a,b\}$ and $a\wedge b=\min\{a,b\}$.

\section{$\myalg$}\label{sec:2}

We now introduce $\myalg$ algorithm. For Problem \eqref{eq:opt_main}, under certain constraint qualifications, a necessary condition for $\tx$ being a local solution is to have dual multipliers $\tlambda$ such that
\begin{equation}\label{eq:kkt}
\resizebox{0.89\columnwidth}{!}{$\begin{pmatrix}\nabla_{\vx}\mathcal{L}(\vx^{\star},\vlambda^{\star})\\\nabla_{\vlambda}\mathcal{L}(\vx^{\star},\vlambda^{\star})\end{pmatrix}=\begin{pmatrix}\nabla f(\vx^{\star})+G(\vx^{\star})^T\vlambda^{\star}\\c(\vx^{\star})\end{pmatrix}=\begin{pmatrix}\vzero\\\vzero\end{pmatrix}$}.
\end{equation}
The starting point of our~\mbox{algorithm is to apply the Newton} method to \eqref{eq:kkt}. At each iteration $k$, we consider solving the following Newton system of the Lagrangian
\begin{equation}\label{eq:kkt_sys_reduced}
\KKT\inczk = - \nablal,
\end{equation}	
where \vskip-20pt
\begin{equation*}
\resizebox{0.89\columnwidth}{!}{$\KKT=\kkt,\;\inczk=\begin{pmatrix}
\Delta\vx_k\\
\Delta\vlambda_k
\end{pmatrix},\; \nablal = \begin{pmatrix}
\nabla_{\vx}\mL_k\\
c_k
\end{pmatrix}$}.
\end{equation*}\vskip-10pt
Here, $\KKT$ approximates the Lagrangian Hessian $\nabla^2\mL_k$ with $B_k$ being a modification of $H_k$ to ensure \eqref{eq:kkt_sys_reduced} is solvable (cf. Assumption \ref{assum2}). We note that solving the Newton~system \eqref{eq:kkt_sys_reduced} is equivalent to solving a constrained quadratic~program
\begin{align*}
\min_{\Delta\vx} \;\; & \frac{1}{2}\Delta\vx^TB_k\Delta\vx + \nabla f_k^T\Delta\vx,\\
\text{s.t.}\;\; & c_k + G_k\Delta\vx = \vzero.
\end{align*}
The objective is a quadratic approximation of $f$ with Hessian coming from the Lagrangian function, while the constraint is a linear approximation of $c$. Instead of solving \eqref{eq:kkt_sys_reduced} exactly and having $\inczk$, we apply a randomized iterative sketching solver on \eqref{eq:kkt_sys_reduced} and derive an (random) inexact solution $\incztk$.

\subsection{Randomized Iterative Sketching} 

For every \textit{outer iteration} $k$, we run multiple \textit{inner iterations} of the sketching solver, indexed by $j$. The solver will stop at some $j$ and output $\incztk\coloneqq \incztkj = (\tilde{\Delta}\vx_{k,j},\tilde{\Delta}\vlambda_{k,j})$ as the inexact direction, once the step $\incztkj$ triggers an adaptive step acceptance condition (introduced in Section~\ref{sec:2.2}). As~we will prove in Lemmas \ref{lemma10} and \ref{lemma11}, the step acceptance condition is always triggered with a finite $j$ (with~probability~one).

We now provide details of the sketch-and-project framework. With different sketching matrices, the framework recovers different randomized methods including, randomized Kaczmarz method, randomized coordinate descent, etc. \citep{Gower2015Randomized, Strohmer2008Randomized, Richtarik2014Iteration}. Let $S\in\mR^{(n+m)\times d}\sim\mathcal{P}$ be a random sketching matrix following the distribution $\mathcal{P}$ (its sketching dimension $d$ can also be random). For each outer iteration $k$ and inner iteration $j$, we generate a copy of $S$ denoted by $S_{k,j}$, and solve the following sketched Newton system
\begin{equation}\label{eq:sketch_system}
S_{k,j}^T\KKT\vu=-S_{k,j}^T\nablal.
\end{equation}
Since \eqref{eq:sketch_system} has multiple solutions including the exact direction $\inczk$, the $j$-th iteration of the solver selects the step $\incztkjn$ to be the one that is closest to the current step $\incztkj$,~i.e., \vskip-20pt
\begin{equation}
\incztkjn=\arg\underset{\vu}{\min}\; \Vert \vu-\incztkj\Vert^2,\quad\text{subject to}\;\;\eqref{eq:sketch_system}.\nonumber
\end{equation} \vskip-10pt
An explicit updating rule is given by \vskip-17pt
\begin{equation}\label{eq:irs}	\resizebox{0.89\columnwidth}{!}{$\incztkjn=\incztkj-\KKT S_{k,j}(S_{k,j}^T\KKT^2S_{k,j})^{\dagger}S_{k,j}^T\vr_{k,j}$},
\end{equation}\vskip-7pt
where $(\cdot)^{\dagger}$ is the Moore–Penrose pseudoinverse and $\vr_{k,j}$ is the residual of $\incztkj$, defined as \vskip-13pt
\begin{equation}\label{eq:residual}
\vr_{k,j}=\KKT\incztkj+\nablal.
\end{equation}
We initialize the solver with $\tilde{\Delta}\vz_{k,0} = \vzero$, and below we~introduce when we should stop the inner iteration \eqref{eq:irs}.

\subsection{Adaptive Step Acceptance Condition}\label{sec:2.2}

Our step acceptance \mbox{condition}~\mbox{consists}~of~two~\mbox{subconditions}: (i) adaptive accuracy condition enforced on the residual~$\vr_{k,j}$ (cf. \eqref{eq:cond1}); and (ii) descent direction condition~\mbox{enforced} on the step $\incztkj$ (cf.~\eqref{eq:cond2}). Both~\mbox{conditions}~rely~on~the~penalty parameters $\veta_k = (\eta_{1,k}, \eta_{2,k})$ of the augmented Lagrangian \eqref{eq:opt_merit}, which we also adaptively choose. Once the step~acceptance condition is triggered, we derive an inexact direction $\incztk = \incztkj$ with the chosen penalty parameters $\veta_k$. With these chosen quantities, we do line search (in Section \ref{sec:2.3}).

In particular, given the penalty parameters $\veta_k$, we want to simultaneously enforce two (sub)conditions. 

\noindent\textbf{Adaptive accuracy condition.} We first compute a threshold \vskip-10pt
\begin{equation}\label{eq:delta_trial}
\delta_k^{\text{trial}}:=\dfrac{(0.5-\beta)\eta_{2,k}}{(1+\eta_{1,k}+\eta_{2,k})\Upsilon_k^2\Psi_k^2},
\end{equation}
where $\beta\in(0,0.5)$ is an input parameter of the algorithm used in the Armijo condition in line search (cf. \eqref{eq:armijo}), and
\begin{equation}\label{eq:aux}
\Psi_k=\frac{20(\normB^2\vee 1)}{(\xi_{B}\wedge1)(\sigma_{1,k}^2\wedge1)},\quad\Upsilon_k=\normG\vee\normH\vee 1,
\end{equation}
with $\sigma_{1,k}$ being the least singular value of $G_k$ and $\xi_{B}>0$ being another input parameter used for constructing $B_k$ (cf. Assumption \ref{assum2}). With any $0<\delta_k\leq \delta_k^{\text{trial}}$ and any prespecified sequence $\{\theta_k\}\subseteq(0,1]$, we require $\vr_{k,j}$ to~\mbox{satisfy} \vskip-15pt
\begin{equation}\label{eq:cond1}
\left\Vert\vr_{k,j}\right\Vert\le\theta_k\delta_k\normnablal/(\normA\Psi_k).
\end{equation} 
We explain four aspects of the above condition. First, $\normA$ in \eqref{eq:cond1} makes the approximation error $\Vert\incztkj-\inczk\Vert$ relating to $\|\vr_{k,j}\|$ be bounded by a factor of the exact direction $\normexactz$, seen from $\|\nabla\mL_k\| = \|\Gamma_k\Delta\vz_k\| \leq \|\Gamma_k\|\cdot\|\Delta\vz_k\|$. Second, the factor 20 in $\Psi_k$ is a conservative, artificial constant coming from the proof, which may be significantly~reduced by finer analysis. Third, the condition $\delta_k\leq \delta_k^{\text{trial}}$ is to ensure that the algorithm selects a unit stepsize locally. In fact, under standard assumptions, we show in Lemma \ref{lemma12} that $\delta_k^{\text{trial}}$ is uniformly lower bounded away from zero, so that $\delta_k$ does not have to converge to zero; instead, for $k$ large enough, $\delta_k$ stabilizes. Fourth, $\theta_k$ is introduced to enhance the flexibility of the accuracy condition. When $\theta_k = \theta\in(0,1]$, $\forall k$,~we show the algorithm exhibits local linear convergence.~When $\theta_k$ decays to zero, we actually gradually sharpen the accuracy condition \eqref{eq:cond1} (as $\delta_k$ finally stabilizes); and the algorithm exhibits local superlinear \mbox{convergence} with a~rate~depending on $\theta_k$. See Section \ref{sec:4} for the results of different~$\theta_k$.

\noindent\textbf{Descent direction condition.} We require $\incztkj$ to be a descent direction of $\mL_{\veta_k}$ at $(\vx_k,\vlambda_k)$. Specifically, we require \vskip-15pt
\begin{equation}\label{eq:cond2}
(\nablaL)^T\incztkj\le-\eta_{2,k}\normnablal^2/2.
\end{equation}

\vskip-5pt
In the above condition, the left-hand side is the reduction of the exact augmented Lagrangian, while the right-hand side regards the KKT residual. We note that both \eqref{eq:cond1} and \eqref{eq:cond2} involve the penalty parameters $\veta_k$, and not every $\veta_k$ can make \eqref{eq:cond1} and \eqref{eq:cond2} satisfied simultaneously. For example, \eqref{eq:cond2} may not hold for some $\veta_k$ even with exact direction $\Delta\vz_k$. This illustrates the necessity for choosing suitable $\veta_k$. We resolve this difficulty using double While loops. 
\vskip-2pt
\underline{Outer While loop:} we check if the step acceptance condition holds, i.e., \eqref{eq:cond1} and \eqref{eq:cond2} hold simultaneously. Thus, after we break out the outer While loop, we always have a favorable direction $\tilde{\Delta}\vz_k$ and suitable penalty parameters $\veta_k$.
\vskip-2pt
\underline{Inner While loop:} with $\veta_k$ and $\delta_k\leq \delta_k^{\text{trial}}$, we repeat \eqref{eq:irs}~until \eqref{eq:cond1} is triggered. Then, we check if \eqref{eq:cond2} holds. If \eqref{eq:cond2} holds, we break the outer While loop. Otherwise, we update the parameters in an adaptive way: \vskip-20pt
\begin{equation}\label{eq:parameter_update}
\begin{aligned}
& \eta_{1,k}\gets\eta_{1,k}\nu^2,\quad \quad\quad \eta_{2,k}\gets\eta_{2,k}/\nu,\\
& \text{compute } \delta_k^{\text{trial}} \text{ as \eqref{eq:delta_trial}}, \quad\;\; \delta_k\gets(\delta_k/\nu^4\wedge\delta_k^{\text{trial}}),
\end{aligned}
\end{equation}\vskip-10pt
with $\nu>1$ being any factor larger than 1, and go back~to~reiterating \eqref{eq:irs} from the latest inexact direction with new parameters $\veta_k, \delta_k$, which lead to a new condition \eqref{eq:cond1}. The motivation of \eqref{eq:parameter_update} is to decrease $\eta_{2,k}$ and the ratio $\delta_k\eta_{1,k}/\eta_{,k2}$, but increase the product $\eta_{1,k}\eta_{2,k}$, and have $\delta_k\leq \delta_k^{\text{trial}}$ (see Lemma \ref{lemma11}). Thus, we use different powers of $\nu$ for the update. See Algorithm \ref{alg:alg1} Lines 6-14 for double While loops.

We now discuss the computational complexity of the double While loops. Once the adaptive accuracy condition~\eqref{eq:cond1} is triggered by the inner While loop, only checking the descent direction condition \eqref{eq:cond2} is left for the outer While loop.~Thus, the double While loops cost $O(\text{\# of outer loop iterations}\times(\text{cost of inner While loop}+(n+m)))$. Furthermore, the cost of inner While loop (i.e., performing \eqref{eq:irs} until \eqref{eq:cond1} is satisfied) is $O(\text{\# of inner loop iterations}\times(n+m)^2)$ with dense sketching vectors (e.g., Gaussian) and $O($\# of inner loop iterations$\times(n+m))$ with sparse sketching vectors~(e.g., Kaczmarz). As we will see in Lemmas~\ref{lemma3.4}, \ref{lemma3.5}, and \ref{lemma10}, the number of inner loop iterations is random. As for the outer loop iterations, there is no precise count when the algorithm is in the phase of adaptively selecting the parameters ($\veta_k$, $\delta_{k}$). However, when $k$ is large enough, there will be only one outer loop iteration; as shown in Lemmas \ref{lemma11} and \ref{lemma12}, all parameters ($\veta_k$, $\delta_{k}$) will be stabilized after large $k$, hence, \eqref{eq:cond2} is always satisfied as long as \eqref{eq:cond1} is satisfied.

\vspace{-2.5pt}
\subsection{Line Search and Iterate Update}\label{sec:2.3}

We select the stepsize $\alpha_k$ by doing line search and enforcing the Armijo condition on the \textit{exact augmented Lagrangian}:
\begin{equation}\label{eq:armijo}
\Lagtrial\le\Lag+\alpha_k\beta(\nablaL)^{T}\incztk,
\end{equation}
where $\vz_k=(\vx_k,\vlambda_k)$. Then, the iterate is updated as \vskip-20pt
\begin{equation}\label{eq:iterate_update}
\vz_{k+1}=\vz_k+\alpha_k\incztk.
\end{equation}\vskip-15pt
The full design of $\myalg$ is in Algorithm \ref{alg:alg1}.

\begin{algorithm}[t]
\caption{$\myalg$ Method}\label{alg:alg1}
\begin{algorithmic}[1]
\STATE{\bfseries Input:} initial \mbox{iterate}~$\vz_0$;~\mbox{sequence}~$\{\theta_k\}\subseteq~(0,1]$;~scalars $\eta_{1,0}$, $\eta_{2,0},\xi_B>0$, $\delta_0\in(0,1)$, $\beta\in(0,0.5)$, $\nu>1$;
\FOR{$k=0,1,2,\dots$}
\STATE Compute $f_k$, $\nablaxf$, $c_k$, $G_k$, $H_k$, and generate $B_k$;
\STATE Compute $\Psi_k$, $\Upsilon_k$ by~\eqref{eq:aux} and $\delta_k^{\text{trial}}$ by~\eqref{eq:delta_trial};
\STATE Set $\delta_k\gets\delta_k \wedge \delta_k^{\text{trial}}$, $\incztk\gets~\vzero$; compute $\vr_{k}$ by \eqref{eq:residual};
\WHILE{\textit{Step Acceptance Condition} does not hold}
\WHILE{$\left\Vert\vr_{k}\right\Vert>\theta_k\delta_k\normnablal/(\normA\Psi_k)$}
\STATE Generate $S\sim \mathcal{P}$ and update $\incztk$ by~\eqref{eq:irs};
\STATE Compute $\vr_{k}$ by~\eqref{eq:residual};
\ENDWHILE
\IF {$(\nablaL)^T\incztk>-\eta_{2,k}\normnablal^2/2$}
\STATE Update $\eta_{1,k}, \eta_{2,k}, \delta_k^{\text{trial}}, \delta_k$ as \eqref{eq:parameter_update};
\ENDIF
\ENDWHILE
\STATE Select $\alpha_k$ to satisfy~\eqref{eq:armijo}; update the iterate by~\eqref{eq:iterate_update};
\STATE	Set $\eta_{1,k+1}\gets \eta_{1,k}$, $\eta_{2,k+1}\gets \eta_{2,k}$, and $\delta_{k+1}\gets\delta_k$;
\ENDFOR
\end{algorithmic}
\end{algorithm}
\setlength{\textfloatsep}{0.15cm}

\section{Well-posedness and Global Convergence}\label{sec:3}

In this section, we first study the well-posedness of Algorithm \ref{alg:alg1} by showing that the step acceptance condition in~Section \ref{sec:2.2} is always triggered with a finite $j$; thus, the double While loops in Algorithm \ref{alg:alg1} Lines 6-14 terminate~in~finite time. We then show a global almost sure convergence guarantee --- starting from any initial point, the KKT residual $\|\nabla\mL_k\|$ converges to zero \textit{almost surely}. Here, the randomness plays a key role in the analysis, because the inexact direction is calculated by a randomized solver and all the algorithmic components that are affected by the direction are also random. For example, the step acceptance condition and Armijo condition are governed by random sketching matrices. We begin by stating the assumptions.

\begin{assumption}\label{assum1}
The iterates $\{\vx_k, \vlambda_k\}_{k\ge0}$ are contained~in a convex compact set $\mathcal{X}\times\Lambda$ such that over $\mX$, the objective $f$ and constraint $c$ are twice continuously differentiable with Hessians being Lipschitz continuous.
\end{assumption}


\begin{assumption}\label{assum2}
	There exist absolute constants $\xi_{G}, \xi_{B}, \Upsilon_B$ such that (i) the Jacobian $G_k$ has full row rank with $G_kG_k^T \\ \succeq \xi_{G} I$, $\forall k$; (ii) the modified Hessian $B_k$ satisfies $\vu^TB_k\vu \\ \ge\xi_B \Vert \vu\Vert^2$ for any $\vu\in\{\vu:G_k\vu=0\}$ and $\normB\le\Upsilon_B$.
\end{assumption}

\begin{assumption}\label{assum3}
There exists a constant $\pi\in(0,1]$ such that the sketching matrices $S_{k,j}\sim S$, $iid$, satisfy $P(S^T\vu\neq \vzero)\geq \pi$ for any $\vu\in\mR^{n+m}\backslash\{\vzero\}$.
\end{assumption}

All three assumptions are mild, standard, and commonly imposed in the literature. In particular, Assumptions \ref{assum1} and \ref{assum2} are required for the analysis of exact, inexact, stochastic SQP methods \citep{Bertsekas1982Constrained, Boggs1995Sequential, Nocedal2006Numerical, Byrd2008Inexact, Byrd2010inexact, Na2023Inequality, Na2022adaptive}. An alternative statement of Assumption \ref{assum1} is to assume the iterates lie in a convex open set, and the objective $f$, the constraint $c$, together with their gradients and Hessians, are Lipschitz continuous and bounded over that set \citep{Curtis2021Inexact}. Assumption~\ref{assum2} is required to ensure that the Newton system \eqref{eq:kkt_sys_reduced} has a unique solution. Assumption \ref{assum3} is a condition on~the~sketching~distribution.~It~holds~for~\mbox{various}~choices~of~sketching~\mbox{matrices}. For example, in randomized Kaczmarz method, $S=\ve^i\in\mR^{n+m}$, the $i$-th canonical basis, with equal probability. This choice satisfies Assumption \ref{assum3} with $\pi = 1/(n+m)$. It is also~immediate to see that Assumption \ref{assum3} holds with $\pi = 1$ for any continuous sketching distribution (e.g., Gaussian sketching). We do not impose any conditions on the sketching dimension $d$; thus, we can set $d=1$ (i.e., pseudoinverse in \eqref{eq:irs} becomes reciprocal) for sake of low flops complexity.

We first present Lemmas \ref{lemma3.4} and \ref{lemma3.5}, which are the building blocks for showing the well-posedness of Algorithm \ref{alg:alg1}.

\begin{lemma}\label{lemma3.4}
For any outer iteration $k$ and inner iteration $j$, let $Q_{k,j}\in\mR^{(n+m)\times d}$ be a matrix that has orthonormal columns spanning the space $row(S_{k,j}^T\KKT)$; and let $\{j_l^{k}\}_{l\ge 0}$ be a subsequence of the inner iteration $j$, where $j_0^k=0$ and $j_l^{k}$, $l\geq 1$, is recursively defined to satisfy \vskip-15pt
\begin{equation*}
col(Q_{k,j_{l-1}^k})\cup\cdots\cup col(Q_{k,j_{l}^k-1})=\mR^{n+m}.
\end{equation*}\vskip-10pt
Let $L$ be any positive integer. Under Assumptions \ref{assum2}, \ref{assum3}, and for any $k$, let us suppose Algorithm \ref{alg:alg1} reaches $\citer$. Then, the event \vskip-15pt
\begin{equation}\label{eq:subseq}
\mathcal{A}_{k}=\cap_{l=1}^L\{j_l^k<\infty\}
\end{equation} \vskip-5pt
happens with probability one.
\end{lemma}

By Lemma \ref{lemma3.4}, we know that the inner iteration $j$ has a~subsequence $j_l^k$ such that the union of the space $row(S_{k,j}^T\KKT)$ from $j = j_l^k$ to $j = j_{l+1}^k-1$, $\forall l\geq 0$, is the full space~$\mR^{n+m}$. Such a full space expansion property is critical to show that the random approximation error, although does not monotonically decrease, has a decreasing subsequence.

\begin{lemma}[A subsequence of error decays linearly]\label{lemma3.5}
Under Assumptions \ref{assum2}, \ref{assum3}, let us suppose the event $\mathcal{A}_{k}$ in \eqref{eq:subseq} happens. Then, there exists a sequence of scalars $\{\gamma_{k,l}\}_{l=1}^L$ such that $\gamma_{k,l}\sim\gamma_k$ is an iid realization of a random variable $\gamma_k\in[0,1)$, and we have for $1\leq l\le L$, \vskip-12pt
\begin{equation*}
\Vert\edzkjlk\Vert\le\gamma_{k,l}\Vert\edzkjlbk\Vert.
\end{equation*}
\end{lemma}

\vskip-5pt
Lemma \ref{lemma3.5} suggests that a subsequence of the approximation error, $\{\|\edzkjlk\|\}_{l\geq 0}$, decays linearly with the rate being an iid copy of some random variable $\gamma_k\in[0,1)$. We should mention that the statement of Lemma \ref{lemma3.5} is deterministic (i.e., $\{\gamma_{k,l}\}_l$ are realized) since we suppose the event $\mA_k$ happens in the statement.

Lemma \ref{lemma3.5} directly leads to the result that the adaptive~accuracy condition \eqref{eq:cond1} can be satisfied with finite inner iterations, although the iteration number may be random.

\begin{lemma}[Well-posedness of accuracy condition]\label{lemma10}

Under Assumptions \ref{assum1}, \ref{assum2}, \ref{assum3}, and for any outer iteration $k$ and any $\delta_k,\theta_k>0$, let us suppose Algorithm \ref{alg:alg1} reaches $\citer$. Then, with probability one, there exists a finite number $J_k<\infty$ such that the accuracy condition \eqref{eq:cond1} with $\delta_k,\theta_k>0$ can~be satisfied by iterating \eqref{eq:irs} for $J_k$ times.

\end{lemma}

From the proof of Lemma \ref{lemma10}, we know the right-hand side of the accuracy condition \eqref{eq:cond1} can be replaced by any~positive upper bound, and the condition is still satisfied with~(random) finite inner iterations. Lemma \ref{lemma10} suggests that the~inner while loop (cf. Algorithm \ref{alg:alg1}, Lines 7-10)~\mbox{always}~terminates properly. We next investigate the descent direction condition \eqref{eq:cond2} to complete the well-posedness study.

\begin{lemma}[Well-posedness of descent direction condition]\label{lemma11}

Under Assumptions \ref{assum1}, \ref{assum2}, \ref{assum3}, and for any outer iteration $k$, we let $\incztkj$ be the inexact solution to \eqref{eq:kkt_sys_reduced} that satisfies \eqref{eq:cond1}. Then, there exists a constant $\Upsilon = \Upsilon(\xi_{G},\xi_{B},\Upsilon_B)>0$ large enough such that the descent direction condition \eqref{eq:cond2} is satisfied as long as
\begin{equation}\label{eq:thres}
\eta_{1,k}\eta_{2,k}\geq \Upsilon\quad \text{ and }\quad \eta_{2,k} \vee \delta_k\eta_{1,k}/\eta_{2,k}\leq 1/\Upsilon.
\end{equation}
	
\end{lemma}

Lemma \ref{lemma11} suggests that the descent direction condition \eqref{eq:cond2} holds as long as $\eta_{1,k}$ is large enough and $\eta_{2,k}, \delta_k$ are small enough. By our updating rule \eqref{eq:parameter_update}, we increase the quantity $\eta_{1,k}\eta_{2,k}$ and decrease the quantities $\eta_{2,k}$ and $\delta_k\eta_{1,k}/\eta_{2,k}$ by a factor of $\nu>1$ whenever \eqref{eq:cond2} is not satisfied. Thus, \eqref{eq:thres} (and hence \eqref{eq:cond2}) will be finally satisfied. Combining~Lemmas \ref{lemma10} and \ref{lemma11}, we have now shown the double While~loops (cf. Algorithm \ref{alg:alg1}, Lines 6-14) terminate properly.

We next study the behavior of adaptive penalty parameters.

\begin{lemma}[Stability of adaptive parameters]\label{lemma12}
Under Assumptions \ref{assum1}, \ref{assum2}, \ref{assum3}, with probability one, there exists an iteration threshold $K$ such that the parameters $(\eta_{1,k},\eta_{2,k}, \\ \delta_{k})$ are stabilized after $K$ iterations, that is, $(\eta_{1,k},\eta_{2,k}, \delta_{k})\\ = (\eta_{1,K},\eta_{2,K},\delta_{K})$, $\forall k\geq K$. 
\end{lemma}

The stability of penalty parameters is crucial for global~convergence. Due to our adaptive design, the algorithm chooses suitable parameters automatically. Thus, the augmented~Lagrangian merit function for line search may differ from step by step. Then, we cannot accumulate the decreases across the steps (since each step may decrease a different function). Lemma \ref{lemma12} suggests that our adaptive design leads to a stabilized augmented Lagrangian in the end; thus we can accumulate all the decreases on the tail.

We then show the stepsize $\alpha_k$ has a uniform lower bound.

\begin{lemma}[Armijo condition]\label{lemma13}
Under Assumptions \ref{assum1}, \ref{assum2}, \ref{assum3}, with probability one, there exists $\alpha_{\min}>0$ such that $\alpha_k\ge\alpha_{\min}$, $\forall k\geq 0$.
\end{lemma}

With all of the above lemmas, we now establish the global almost sure convergence of Algorithm \ref{alg:alg1} in Theorem \ref{theorem1}.

\begin{theorem}[Global convergence]\label{theorem1}
Under Assumptions \mbox{\ref{assum1}--\ref{assum3}},~with~\mbox{probability}~one,~$\normnablal \rightarrow~0$~as~\mbox{$k\rightarrow \infty$.}
\end{theorem}

Compared to randomized inexact Newton methods for unconstrained convex optimization \citep{Qu2016SDNA, Doikov2018Randomized, Luo2016Efficient, Gower2019Rsn, Hanzely2020Stochastic}, we establish the convergence of KKT residual that interprets the constraint violation. Compared to deterministic inexact Newton methods for constrained optimization \citep{Byrd2008Inexact, Byrd2010inexact}, our global result~holds almost surely, instead of deterministically; and our next~local analysis also complements the missing part in their studies.

\section{Local Convergence}\label{sec:4}

In this section, we establish the local convergence rate of~Algorithm \ref{alg:alg1}. We first show that when we set $\theta_k = \theta\in(0,1]$, Algorithm \ref{alg:alg1} exhibits local linear convergence. We~then~show that when we let $\theta_k$ decay to zero, i.e., when we gradually sharpen the adaptive accuracy condition \eqref{eq:cond1}, Algorithm \ref{alg:alg1} exhibits local superlinear convergence.

We first present two additional assumptions that are necessary for local analysis.

\begin{assumption}\label{assum4}
We assume $f$ and each coordinate of $c$~are thrice continuously differentiable over $\mathcal{X}$.
\end{assumption}

\begin{assumption}[Hessian modification vanishes]\label{assum5}
We assume~$\Vert H_k-~B_k\Vert=~O(\tau_k)$~for~a~\mbox{sequence}~$\tau_k\rightarrow~0$~as~$k\rightarrow~\infty$.
\end{assumption}

Assumption \ref{assum4} strengthens Assumption \ref{assum1} by requiring one more derivative --- the third derivative --- of $f$ and $c$ to exist. This condition is standard for local \mbox{analysis}~when using the exact augmented Lagrangian merit function in~the algorithm \citep{Bertsekas1982Constrained, Zavala2014Scalable, Na2021fast, Na2021Global}, because the Hessian of the augmented Lagrangian $\nabla_{\vx}^2\mathcal{L}_{\veta}$ requires the existence of $\nabla^3f$ and $\nabla^3 c$. Fortunately, the third derivatives are never computed in the algorithm. Assumption \ref{assum5} assumes the Hessian modification gradually vanishes, which is also standard in the SQP literature \citep{Boggs1995Sequential,Nocedal2006Numerical}. It is worth mentioning that Assumption \ref{assum5} implies that $H_k$ satisfies Assumption~\ref{assum2}-(ii) in the limit, which is known as the second-order sufficient conditions. There are multiple ways to generate $B_k$ that satisfies Assumptions~\ref{assum2} and \ref{assum5}. One example is to test the positiveness of $Z_k^TH_kZ_k$, where the columns of $Z_k\in\mR^{n\times(n-m)}$ span the null space of Jacobian $G_k$. If $Z_k^TH_kZ_k$ is positive definite, then we set $B_k=H_k$; otherwise we set $B_k=H_k+(\xi_B+\Vert H_k\Vert)I$. By this way,~we~have $\tau_k = 0$ for all large enough $k$.

\begin{theorem}[Local linear convergence]\label{theorem2}
Let $\vz^\star$ be a local solution to \eqref{eq:opt_main} and $\theta_k=\theta\in(0,1]$, $\forall k$. Under Assumptions \ref{assum1}--\ref{assum3}, \ref{assum4}, \ref{assum5} and suppose $\citer\rightarrow\vz^\star$, for~all sufficiently large $k$, we have $\alpha_k=1$ and (noting that $\theta\delta_K<1$) \vskip-0.7cm
\begin{equation*}
\Vert\niter-\vz^\star\Vert\le (1+\varphi)\theta\delta_K\Vert\citer-\vz^\star\Vert, \quad \text{for any $\varphi>0$.}
\end{equation*}
\end{theorem}

\begin{corollary}[Local superlinear convergence]\label{collorary1}
Let $\vz^\star$ be a local solution to \eqref{eq:opt_main} and $\theta_k$ be any input sequence such that $\theta_k\rightarrow 0$ as $k\rightarrow\infty$. Under~Assumptions \ref{assum1}--\ref{assum3}, \ref{assum4}, \ref{assum5} and suppose $\citer\rightarrow\vz^\star$, for all sufficiently large $k$, we have $\alpha_k=1$ and that\vskip-22pt
\begin{equation*}
\left\Vert\niter-\vz^\star\right\Vert\le O(\theta_k+\tau_k)\left\Vert\citer-\vz^\star\right\Vert+O\left(\|\citer-\vz^\star\|^2\right).
\end{equation*}
\end{corollary}
\vskip-10pt
Since $\theta_k$ is a factor of accuracy condition on the sketching solver (cf. \eqref{eq:cond1}), and $\delta_k$ is stabilized to $\delta_K$, we know~a~decaying input sequence $\theta_k$ suggests Algorithm \ref{alg:alg1} performs more inner iterations in expense of a faster local rate.

From the global convergence in Theorem~\ref{theorem1}, we know that Algorithm~\ref{alg:alg1} generates iterates that converge to any stationary points. In contrast, as we mentioned earlier, we assume that the second-order sufficient conditions hold at $\vz^\star$ for the local convergence results. Thus, Theorem~\ref{theorem2}~and~Corollary \ref{collorary1} indicate that the iterates generated by Algorithm~\ref{alg:alg1} will exhibit linear/superlinear local rates, provided the stationary point is a second-order stationary point.

\section{Experiments}\label{sec:5}

We benchmark $\myalg$ (Algorithm \ref{alg:alg1}) on nonlinear problems in CUTEst collection set \citep{Gould2014CUTEst}, on constrained logistic regression with data from LIBSVM \citep{Chang2011LIBSVM}, and on a PDE-constrained problem \citep{Hintermueller2002Primal}. We compare the performance of
Algorithm \ref{alg:alg1} with that of two inexact SQP methods designed for constrained problems with deterministic solvers: Algorithm B of \citet{Byrd2008Inexact} and its adaptive modification. The two methods are detailed in Algorithms \ref{alg:alg2} and \ref{alg:alg3} in Appendix \ref{sec:appenB}. Compared to Algorithm \ref{alg:alg2} that uses a fixed 
bound throughout all iterations, Algorithm \ref{alg:alg3} adaptively controls the accuracy of a deterministic solver and can be seen as a deterministic version of Algorithm \ref{alg:alg1}. However, both methods employ an nonsmooth $\ell_1$ penalized merit function that differs from \eqref{eq:opt_merit}. As another baseline for the comparison, we also consider standard augmented Lagrangian method with an inexact Newton subproblem solver (see Algorithm \ref{alg:alg4} in Appendix \ref{sec:appenB}). For Algorithm \ref{alg:alg1},~we apply two sketching
distributions: (1) Gaussian vector sketch and (2) Randomized Kaczmarz sketch, referred as $\myalg$-GV and $\myalg$-RK, respectively. As suggested in \citet{Byrd2008Inexact} and \citet{Nocedal2006Numerical}, we use GMRES \citep{Saad1986GMRES} as~the deterministic solver for Algorithms \ref{alg:alg2}, \ref{alg:alg3}, and \ref{alg:alg4}. We evaluate each algorithm with the following three criteria: (1) the~KKT residual $\normnablal$; (2) the number of objective and constraints evaluations; and (3) the number of gradient and Jacobian evaluations. Further, we assess each algorithm with the performance profile \citep{Dolan2002Benchmarking} for CUTEst and constrained logistic regression, in which the total number of flops is used as a performance measure. For all methods, we stop iterating~if \vskip-15pt
\begin{equation*}
\normnablal \le 10^{-4}\quad\text{OR}\quad k\ge 10^4.    
\end{equation*}\vskip-10pt
The parameters of each algorithm are specified as follows. (We further test the sensitivity to parameters for Algorithm~\ref{alg:alg1} in Section \ref{subsec:robust}).

\vskip-5pt
\underline{Alg. \ref{alg:alg1}}: {\small $\eta_{2,0}=\delta_0=\xi_B=\beta=0.1$, $\eta_{1,0}= \theta_k= 1$,~$\nu=1.5$.}
\vskip-5pt
\underline{Alg. \ref{alg:alg2}}: we follow the exact same setup as \citet{Byrd2008Inexact}. In particular, with their notation, $\eta=10^{-8}$, $\kappa_1=\epsilon=\tau=\xi_B=0.1$, $\pi_{0} = \kappa=1$, $\beta = 1\vee\| \nabla\mathcal{L}_0\|_1/(\| c_0\|_1 + 1)$.
\vskip-5pt
\underline{Alg. \ref{alg:alg3}}: we follow the same setup as Algorithm \ref{alg:alg2} and let $\eta= 10^{-8}$, $\kappa_{0}=\xi_B=0.1$, $\pi_{0}=1$, $\nu=1.5$.
\vskip-5pt
\underline{Alg. \ref{alg:alg4}}: {\small$\kappa=10^{-4}$, $\tau_0=\eta=0.1$, $\mu_0=1$, $\nu_{\mu}=1.5$, $\nu_{\tau}=0.5$.}

Our code for the implementation is available at \url{https://github.com/IlgeeHong/AdaSketch-Newton}.

\begin{figure*}[t]
	\begin{center}
		\begin{minipage}[ht]{.33\linewidth}
			\centerline{\includegraphics[scale=0.23]{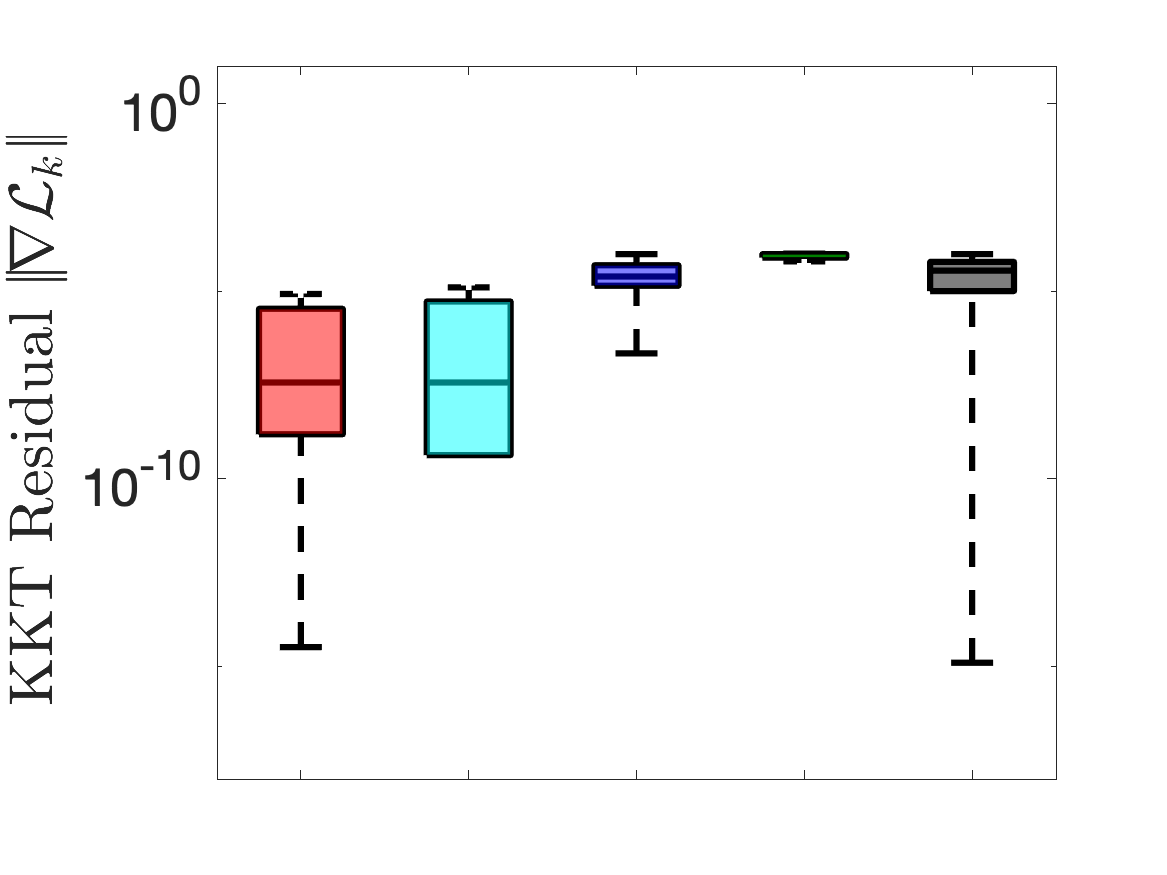}}
		\end{minipage}
		\begin{minipage}[ht]{.33\linewidth}
			\centerline{\includegraphics[scale=0.23]{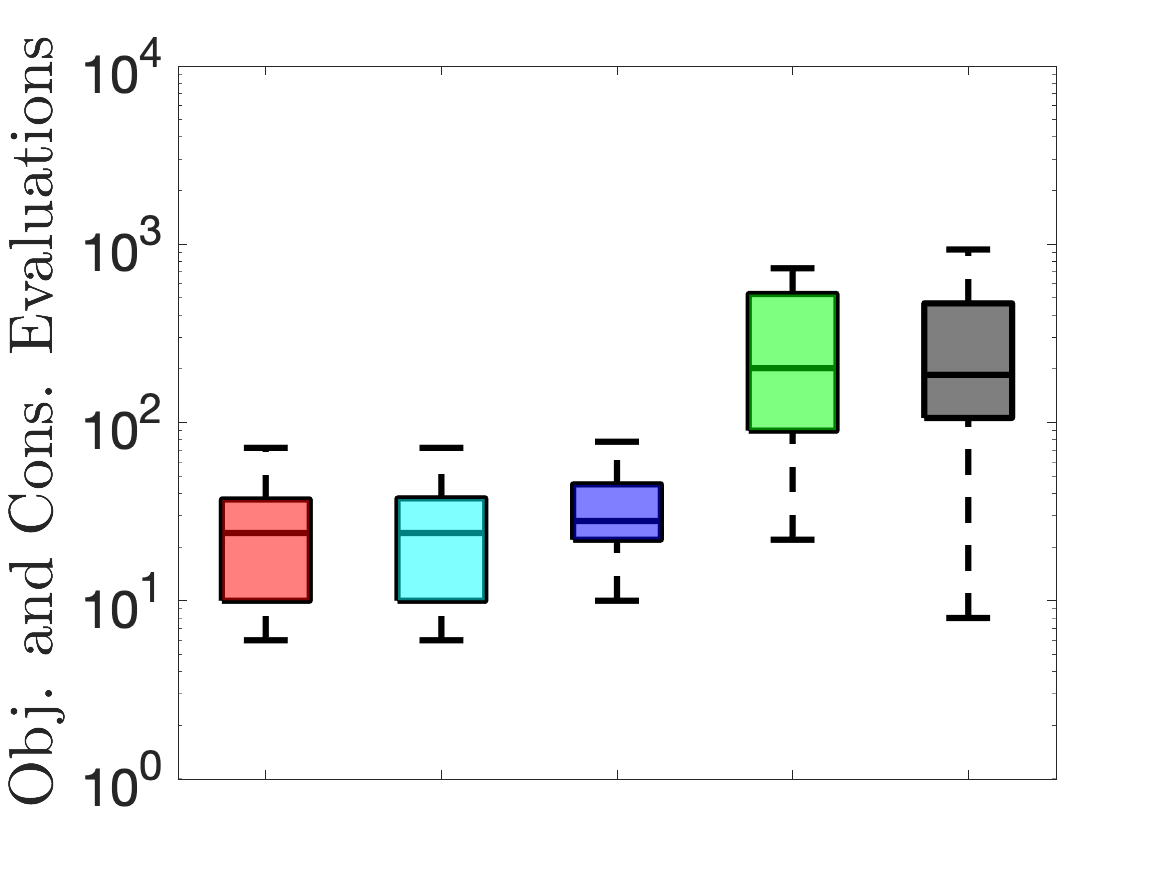}}
		\end{minipage}
		\begin{minipage}[ht]{.33\linewidth}
			\centerline{\includegraphics[scale=0.23]{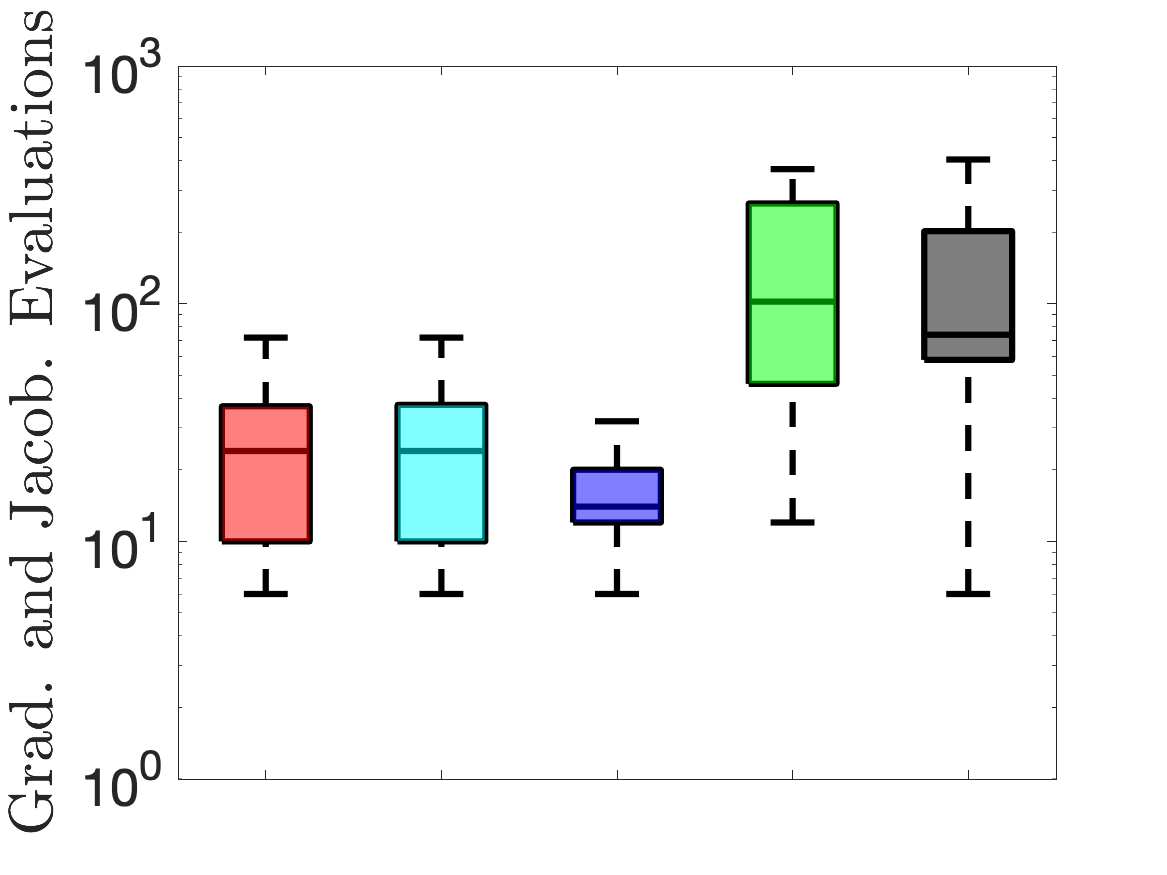}}
		\end{minipage}\\
		\begin{minipage}[ht]{.9\linewidth}
			\centerline{\includegraphics[scale=0.60]{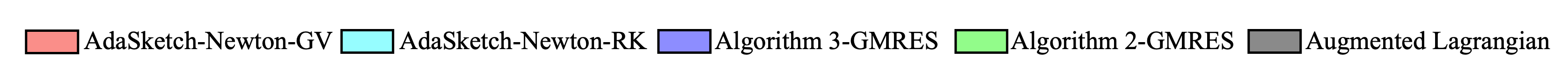}}
		\end{minipage}
		\caption{The boxplots of the KKT residual, the number of objective and constraints evaluations, and the number of gradient and Jacobian evaluations for $\myalg$, Algorithm~\ref{alg:alg2}, Algorithm~\ref{alg:alg3}, and Algorithm \ref{alg:alg4} on CUTEst problems.}\label{fig:fig1}
	\end{center}
\end{figure*}

\begin{figure*}[h]
	\begin{center}
		\begin{minipage}[ht]{.24\linewidth}
			\centerline{\includegraphics[scale=0.21]{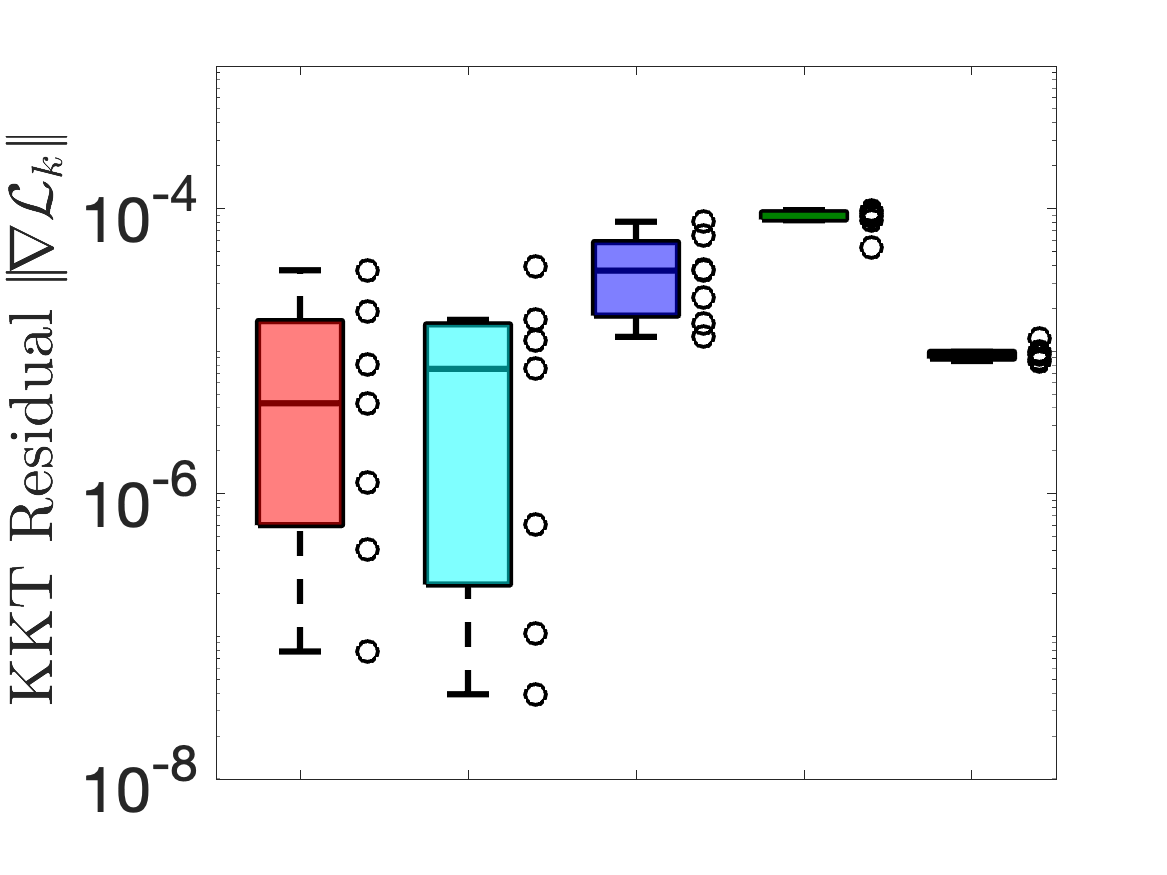}}
		\end{minipage}%
		\begin{minipage}[ht]{.24\linewidth}
			\centerline{\includegraphics[scale=0.21]{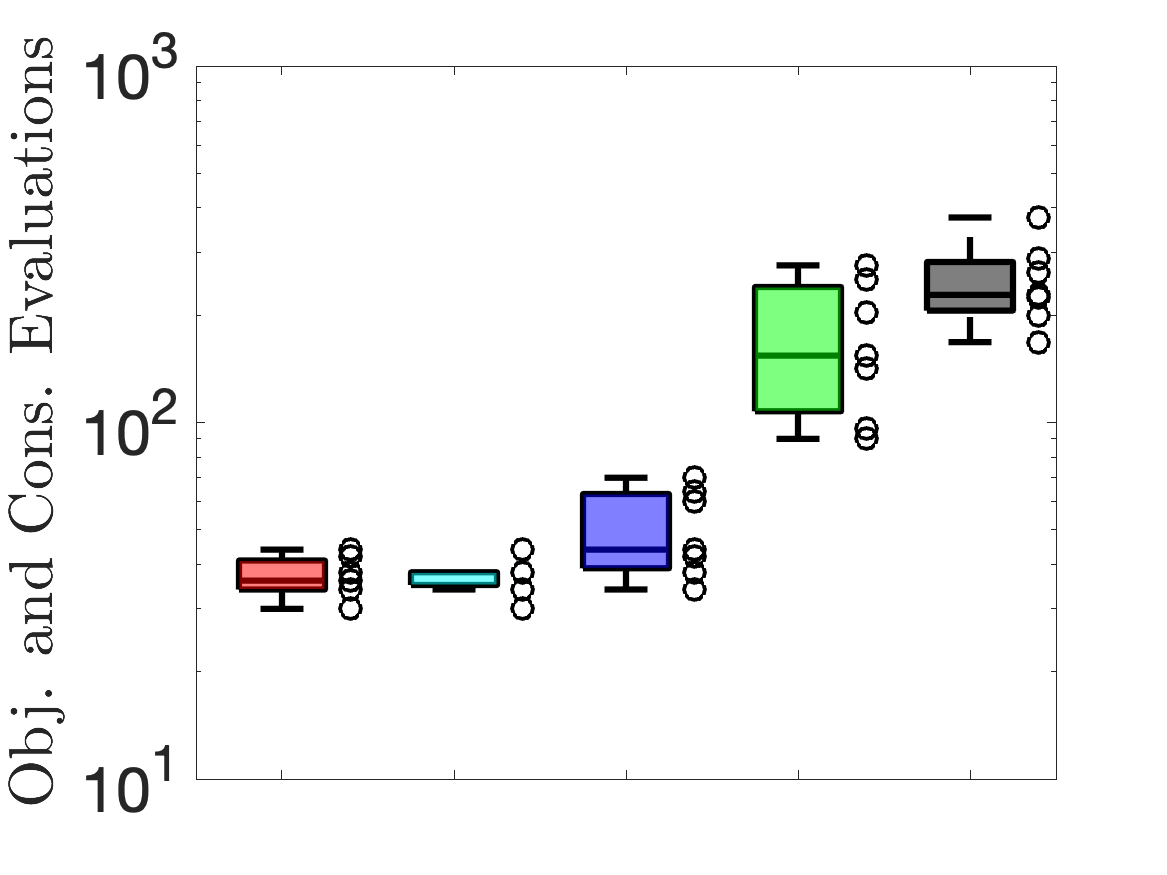}}
		\end{minipage}
		\begin{minipage}[ht]{.24\linewidth}
			\centerline{\includegraphics[scale=0.21]{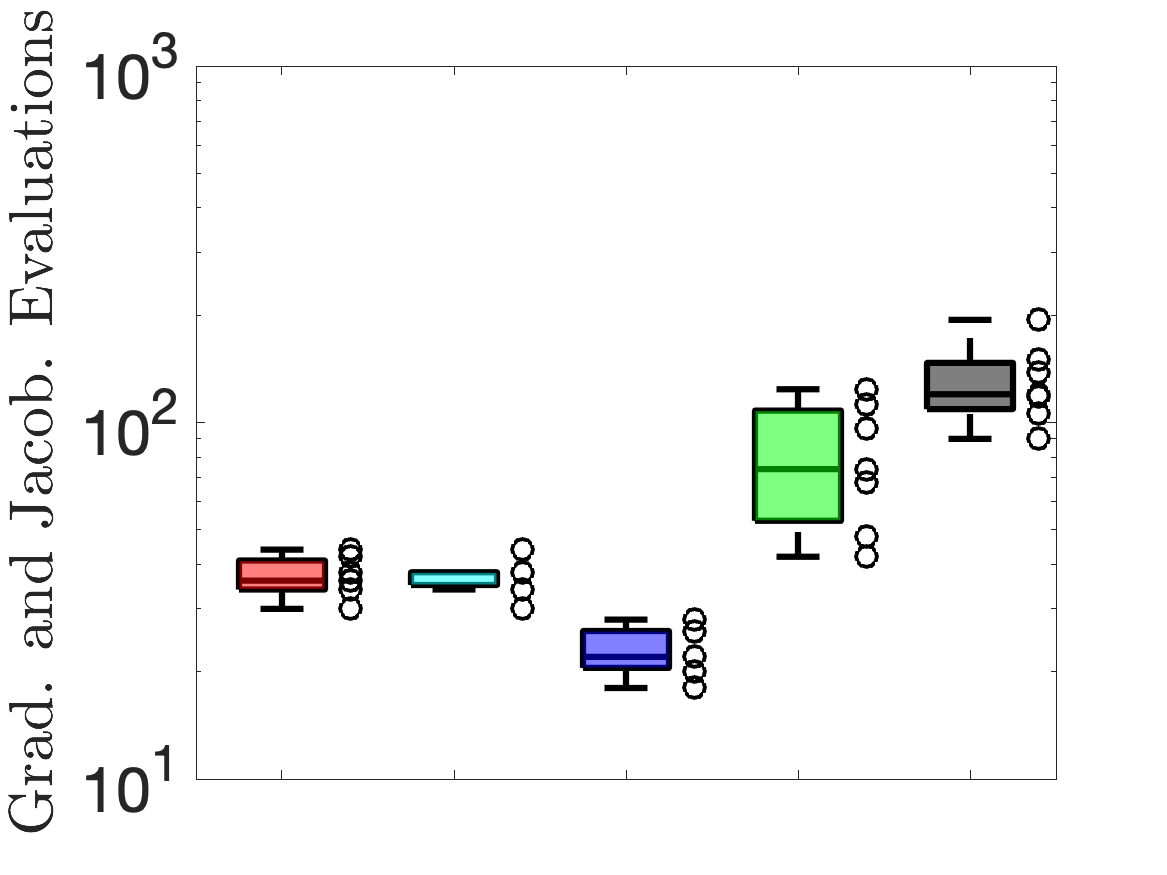}}
		\end{minipage}
		\begin{minipage}[ht]{.24\linewidth}
			\centerline{\includegraphics[scale=0.21]{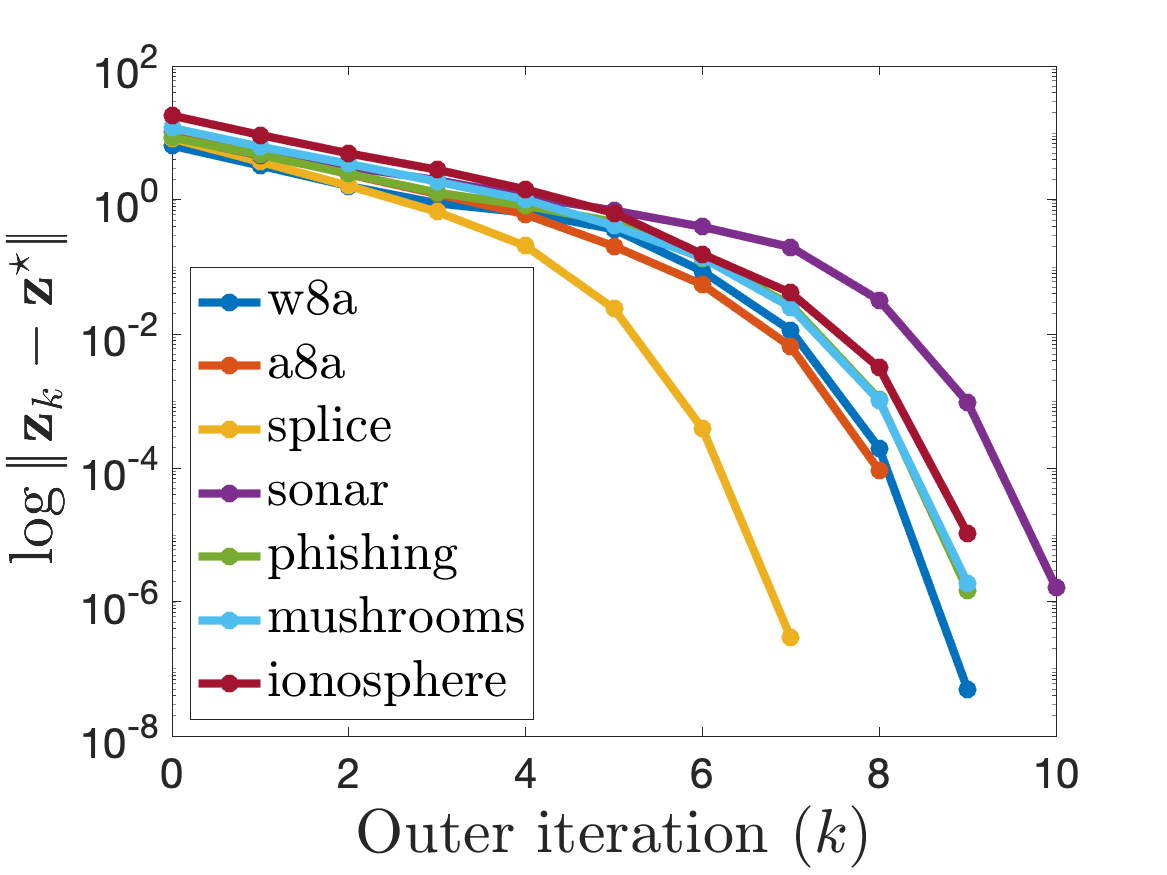}}
		\end{minipage}\\
		\begin{minipage}[ht]{.24\linewidth}
			\centerline{\includegraphics[scale=0.60]{Figure/Legend.png}}
		\end{minipage}
		\caption{The  boxplots of the KKT residual, the number of objective and constraints evaluations, and the number of gradient and~Jacobian evaluations for $\myalg$, Algorithms \ref{alg:alg2}, Algorithm \ref{alg:alg3}, and Algorithm \ref{alg:alg4} on 7 LIBSVM datasets. The right figure shows the decay trajectories of the log error $\log\|\vz_k - \vz^\star\|$, where each line corresponds to each LIBSVM dataset ($\vz^\star$ is estimated by the last iterate of the IPOPT solver \citep{Waechter2006implementation}).}\label{fig:fig2}
	\end{center}
\end{figure*}

\subsection{CUTEst}\label{subsec:cutest}

Among all the problems collected in CUTEst test set, we~select the problems that have a non-constant $f$ with $n < 1000$, and contain only equality constraints with Jacobian $G_k$ being full rank at each step for all three algorithms. This leads to a total of 47 problems. For each problem, the initial iterate $\vz_0=(\vx_0,\vlambda_0)$ is provided by CUTEst package, and we average results over 10 independent runs for Algorithm \ref{alg:alg1}.

The boxplots over 47 problems of the three criteria for the five methods are shown in Figure \ref{fig:fig1}. From Figure \ref{fig:fig1}, we~observe that AdaSketch-Newton (both GV and RK) outperforms Algorithms \ref{alg:alg2} and \ref{alg:alg3} in terms of the KKT residual and the number of objective and constraints evaluations.~This~observation is expected since $\myalg$ employs a smooth exact augmented Lagrangian as the merit function, with an adaptive step acceptance condition to ensure a fast local convergence. Our design leads to steeper decreases on the merit function and fewer outer iterations. Compared to Algorithms \ref{alg:alg2} and \ref{alg:alg3}, $\myalg$ requires gradient and Jacobian evaluations to evaluate the merit function (cf. \eqref{eq:opt_merit}), which are not needed for evaluating the $\ell_1$ penalized merit function used in the other two methods. Despite this fact, we observe in Figure \ref{fig:fig1} that AdaSketch-Newton enjoys a competitive performance on the number of gradient and Jacobian evaluations with Algorithm \ref{alg:alg3} and even performs better than Algorithm \ref{alg:alg2}. In addition, we also note that our adaptive modification of Algorithm \ref{alg:alg2}~(i.e.,~Algorithm \ref{alg:alg3}) is superior to Algorithm \ref{alg:alg2} in all three criteria, which illustrates the effectiveness of our adaptive technique. Finally, we see in Figure~\ref{fig:fig1} that the augmented Lagrangian method (Algorithm~\ref{alg:alg4}) underperforms AdaSketch-Newton (both GV and RK) in terms of all three criteria.

\vskip-3pt

The performance profiles over 47 problems of the total number of flops for the five methods are presented in the left-hand side of Figure~\ref{fig:fig3} in Appendix~\ref{sec:appenB}. From the figure, we see that the two AdaSketch-Newton methods are the fastest, followed by Algorithm \ref{alg:alg3}. Algorithm \ref{alg:alg2} and the augmented Lagrangian method are much slower than all three adaptive methods. These results show that AdaSketch-Newton (both GV and RK) is the most effective among Algorithms \ref{alg:alg2}, \ref{alg:alg3}, and \ref{alg:alg4} in terms of the total number of flops.

\subsection{Constrained Logistic Regression}\label{subsec:libsvm}

We consider equality-constrained logistic regression problems of the form\vskip-25pt
\begin{align*}
\underset{\vx\in\mR^n}{\min} \;\; & f(\vx)=\frac{1}{N}\sum_{i=1}^N \log\left(1+\exp\left({-y_{i}\cdot\vd_i^T\vx}\right)\right)\\
\text{s.t.}\;\; & A\vx=\vb,\;\Vert\vx\Vert^2=1,
\end{align*}\vskip-8pt
where $(y_i, \vd_i)\in\{-1,1\}\times\mR^n$ are the $i$-th data point and $A\in\mR^{m\times n}$ and $\vb\in\mR^{m}$ are linear constraint coefficients. We consider 7~binary classification datasets from LIBSVM for which $12\le n\le 10^3$~and~$256\le N\le 10^5$. Details of the datasets are given in Table~\ref{tab:table1}. For linear constraint, we set $m=10$ and randomly generate each entry of $A$ and $\vb$ from a standard normal distribution. Combining with unit norm constraint, we have 11 equality~constraints in total. For all datasets, we set $\vz_0$ to be all-one vector and, again, average results over 10 independent runs for Algorithm \ref{alg:alg1}.

The boxplots (with points) over 7 datasets for the five~methods and the decay trajectory of the log error ($\log\|\vz_k-\vz^\star\|$) for $\myalg$-RK are shown in Figure \ref{fig:fig2}. From the figure, we have similar observations to CUTEst. AdaSketch-Newton (both GV and RK) outperforms Algorithm \ref{alg:alg3} in terms of the KKT residual and the number of objective and constraints evaluations, and it performs better than Algorithm \ref{alg:alg3} and the augmented Lagrangian method in terms of all three criteria. Algorithm \ref{alg:alg3} is superior to Algorithm \ref{alg:alg2} in all three criteria. The trajectory plots of all 7 datasets show that the log error decays at least linearly, that is consistent with our theoretical analysis in Section \ref{sec:4}.

The performance profiles over 7 datasets of the total number of flops for the five methods are presented in the right-hand side of Figure~\ref{fig:fig3} in Appendix~\ref{sec:appenB}. From the figure, we see~all three adaptive methods (two $\myalg$, Algorithm~\ref{alg:alg3}) are more effective than Algorithm \ref{alg:alg2} and the augmented Lagrangian in terms of the total number of flops.

\begin{table}[H]
\caption{Dataset Statistics}\label{tab:table1}
\begin{center}
\resizebox{\columnwidth}{!}{
\begin{small}
\begin{sc}
\begin{tabular}{lcccr}
\toprule
Dataset & feature dimension ($n$) & \# data points ($N$)\\ 
\midrule
w8a & 300 & 49,749 \\
a9a & 123  & 32,561  \\
splice & 60 & 1,000  \\ 
sonar & 60 & 208  \\ 
phishing & 68 & 11,055  \\ 
mushrooms & 112 & 8,124 \\ 
ionosphere & 34  & 351  \\ 
\bottomrule
\end{tabular}
\end{sc}
\end{small}}
\end{center}
\vskip -0.2in
\end{table}

\subsection{PDE-constrained Problem}\label{subsec:pde}

We consider optimal control problem with Dirichlet~boundary conditions of the form
\begin{align*}
\underset{x,y}{\min} \;\; & \frac{1}{2}\Vert x-u\Vert_{L^2(\Omega)}^2+\frac{\zeta}{2}\Vert y\Vert_{L^2(\Omega)}^2\\
\text{s.t.}\;\; & -\Delta x = y\text{ in } \Omega,\;x=\vzero\text{ on }\partial\Omega,
\end{align*}\vskip-8pt
where $u\in L^2(\Omega)$ is a reference function and $\zeta>0$ is~a regularization parameter. We discretize $\Omega$ on an evenly spaced $N\times N$ grid. For reference function $u=[u_{ij}]^{N,N}_{i,j=1}$, we follow \citet{Curtis2021Inexact} and choose for all $(i,j)\in\{1,\dots,N\}\times\{1,\dots,N\}$ the following:
\begin{align*}
\resizebox{0.99\columnwidth}{!}{$u_{ij} = \sin\left(4+\frac{\epsilon_N}{\epsilon_S}\left(i-\frac{N+1}{2}\right)\right)+\cos\left(3+\frac{\epsilon_N}{\epsilon_S}\left(j-\frac{N+1}{2}\right)\right)$},
\end{align*}
where $\epsilon_N,\epsilon_S>0$. We select $N=3$ and $\zeta=10^{-1}$, and as suggested in \citet{Curtis2021Inexact}, we choose $\epsilon_N=10^{-1}$ and  $\epsilon_S=\sqrt{15}$. We set $\vz_0=(x_0,y_0,\lambda_0)$ to be all-one matrix and average over 10 independent runs for Algorithm~\ref{alg:alg1}. 

The numerical results for the five methods are summarized in Table~\ref{tab:table2} in Appendix~\ref{sec:appenB}. From Table~\ref{tab:table2}, we reconfirm our observations from CUTEst and constrained logistic regression; AdaSketch-Newton (both GV and RK) outperforms Algorithms \ref{alg:alg2}, \ref{alg:alg3}, and \ref{alg:alg4} in terms of all three criteria.

\subsection{Sensitivity to Parameters}\label{subsec:robust}

We test the sensitivity of $\myalg$-GV to four input parameters, $(\eta_{1,0}, \eta_{2,0}, \delta_0, \beta)$, on 47 CUTEst problems used in Section \ref{subsec:cutest}. Here, $(\eta_{1,0}, \eta_{2,0})$ are the initial~penalty parameters of the augmented Lagrangian; $\delta_0$ is the initial threshold parameter that controls the accuracy of the randomized iterative sketching solver; and $\beta$ is a parameter~in the Armijo condition in line search. As set in Section \ref{subsec:cutest}, the default values are $(\eta_{1,0},\eta_{2,0},\delta_0,\beta)=(1,0.1,0.1,0.1)$. Here, we choose a larger $\eta_{1,0}$ and a smaller $\eta_{2,0}$ to make \eqref{eq:opt_merit} close to a standard augmented~Lagrangian. We vary the four parameters in ranges as follows: $\eta_{1,0}\in \{0.1,1,10\}$,~$\eta_{2,0}\in \{0.01,0.1,1\}$, $\delta_0\in \{0.01, 0.1,0.9\}$, {\small$\beta\in\{10^{-7},10^{-3}, 0.1\}$}. When we change one parameter, the other three are set as default. The results are summarized in Figure \ref{fig:fig4} in Appendix \ref{sec:appenB}. From Figure \ref{fig:fig4}, we note that, for all three~\mbox{criteria},~there are only very marginal differences in the performance of $\myalg$ for different parameter settings. Thus, we conclude that, by adaptively choosing suitable penalty parameters within the algorithm, our method enjoys a robust performance to the tuning parameters.

\section{Conclusion}\label{sec:6}

We proposed an adaptive inexact Newton method, called $\myalg$, for solving constrained nonlinear optimization problems. At each step, the method applies a randomized iterative sketching solver to solve the Lagrangian Newton system inexactly, and employs an exact augmented Lagrangian merit function for selecting the stepsize via line search. The method adaptively controls the accuracy of the sketching solver, and selects suitable penalty parameters of the augmented Lagrangian. Under mild assumptions, we established the global almost sure convergence guarantee with a local linear rate for the method; and we proved~that the linear rate can be accelerated to a superlinear rate if we gradually sharpen the accuracy condition, which is achievable by simply specifying a decreasing~input~sequence. We demonstrated the superior performance and robustness to parameters of our method via experiments on benchmark nonlinear problems, constrained logistic regression, and a PDE-constrained problem. The future directions include: (i) conducting complexity analysis, (ii) applying quasi-Newton updates for Hessians, and (iii) utilizing trust-region techniques to select the inexact directions and stepsizes jointly.

\textbf{Acknowledgements.}  {\small This work was supported by Laboratory Directed Research and Development (LDRD) funding from Berkeley Lab, provided by the Director, Office of Science, of the U.S. Department of Energy under Contract No. DE-AC02-05CH11231. This work was also supported by the Intelligence Advanced Research Projects Agency (IARPA) and Army Research Office (ARO) under Contract No. W911NF-20-C-0035. The research of MK is supported in part by NSF Grant ECCS-2216912.

}
\bibliography{ref}

\begin{thebibliography}{63}
\providecommand{\natexlab}[1]{#1}
\providecommand{\url}[1]{\texttt{#1}}
\expandafter\ifx\csname urlstyle\endcsname\relax
  \providecommand{\doi}[1]{doi: #1}\else
  \providecommand{\doi}{doi: \begingroup \urlstyle{rm}\Url}\fi

\bibitem[Berahas et~al.(2020)Berahas, Bollapragada, and
  Nocedal]{Berahas2020investigation}
Berahas, A.~S., Bollapragada, R., and Nocedal, J.
\newblock An investigation of newton-sketch and subsampled newton methods.
\newblock \emph{Optimization Methods and Software}, 35\penalty0 (4):\penalty0
  661--680, feb 2020.
\newblock \doi{10.1080/10556788.2020.1725751}.

\bibitem[Berahas et~al.(2021{\natexlab{a}})Berahas, Curtis, O'Neill, and
  Robinson]{Berahas2021stochastic}
Berahas, A.~S., Curtis, F.~E., O'Neill, M.~J., and Robinson, D.~P.
\newblock A stochastic sequential quadratic optimization algorithm for
  nonlinear equality constrained optimization with rank-deficient jacobians.
\newblock \emph{arXiv preprint arXiv:2106.13015}, 2021{\natexlab{a}}.

\bibitem[Berahas et~al.(2021{\natexlab{b}})Berahas, Curtis, Robinson, and
  Zhou]{Berahas2021Sequential}
Berahas, A.~S., Curtis, F.~E., Robinson, D., and Zhou, B.
\newblock Sequential quadratic optimization for nonlinear equality constrained
  stochastic optimization.
\newblock \emph{{SIAM} Journal on Optimization}, 31\penalty0 (2):\penalty0
  1352--1379, jan 2021{\natexlab{b}}.
\newblock \doi{10.1137/20m1354556}.

\bibitem[Bertsekas(1982)]{Bertsekas1982Constrained}
Bertsekas, D.~P.
\newblock \emph{Constrained Optimization and Lagrange Multiplier Methods}.
\newblock Elsevier, 1982.
\newblock \doi{10.1016/c2013-0-10366-2}.

\bibitem[Boggs \& Tolle(1995)Boggs and Tolle]{Boggs1995Sequential}
Boggs, P.~T. and Tolle, J.~W.
\newblock Sequential quadratic programming.
\newblock \emph{Acta Numerica}, 4:\penalty0 1--51, jan 1995.
\newblock \doi{10.1017/s0962492900002518}.

\bibitem[Bollapragada et~al.(2018)Bollapragada, Byrd, and
  Nocedal]{Bollapragada2018Exact}
Bollapragada, R., Byrd, R.~H., and Nocedal, J.
\newblock Exact and inexact subsampled newton methods for optimization.
\newblock \emph{{IMA} Journal of Numerical Analysis}, 39\penalty0 (2):\penalty0
  545--578, apr 2018.
\newblock \doi{10.1093/imanum/dry009}.

\bibitem[Burke et~al.(2020)Burke, Curtis, Wang, and Wang]{Burke2020Inexact}
Burke, J.~V., Curtis, F.~E., Wang, H., and Wang, J.
\newblock Inexact sequential quadratic optimization with penalty parameter
  updates within the {QP} solver.
\newblock \emph{{SIAM} Journal on Optimization}, 30\penalty0 (3):\penalty0
  1822--1849, jan 2020.
\newblock \doi{10.1137/18m1176488}.

\bibitem[Byrd et~al.(2008)Byrd, Curtis, and Nocedal]{Byrd2008Inexact}
Byrd, R.~H., Curtis, F.~E., and Nocedal, J.
\newblock An inexact {SQP} method for equality constrained optimization.
\newblock \emph{{SIAM} Journal on Optimization}, 19\penalty0 (1):\penalty0
  351--369, jan 2008.
\newblock \doi{10.1137/060674004}.

\bibitem[Byrd et~al.(2010)Byrd, Curtis, and Nocedal]{Byrd2010inexact}
Byrd, R.~H., Curtis, F.~E., and Nocedal, J.
\newblock An inexact newton method for nonconvex equality constrained
  optimization.
\newblock \emph{Mathematical Programming}, 122\penalty0 (2):\penalty0 273--299,
  oct 2010.
\newblock \doi{10.1007/s10107-008-0248-3}.

\bibitem[Chang \& Lin(2011)Chang and Lin]{Chang2011LIBSVM}
Chang, C.-C. and Lin, C.-J.
\newblock {LIBSVM}: a library for support vector machines.
\newblock \emph{{ACM} Transactions on Intelligent Systems and Technology},
  2\penalty0 (3):\penalty0 1--27, apr 2011.
\newblock \doi{10.1145/1961189.1961199}.

\bibitem[Chen et~al.(2018)Chen, Tung, Vedula, and Mori]{Chen2018Constraint}
Chen, C., Tung, F., Vedula, N., and Mori, G.
\newblock Constraint-aware deep neural network compression.
\newblock In \emph{Proceedings of the European Conference on Computer Vision
  (ECCV)}, pp.\  400--415, 2018.

\bibitem[Cuomo et~al.(2022)Cuomo, Di~Cola, Giampaolo, Rozza, Raissi, and
  Piccialli]{Cuomo2022Scientific}
Cuomo, S., Di~Cola, V.~S., Giampaolo, F., Rozza, G., Raissi, M., and Piccialli,
  F.
\newblock Scientific machine learning through physics-informed neural networks:
  Where we are and what's next.
\newblock \emph{arXiv preprint arXiv:2201.05624}, 2022.

\bibitem[Curtis et~al.(2014{\natexlab{a}})Curtis, Jiang, and
  Robinson]{Curtis2014adaptive}
Curtis, F.~E., Jiang, H., and Robinson, D.~P.
\newblock An adaptive augmented lagrangian method for large-scale constrained
  optimization.
\newblock \emph{Mathematical Programming}, 152\penalty0 (1-2):\penalty0
  201--245, apr 2014{\natexlab{a}}.
\newblock \doi{10.1007/s10107-014-0784-y}.

\bibitem[Curtis et~al.(2014{\natexlab{b}})Curtis, Johnson, Robinson, and
  Wächter]{Curtis2014Inexact}
Curtis, F.~E., Johnson, T.~C., Robinson, D.~P., and Wächter, A.
\newblock An inexact sequential quadratic optimization algorithm for nonlinear
  optimization.
\newblock \emph{{SIAM} Journal on Optimization}, 24\penalty0 (3):\penalty0
  1041--1074, jan 2014{\natexlab{b}}.
\newblock \doi{10.1137/130918320}.

\bibitem[Curtis et~al.(2018)Curtis, Robinson, and Samadi]{Curtis2018inexact}
Curtis, F.~E., Robinson, D.~P., and Samadi, M.
\newblock An inexact regularized newton framework with a worst-case iteration
  complexity of {\textdollar} $\lbrace${\textbackslash}mathscr
  o$\rbrace$({\textbackslash}varepsilon{\^{}}$\lbrace$-3/2$\rbrace$)
  {\textdollar} for nonconvex optimization.
\newblock \emph{{IMA} Journal of Numerical Analysis}, 39\penalty0 (3):\penalty0
  1296--1327, may 2018.
\newblock \doi{10.1093/imanum/dry022}.

\bibitem[Curtis et~al.(2021{\natexlab{a}})Curtis, Robinson, Royer, and
  Wright]{Curtis2021Trust}
Curtis, F.~E., Robinson, D.~P., Royer, C.~W., and Wright, S.~J.
\newblock Trust-region newton-{CG} with strong second-order complexity
  guarantees for nonconvex optimization.
\newblock \emph{{SIAM} Journal on Optimization}, 31\penalty0 (1):\penalty0
  518--544, jan 2021{\natexlab{a}}.
\newblock \doi{10.1137/19m130563x}.

\bibitem[Curtis et~al.(2021{\natexlab{b}})Curtis, Robinson, and
  Zhou]{Curtis2021Inexact}
Curtis, F.~E., Robinson, D.~P., and Zhou, B.
\newblock Inexact sequential quadratic optimization for minimizing a stochastic
  objective function subject to deterministic nonlinear equality constraints.
\newblock \emph{arXiv preprint arXiv:2107.03512}, 2021{\natexlab{b}}.

\bibitem[Derezi{\'n}ski \& Rebrova(2022)Derezi{\'n}ski and
  Rebrova]{Derezinski2022Sharp}
Derezi{\'n}ski, M. and Rebrova, E.
\newblock Sharp analysis of sketch-and-project methods via a connection to
  randomized singular value decomposition.
\newblock \emph{arXiv preprint arXiv:2208.09585}, 2022.

\bibitem[Derezinski et~al.(2020{\natexlab{a}})Derezinski, Bartan, Pilanci, and
  Mahoney]{Derezinski2020Debiasing}
Derezinski, M., Bartan, B., Pilanci, M., and Mahoney, M.~W.
\newblock Debiasing distributed second order optimization with surrogate
  sketching and scaled regularization.
\newblock \emph{Advances in Neural Information Processing Systems},
  33:\penalty0 6684--6695, 2020{\natexlab{a}}.

\bibitem[Derezinski et~al.(2020{\natexlab{b}})Derezinski, Liang, Liao, and
  Mahoney]{Derezinski2020Precise}
Derezinski, M., Liang, F.~T., Liao, Z., and Mahoney, M.~W.
\newblock Precise expressions for random projections: Low-rank approximation
  and randomized newton.
\newblock \emph{Advances in Neural Information Processing Systems},
  33:\penalty0 18272--18283, 2020{\natexlab{b}}.

\bibitem[Derezinski et~al.(2021)Derezinski, Lacotte, Pilanci, and
  Mahoney]{Derezinski2021Newton}
Derezinski, M., Lacotte, J., Pilanci, M., and Mahoney, M.~W.
\newblock Newton-less: Sparsification without trade-offs for the sketched
  newton update.
\newblock \emph{Advances in Neural Information Processing Systems},
  34:\penalty0 2835--2847, 2021.

\bibitem[Doikov \& Richt{\'a}rik(2018)Doikov and
  Richt{\'a}rik]{Doikov2018Randomized}
Doikov, N. and Richt{\'a}rik, P.
\newblock Randomized block cubic newton method.
\newblock In \emph{International Conference on Machine Learning}, pp.\
  1290--1298. PMLR, 2018.

\bibitem[Dolan \& Mor{\'{e}}(2002)Dolan and Mor{\'{e}}]{Dolan2002Benchmarking}
Dolan, E.~D. and Mor{\'{e}}, J.~J.
\newblock Benchmarking optimization software with performance profiles.
\newblock \emph{Mathematical Programming}, 91\penalty0 (2):\penalty0 201--213,
  jan 2002.
\newblock \doi{10.1007/s101070100263}.

\bibitem[Fang et~al.(2022)Fang, Na, Mahoney, and Kolar]{Fang2022Fully}
Fang, Y., Na, S., Mahoney, M.~W., and Kolar, M.
\newblock Fully stochastic trust-region sequential quadratic programming for
  equality-constrained optimization problems.
\newblock \emph{arXiv preprint arXiv:2211.15943}, 2022.

\bibitem[Fong \& Saunders(2012)Fong and Saunders]{Fong2012CG}
Fong, D. C.-L. and Saunders, M.
\newblock {CG} versus {MINRES}: An empirical comparison.
\newblock \emph{Sultan Qaboos University Journal for Science [{SQUJS}]},
  16\penalty0 (1):\penalty0 44, apr 2012.
\newblock \doi{10.24200/squjs.vol17iss1pp44-62}.

\bibitem[Fuji et~al.(2022)Fuji, Poirion, and Takeda]{Fuji2022Randomized}
Fuji, T., Poirion, P.-L., and Takeda, A.
\newblock Randomized subspace regularized newton method for unconstrained
  non-convex optimization.
\newblock \emph{arXiv preprint arXiv:2209.04170}, 2022.

\bibitem[Gould et~al.(2014)Gould, Orban, and Toint]{Gould2014CUTEst}
Gould, N. I.~M., Orban, D., and Toint, P.~L.
\newblock {CUTEst}: a constrained and unconstrained testing environment with
  safe threads for mathematical optimization.
\newblock \emph{Computational Optimization and Applications}, 60\penalty0
  (3):\penalty0 545--557, aug 2014.
\newblock \doi{10.1007/s10589-014-9687-3}.

\bibitem[Gower et~al.(2019)Gower, Kovalev, Lieder, and
  Richt{\'a}rik]{Gower2019Rsn}
Gower, R., Kovalev, D., Lieder, F., and Richt{\'a}rik, P.
\newblock Rsn: Randomized subspace newton.
\newblock \emph{Advances in Neural Information Processing Systems}, 32, 2019.

\bibitem[Gower \& Richt{\'{a}}rik(2015)Gower and
  Richt{\'{a}}rik]{Gower2015Randomized}
Gower, R.~M. and Richt{\'{a}}rik, P.
\newblock Randomized iterative methods for linear systems.
\newblock \emph{{SIAM} Journal on Matrix Analysis and Applications},
  36\penalty0 (4):\penalty0 1660--1690, jan 2015.
\newblock \doi{10.1137/15m1025487}.

\bibitem[Gower et~al.(2021)Gower, Molitor, Moorman, and
  Needell]{Gower2021Adaptive}
Gower, R.~M., Molitor, D., Moorman, J., and Needell, D.
\newblock On adaptive sketch-and-project for solving linear systems.
\newblock \emph{{SIAM} Journal on Matrix Analysis and Applications},
  42\penalty0 (2):\penalty0 954--989, jan 2021.
\newblock \doi{10.1137/19m1285846}.

\bibitem[Gu et~al.(2017)Gu, Zhu, and Pei]{Gu2017new}
Gu, C., Zhu, D., and Pei, Y.
\newblock A new inexact {SQP} algorithm for nonlinear systems of mixed
  equalities and inequalities.
\newblock \emph{Numerical Algorithms}, 78\penalty0 (4):\penalty0 1233--1253,
  sep 2017.
\newblock \doi{10.1007/s11075-017-0421-y}.

\bibitem[Hanzely et~al.(2020)Hanzely, Doikov, Nesterov, and
  Richtarik]{Hanzely2020Stochastic}
Hanzely, F., Doikov, N., Nesterov, Y., and Richtarik, P.
\newblock Stochastic subspace cubic newton method.
\newblock In \emph{International Conference on Machine Learning}, pp.\
  4027--4038. PMLR, 2020.

\bibitem[Hintermüller et~al.(2002)Hintermüller, Ito, and
  Kunisch]{Hintermueller2002Primal}
Hintermüller, M., Ito, K., and Kunisch, K.
\newblock The primal-dual active set strategy as a semismooth newton method.
\newblock \emph{{SIAM} Journal on Optimization}, 13\penalty0 (3):\penalty0
  865--888, jan 2002.
\newblock \doi{10.1137/s1052623401383558}.

\bibitem[Kouri et~al.(2014)Kouri, Heinkenschloss, Ridzal, and van
  Bloemen~Waanders]{Kouri2014Inexact}
Kouri, D.~P., Heinkenschloss, M., Ridzal, D., and van Bloemen~Waanders, B.~G.
\newblock Inexact objective function evaluations in a trust-region algorithm
  for {PDE}-constrained optimization under uncertainty.
\newblock \emph{{SIAM} Journal on Scientific Computing}, 36\penalty0
  (6):\penalty0 A3011--A3029, jan 2014.
\newblock \doi{10.1137/140955665}.

\bibitem[Krishnapriyan et~al.(2021)Krishnapriyan, Gholami, Zhe, Kirby, and
  Mahoney]{Krishnapriyan2021Characterizing}
Krishnapriyan, A., Gholami, A., Zhe, S., Kirby, R., and Mahoney, M.~W.
\newblock Characterizing possible failure modes in physics-informed neural
  networks.
\newblock \emph{Advances in Neural Information Processing Systems},
  34:\penalty0 26548--26560, 2021.

\bibitem[Lacotte et~al.(2021)Lacotte, Wang, and Pilanci]{Lacotte2021Adaptive}
Lacotte, J., Wang, Y., and Pilanci, M.
\newblock Adaptive newton sketch: Linear-time optimization with quadratic
  convergence and effective hessian dimensionality.
\newblock In \emph{International Conference on Machine Learning}, pp.\
  5926--5936. PMLR, 2021.

\bibitem[Lewis et~al.(2012)Lewis, Vrabie, and Syrmos]{Lewis2012Optimal}
Lewis, F.~L., Vrabie, D.~L., and Syrmos, V.~L.
\newblock \emph{Optimal Control}.
\newblock John Wiley {\&} Sons, Inc., jan 2012.
\newblock \doi{10.1002/9781118122631}.

\bibitem[Liu \& Roosta(2022)Liu and Roosta]{Liu2022Newton}
Liu, Y. and Roosta, F.
\newblock A newton-mr algorithm with complexity guarantees for nonconvex smooth
  unconstrained optimization.
\newblock \emph{arXiv preprint arXiv:2208.07095}, 2022.

\bibitem[Luo et~al.(2016)Luo, Agarwal, Cesa-Bianchi, and
  Langford]{Luo2016Efficient}
Luo, H., Agarwal, A., Cesa-Bianchi, N., and Langford, J.
\newblock Efficient second order online learning by sketching.
\newblock \emph{Advances in Neural Information Processing Systems}, 29, 2016.

\bibitem[Murray et~al.(2023)Murray, Demmel, Mahoney, Erichson, Melnichenko,
  Malik, Grigori, Luszczek, Derezi{\'n}ski, Lopes,
  et~al.]{Murray2023Randomized}
Murray, R., Demmel, J., Mahoney, M.~W., Erichson, N.~B., Melnichenko, M.,
  Malik, O.~A., Grigori, L., Luszczek, P., Derezi{\'n}ski, M., Lopes, M.~E.,
  et~al.
\newblock Randomized numerical linear algebra: A perspective on the field with
  an eye to software.
\newblock \emph{arXiv preprint arXiv:2302.11474}, 2023.

\bibitem[Mutny et~al.(2020)Mutny, Derezinski, and Krause]{Mutny2020Convergence}
Mutny, M., Derezinski, M., and Krause, A.
\newblock Convergence analysis of block coordinate algorithms with
  determinantal sampling.
\newblock In \emph{International Conference on Artificial Intelligence and
  Statistics}, pp.\  3110--3120. PMLR, 2020.

\bibitem[Na(2021)]{Na2021Global}
Na, S.
\newblock Global convergence of online optimization for nonlinear model
  predictive control.
\newblock \emph{Advances in Neural Information Processing Systems},
  34:\penalty0 12441--12453, 2021.

\bibitem[Na \& Mahoney(2022)Na and Mahoney]{Na2022Asymptotic}
Na, S. and Mahoney, M.~W.
\newblock Asymptotic convergence rate and statistical inference for stochastic
  sequential quadratic programming.
\newblock \emph{arXiv preprint arXiv:2205.13687}, 2022.

\bibitem[Na et~al.(2021)Na, Anitescu, and Kolar]{Na2021fast}
Na, S., Anitescu, M., and Kolar, M.
\newblock A fast temporal decomposition procedure for long-horizon nonlinear
  dynamic programming.
\newblock \emph{arXiv preprint arXiv:2107.11560}, 2021.

\bibitem[Na et~al.(2022{\natexlab{a}})Na, Anitescu, and Kolar]{Na2022adaptive}
Na, S., Anitescu, M., and Kolar, M.
\newblock An adaptive stochastic sequential quadratic programming with
  differentiable exact augmented lagrangians.
\newblock \emph{Mathematical Programming}, pp.\  1--71, jun 2022{\natexlab{a}}.
\newblock \doi{10.1007/s10107-022-01846-z}.

\bibitem[Na et~al.(2022{\natexlab{b}})Na, Derezi{\'{n}}ski, and
  Mahoney]{Na2022Hessian}
Na, S., Derezi{\'{n}}ski, M., and Mahoney, M.~W.
\newblock Hessian averaging in stochastic newton methods achieves superlinear
  convergence.
\newblock \emph{Mathematical Programming}, dec 2022{\natexlab{b}}.
\newblock \doi{10.1007/s10107-022-01913-5}.

\bibitem[Na et~al.(2023)Na, Anitescu, and Kolar]{Na2023Inequality}
Na, S., Anitescu, M., and Kolar, M.
\newblock Inequality constrained stochastic nonlinear optimization via
  active-set sequential quadratic programming.
\newblock \emph{Mathematical Programming}, mar 2023.
\newblock \doi{10.1007/s10107-023-01935-7}.

\bibitem[Nocedal \& Wright(2006)Nocedal and Wright]{Nocedal2006Numerical}
Nocedal, J. and Wright, S.~J.
\newblock \emph{Numerical Optimization (2nd Edition)}.
\newblock Springer-Verlag, 2006.
\newblock \doi{10.1007/b98874}.

\bibitem[Pilanci \& Wainwright(2017)Pilanci and Wainwright]{Pilanci2017Newton}
Pilanci, M. and Wainwright, M.~J.
\newblock Newton sketch: A near linear-time optimization algorithm with
  linear-quadratic convergence.
\newblock \emph{{SIAM} Journal on Optimization}, 27\penalty0 (1):\penalty0
  205--245, jan 2017.
\newblock \doi{10.1137/15m1021106}.

\bibitem[Qu et~al.(2016)Qu, Richt{\'a}rik, Tak{\'a}c, and Fercoq]{Qu2016SDNA}
Qu, Z., Richt{\'a}rik, P., Tak{\'a}c, M., and Fercoq, O.
\newblock Sdna: stochastic dual newton ascent for empirical risk minimization.
\newblock In \emph{International Conference on Machine Learning}, pp.\
  1823--1832. PMLR, 2016.

\bibitem[Richt{\'{a}}rik \& Tak{\'{a}}{\v{c}}(2014)Richt{\'{a}}rik and
  Tak{\'{a}}{\v{c}}]{Richtarik2014Iteration}
Richt{\'{a}}rik, P. and Tak{\'{a}}{\v{c}}, M.
\newblock Iteration complexity of randomized block-coordinate descent methods
  for minimizing a composite function.
\newblock \emph{Mathematical Programming}, 144\penalty0 (1-2):\penalty0 1--38,
  dec 2014.
\newblock \doi{10.1007/s10107-012-0614-z}.

\bibitem[Roosta et~al.(2022)Roosta, Liu, Xu, and Mahoney]{Roosta2022Newton}
Roosta, F., Liu, Y., Xu, P., and Mahoney, M.~W.
\newblock Newton-{MR}: Inexact newton method with minimum residual sub-problem
  solver.
\newblock \emph{{EURO} Journal on Computational Optimization}, 10:\penalty0
  100035, 2022.
\newblock \doi{10.1016/j.ejco.2022.100035}.

\bibitem[Roosta-Khorasani \& Mahoney(2018)Roosta-Khorasani and
  Mahoney]{RoostaKhorasani2018Sub}
Roosta-Khorasani, F. and Mahoney, M.~W.
\newblock Sub-sampled newton methods.
\newblock \emph{Mathematical Programming}, 174\penalty0 (1-2):\penalty0
  293--326, nov 2018.
\newblock \doi{10.1007/s10107-018-1346-5}.

\bibitem[Royer \& Wright(2018)Royer and Wright]{Royer2018Complexity}
Royer, C.~W. and Wright, S.~J.
\newblock Complexity analysis of second-order line-search algorithms for smooth
  nonconvex optimization.
\newblock \emph{{SIAM} Journal on Optimization}, 28\penalty0 (2):\penalty0
  1448--1477, jan 2018.
\newblock \doi{10.1137/17m1134329}.

\bibitem[Royer et~al.(2019)Royer, O'Neill, and Wright]{Royer2019Newton}
Royer, C.~W., O'Neill, M., and Wright, S.~J.
\newblock A newton-{CG} algorithm with complexity guarantees for smooth
  unconstrained optimization.
\newblock \emph{Mathematical Programming}, 180\penalty0 (1-2):\penalty0
  451--488, jan 2019.
\newblock \doi{10.1007/s10107-019-01362-7}.

\bibitem[Saad \& Schultz(1986)Saad and Schultz]{Saad1986GMRES}
Saad, Y. and Schultz, M.~H.
\newblock {GMRES}: A generalized minimal residual algorithm for solving
  nonsymmetric linear systems.
\newblock \emph{{SIAM} Journal on Scientific and Statistical Computing},
  7\penalty0 (3):\penalty0 856--869, jul 1986.
\newblock \doi{10.1137/0907058}.

\bibitem[Scott \& T{\r{u}}ma(2022)Scott and T{\r{u}}ma]{Scott2022Solving}
Scott, J. and T{\r{u}}ma, M.
\newblock Solving large linear least squares problems with linear equality
  constraints.
\newblock \emph{{BIT} Numerical Mathematics}, 62\penalty0 (4):\penalty0
  1765--1787, jul 2022.
\newblock \doi{10.1007/s10543-022-00930-2}.

\bibitem[Strohmer \& Vershynin(2008)Strohmer and
  Vershynin]{Strohmer2008Randomized}
Strohmer, T. and Vershynin, R.
\newblock A randomized kaczmarz algorithm with exponential convergence.
\newblock \emph{Journal of Fourier Analysis and Applications}, 15\penalty0
  (2):\penalty0 262--278, apr 2008.
\newblock \doi{10.1007/s00041-008-9030-4}.

\bibitem[Wächter \& Biegler(2006)Wächter and
  Biegler]{Waechter2006implementation}
Wächter, A. and Biegler, L.~T.
\newblock On the implementation of an interior-point filter line-search
  algorithm for large-scale nonlinear programming.
\newblock \emph{Mathematical Programming}, 106\penalty0 (1):\penalty0 25--57,
  apr 2006.
\newblock \doi{10.1007/s10107-004-0559-y}.

\bibitem[Yao et~al.(2021)Yao, Xu, Roosta, and Mahoney]{Yao2021Inexact}
Yao, Z., Xu, P., Roosta, F., and Mahoney, M.~W.
\newblock Inexact nonconvex newton-type methods.
\newblock \emph{{INFORMS} Journal on Optimization}, 3\penalty0 (2):\penalty0
  154--182, jan 2021.
\newblock \doi{10.1287/ijoo.2019.0043}.

\bibitem[Yao et~al.(2022)Yao, Xu, Roosta, Wright, and Mahoney]{Yao2022Inexact}
Yao, Z., Xu, P., Roosta, F., Wright, S.~J., and Mahoney, M.~W.
\newblock Inexact newton-{CG} algorithms with complexity guarantees.
\newblock \emph{{IMA} Journal of Numerical Analysis}, aug 2022.
\newblock \doi{10.1093/imanum/drac043}.

\bibitem[Yuan et~al.(2022)Yuan, Lazaric, and Gower]{Yuan2022Sketched}
Yuan, R., Lazaric, A., and Gower, R.~M.
\newblock Sketched newton--raphson.
\newblock \emph{{SIAM} Journal on Optimization}, 32\penalty0 (3):\penalty0
  1555--1583, jul 2022.
\newblock \doi{10.1137/21m139788x}.

\bibitem[Zavala \& Anitescu(2014)Zavala and Anitescu]{Zavala2014Scalable}
Zavala, V.~M. and Anitescu, M.
\newblock Scalable nonlinear programming via exact differentiable penalty
  functions and trust-region newton methods.
\newblock \emph{{SIAM} Journal on Optimization}, 24\penalty0 (1):\penalty0
  528--558, jan 2014.
\newblock \doi{10.1137/120888181}.

\end{thebibliography}
\bibliographystyle{icml2023}

\newpage
\appendix
\onecolumn

\begin{center}
\large\textbf{Appendix: Constrained Optimization via Exact Augmented Lagrangian and \\ Randomized Iterative Sketching}
\end{center}

\section{Proofs of Section \ref{sec:3}}\label{sec:appenA}

\subsection{Proof of Lemma \ref{lemma3.4}}

Throughout the proof, we fix $k\ge 0$ and suppose the algorithm reaches $\vz_k=(\vx_k,\vlambda_k)$. We note that
\begin{equation}\label{eq:irs_reduced1}
\incztkjn \stackrel{\eqref{eq:irs}}{=} \incztkj-\KKT S_{k,j}(S_{k,j}^T\KKT^2S_{k,j})^{\dagger}S_{k,j}^T\vr_{k,j} \;\; \stackrel{\mathclap{\substack{\eqref{eq:kkt_sys_reduced}, \eqref{eq:residual}}}}{=}\;\; \incztkj-W_{k,j}(\incztkj-\inczk),
\end{equation}
where $W_{k,j}=\KKT S_{k,j}(S_{k,j}^T\KKT^2S_{k,j})^{\dagger}S_{k,j}^T\KKT$. Since $W_{k,j}$ is an orthogonal projection onto $row(S_{k,j}^T\KKT)$, we rewrite~\eqref{eq:irs_reduced1} as
\begin{equation}\label{eq:irs_reduced2}
\incztkjn-\inczk=\incztkj-\inczk-Q_{k,j}Q^T_{k,j}(\incztkj-\inczk).
\end{equation}
By Assumption \ref{assum2}, $\KKT$ is invertible (see \citet[Lemma 16.1]{Nocedal2006Numerical}). Thus, Assumption \ref{assum3} implies that $\prob(S^T\KKT\vz\neq \vzero)\geq \pi$ for any $\vz\in\mR^{n+m}\backslash\{\vzero\}$. Given the relationship between $row(S_{k,j}^T\KKT)$ and $Q_{k,j}$, we further have $\prob(Q_{k,j}^T\vz\neq \vzero) = \prob(S_{k,j}^T\KKT\vz\neq \vzero)\ge\pi$ for any $\vz\in\mR^{n+m}\backslash\{\vzero\}$. Since $\{Q_{k,j}\}_j$ are independent and identically distributed, conditional on $\{j_{l-1}^{k}<\infty\}$, the probability that $\dim(\cup_{i=0}^{t+1}col(Q_{k,j_{l-1}^k+i}))$ grows relative to $\dim(\cup_{i=0}^{t}col(Q_{k,j_{l-1}^k+i}))$~when $\dim(\cup_{i=0}^{t}col(Q_{k,j_{l-1}^k+i}))<n+m$, is at least $\pi$. As a result, the probability that the event $\{\dim(\cup_{i=0}^{t+1}col(Q_{k,j_{l-1}^k+i}))>  \dim(\cup_{i=0}^{t}col(Q_{k,j_{l-1}^k+i}))\}$ happens $n+m$ times in $N$ iterations with $N\ge n+m$ is dominated by a negative binomial distribution. In particular, for $N\ge n+m$,
\begin{equation*}
\prob(j_l^{k}-1=N+j_{l-1}^{k}|j_{l-1}^{k}<\infty) \le\begin{pmatrix}N-1\\n+m-1\end{pmatrix}(1-\pi)^{N-n-m}\pi^{n+m}.
\end{equation*}
Taking $N\rightarrow\infty$, we have for any $l\geq 1$,
\begin{equation*}
\prob(j_l^{k}=\infty|j_{l-1}^{k}<\infty)=0.
\end{equation*}
Thus, for any $l\in\mN$, $P(j_l^{k}<\infty|j_{l-1}^{k}<\infty)=1$. Furthermore, we have
\begin{align*}
\prob\left(\cap_{l=1}^L\{j_l^{k}<\infty\}\right)&=\prob(j_1^{k}<\infty)\times \prob(j_2^{k}<\infty|j_1^{k}<\infty)\times \cdots\times \prob(j_L^{k}<\infty|j_{L-1}^{k}<\infty,\dots,j_1^{k}<\infty)\\
&=\prob(j_1^{k}<\infty)\times \prob(j_2^{k}<\infty|j_1^{k}<\infty)\times \cdots\times \prob(j_L^{k}<\infty|j_{L-1}^{k}<\infty)=1.
\end{align*}
This completes the proof.

\subsection{Proof of Lemma \ref{lemma3.5}}

Let us denote $\vq_{k,j,h}$ to be the $h$-th column of $Q_{k,j}$ for $1\leq h\leq d$. We know for any $1\leq l\le L$,
\begin{equation*}
\incztkjlk-\inczk \stackrel{\eqref{eq:irs_reduced2}}{=}\left(\Pi_{j=j_{l-1}^k}^{j_l^k-1}\Pi_{h=1}^d\left(I-\vq_{k,j,h}\vq^T_{k,j,h}\right)\right)(\incztkjlbk-\inczk).
\end{equation*}
Taking $\ell_2$ norm on both sides, we obtain
\begin{equation*}
\Vert\incztkjlk-\inczk\Vert\le\Vert\Pi_{j=j_{l-1}^k}^{j_l^k-1}\Pi_{h=1}^d\left(I-\vq_{k,j,h}\vq^T_{k,j,h}\right)\Vert\cdot \Vert\incztkjlbk-\inczk\Vert.
\end{equation*}
Let $\mathcal{F}_{k,l}$ denote the set of all matrices $F_{k,l}$, whose columns $\{f_{k,l,1},\dots,f_{k,l,n+m}\}$ form a maximal linearly independent subset of $\{\vq_{k,j_{l-1}^k,1},\dots,\vq_{k,j_l^k-1,d}\}$. \citet[Theorem 4.1]{Patel2021Implicit} implies that 
\begin{equation*}
\Vert\Pi_{j=j_{l-1}^k}^{j_l^k-1}\Pi_{h=1}^d\left(I-\vq_{k,j,h}\vq^T_{k,j,h}\right)\Vert\le\sqrt{1-\min_{F_{k,l}\in\mathcal{F}_{k,l}}\det(F_{k,l}^TF_{k,l})} \eqqcolon \gamma_{k,l}.
\end{equation*}
Thus, we have $\Vert\incztkjlk-\inczk\Vert\le\gamma_{k,l} \Vert\incztkjlbk-\inczk\Vert$. By Hadamard's inequality and the fact that $F_{k,l}^TF_{k,l}$ is positive definite, we know $0\leq \gamma_{k,l}<1$ (note that $\mathcal{F}_{k,l}$ is a finite set). Since $S_{k,j}\sim S$, $iid$, we know the distribution of $\gamma_{k,l}$ is independent of $l$. In particular, we let $\mathcal{Q}_{k,l}=\{\vq_{k,j_{l-1}^k,1},\dots,\vq_{k,j_l^k-1,d}\}$ and have $\mathcal{Q}_{k,1},\dots,\mathcal{Q}_{k,L}\sim\mathcal{Q}_k$, $iid$. This implies $\gamma_{k,1},\dots,\gamma_{k,L}\sim\gamma_k$, $iid$, and completes the proof.

\subsection{Proof of Lemma~\ref{lemma10}}

Throughout the proof, we fix $k\ge 0$ and suppose that the algorithm reaches $\vz_k = (\vx_k,\vlambda_k)$. From Lemma \ref{lemma3.5}, we know $\prob(\gamma_k=1)=0$; thus, there exists $\tau_k\in(0,1)$ such that $\prob(\gamma_k\le\tau_k)>0$. Let us define $\pi_k\coloneqq \prob(\gamma_k\le\tau_k)$. By Assumptions \ref{assum1} and \ref{assum2}, we know $\normA\le\Upsilon_{\Gamma}$ and $\Psi_k\le\Psi$ for some $\Upsilon_{\Gamma}$ and $\Psi>0$. Thus, $\theta_k\delta_k/(\normA^2\Psi_k^2)>0$. Let $\bar{N}$ be the smallest positive integer such that $\bar{N}\ge\log\{\theta_k\delta_k/(\normA^2\Psi_k^2)\}/\log(\tau_k)$. Then, we have $\tau_k^{\bar{N}}\le\theta_k\delta_k/(\normA^2\Psi_k^2)$. We now consider a procedure where for each iteration $l$, we generate $\gamma_{k,l}\sim \gamma_k$ independently. Let $L_k$ be the iteration such that 
\begin{equation*}
I\{\gamma_{k,1}\le\tau_k\}+\cdots+I\{\gamma_{k,L_k}\le\tau_k\}=\bar{N}.
\end{equation*}
We note the probability that the event $\{\gamma_{k,l}\le\tau_k\}$ happens $\bar{N}$ times in $N$ iterations with $N\ge\bar{N}$ is dominated by a negative binomial distribution. In particular, for $N\ge\bar{N}$,
\begin{equation*}
P(L_k=N)\le\begin{pmatrix} N-1 \\ \bar{N}-1 \end{pmatrix}(1-\pi_k)^{N-\bar{N}}\pi_k^{\bar{N}}.
\end{equation*}
Taking $N\rightarrow\infty$, we have $P(L_k=\infty)=0$. Thus, $L_k$ is finite with probability one. We now apply Lemma \ref{lemma3.5}. We have
\begin{equation}\label{eq:inner}
\Vert\incztkjLk-\inczk\Vert\le \left(\Pi_{l=1}^{L_k}\gamma_{k,l}\right)\Vert\incztkzero-\inczk\Vert=\left(\Pi_{l=1}^{L_k}\gamma_{k,l}\right)\normexactz.    
\end{equation}
By Assumption \ref{assum2} and the fact that $\|(G_kG_k^T)^{-1}\|\le 1/\sigma_{1,k}^2$ with $\sigma_{1,k}$ being the least singular value of $G_k$, we apply \citet[Lemma 5.1]{Na2021fast} and have
\begin{equation}\label{eq:Psi}
\normAI \le\frac{1}{\xi_{B}}+\frac{2}{\sigma_{1,k}}\left(1 + \frac{\normB}{\xi_{B}}\right) + \frac{1}{\sigma_{1,k}^2}\left(\normB+\frac{\normB^2}{\xi_{B}}\right)\le\frac{7(\normB^2\vee 1)}{(\xi_{B}\wedge1)(\sigma_{1,k}^2\wedge1)}\stackrel{\eqref{eq:aux}}{\leq}\Psi_k.
\end{equation}
Thus, we get
\begin{equation*}
\resizebox{0.99\linewidth}{!}{$\Vert\vr_{k,j_{L_k}^k}\Vert\stackrel{\eqref{eq:residual}}{=}\Vert \KKT\incztkjLk+\nablal\Vert\stackrel{\eqref{eq:kkt_sys_reduced}}{\le}\normA\Vert\incztkjLk-\incztk\Vert\stackrel{\eqref{eq:inner}}{\le}\left(\Pi_{l=1}^{L_k}\gamma_{k,l}\right)\normA\normexactz\;\;\stackrel{\mathclap{\substack{\eqref{eq:kkt_sys_reduced}, \eqref{eq:Psi}}}}{\le}\;\;\left(\Pi_{l=1}^{L_k}\gamma_{k,l}\right)\normA\Psi_k\normnablal.$}
\end{equation*}
Furthermore, we have
\begin{align*}
\resizebox{0.99\linewidth}{!}{$\big\{I\{\gamma_{k,1}\le\tau_k\}+\cdots+I\{\gamma_{k,L_k}\le\tau_k\}=\bar{N}\big\}
\Rightarrow \big\{\Pi_{l=1}^{L_k}\gamma_{k,l}\le\theta_k\delta_k/(\normA^2\Psi_k^2)\big\}
\Rightarrow \big\{\Vert\vr_{k,j_{L_k}^k}\Vert\le\theta_k\delta_k\normnablal/(\normA\Psi_k)\big\}.$}
\end{align*}
Finally, we let $J_k=j_{L_k}^k$ and complete the proof.

\subsection{Proof of Lemma~\ref{lemma11}}

Throughout the proof, we fix $k\ge 0$ and suppose that the algorithm reaches $\vz_k = (\vx_k,\vlambda_k)$. We first observe that
\begin{align}
\left\Vert\vr_{k,j}\right\Vert\le\theta_k\delta_k\normnablal/(\normA\Psi_k)
&\stackrel{\mathclap{\substack{\eqref{eq:kkt_sys_reduced}\\\eqref{eq:residual}}}}{\Longrightarrow} \Psi_k\Vert\KKT(\incztkj-\inczk)\Vert\le\theta_k\delta_k\normnablal/\normA\nonumber\\
&\stackrel{\mathclap{\substack{\eqref{eq:Psi}}}}{\Longrightarrow} \Vert\incztkj-\inczk\Vert\le\theta_k\delta_k\normnablal/\normA\nonumber\\
&\stackrel{\mathclap{\substack{\eqref{eq:kkt_sys_reduced}}}}{\Longrightarrow} \Vert\incztkj-\inczk\Vert\le\theta_k\delta_k\normA\normexactz/\normA\nonumber\\
&\Longrightarrow \Vert\incztkj-\inczk\Vert\le\theta_k\delta_k\normexactz\label{eq:errorsteptheta}\\
&\Longrightarrow \Vert\incztkj-\inczk\Vert\le\delta_k\normexactz\label{eq:errorstep}.
\end{align}
By Assumption~\ref{assum1}, we know there exist constants $\Upsilon_H,\Upsilon_G>0$ such that $\|H_k\|\le\Upsilon_H$ and $\|G_k\|\le\Upsilon_G$. We let $\Upsilon>0$ be a constant such that $\Upsilon_H\vee\Upsilon_G\vee\Upsilon_B\le \Upsilon$. We now divide $(\nablaL)^T\incztkj$ into two terms as follows: 
\begin{equation}\label{eq:descent1}
(\nablaL)^T\incztkj=\begin{pmatrix}\nablaxL\\\nablalL\end{pmatrix}^T\incxltkj = \begin{pmatrix}\nablaxL\\\nablalL\end{pmatrix}^T\incxlk + \begin{pmatrix}\nablaxL\\\nablalL\end{pmatrix}^T\incedj.
\end{equation}
For the first term, we have
\begin{align*}
\begin{pmatrix}\nablaxL\\\nablalL\end{pmatrix}^T\incxlk
&\stackrel{\mathclap{\eqref{eq:opt_merit}}}{=}\begin{pmatrix}(I+\eta_{2,k}H_k)\nablaxl+\eta_{1,k}G_k^Tc_k \\ c_k+\eta_{2,k} G_k\nablaxl\end{pmatrix}^T\incxlk\\ 
& =\incxlk^T\begin{pmatrix} I+\eta_{2,k} H_k&\eta_{1,k} G_k^T \\ \eta_{2,k} G_k & I\end{pmatrix}\begin{pmatrix}\nablaxl \\ c_k\end{pmatrix}\\ 
&\stackrel{\mathclap{\eqref{eq:kkt_sys_reduced}}}{=}-\incxlk^T\begin{pmatrix} I+\eta_{2,k} H_k&\eta_{1,k} G_k^T \\ \eta_{2,k} G_k & I\end{pmatrix}\begin{pmatrix} B_k & G_k^T \\ G_k & \vzero\end{pmatrix}\incxlk\\
&=-\incxlk^T\begin{pmatrix}(I+\eta_{2,k}H_k)B_k+\eta_{1,k}G_k^TG_k&(I+\eta_{2,k}H_k)G_k^T \\ G_k(I+\eta_{2,k}B_k) & \eta_{2,k} G_k G_k^T\end{pmatrix}\incxlk\\
&=-\incxk^T\left((I+\eta_{2,k}H_k)B_k+\frac{\eta_{1,k}}{2}G_k^TG_k\right)\incxk-\frac{\eta_{1,k}}{2}\incxk^T G_k^TG_k\incxk\\
&\quad-\eta_{2,k}\inclk^TG_kG_k^T{\inclk}-\inclk^TG_k\left(2I+\eta_{2,k}(B_k+H_k)\right)\incxk\\
&=-\incxk^T\left((I+\eta_{2,k}H_k)B_k+\frac{\eta_{1,k}}{2}G_k^TG_k\right)\incxk-\frac{\eta_{1,k}}{2}\left\Vert G_k\incxk\right\Vert^2-\eta_{2,k}\left\Vert G_k^T{\inclk}\right\Vert^2\\
&\quad -\inclk^TG_k\left(2I+\eta_{2,k}(B_k+H_k)\right)\incxk.
\end{align*}
Using \eqref{eq:kkt_sys_reduced}, we have $G_k\incx=-c_k$ and $G_k^T\inclk=-\left(B_k\incxk+\nablaxl\right)$. Combining these equations with the above~display, we have
\begin{align*}
\begin{pmatrix}\nablaxL\\\nablalL\end{pmatrix}^T&\incxlk\\
&=-\incxk^T\left((I+\eta_{2,k}H_k)B_k+\frac{\eta_{1,k}}{2}G_k^TG_k\right)\incxk-\frac{\eta_{1,k}}{2}\normc^2-\inclk^T G_k\left(2I+\eta_{2,k}(B_k+H_k)\right)\incxk \\
&\quad -\eta_{2,k}\left\Vert B_k{\incxk}+\nablaxl\right\Vert^2 \\
&=-\incxk^T\left((I+\eta_{2,k}H_k)B_k+\frac{\eta_{1,k}}{2} G_k^TG_k\right)\incxk-\frac{\eta_{1,k}}{2}\normc^2-\inclk^T G_k\left(2I+\eta_{2,k}(B_k+H_k)\right)\incxk \\
&\quad -\frac{\eta_{2,k}}{2}\normnablaxl^2+\frac{\eta_{2,k}}{2}\normnablaxl^2-\eta_{2,k}\left\Vert B_k{\incxk}+\nablaxl\right\Vert^2.
\end{align*}
For the last two terms in the above display, we have
\begin{align*}
\frac{\eta_{2,k}}{2}\normnablaxl^2-&\eta_{2,k}\left\Vert B_k{\incxk}+\nablaxl\right\Vert^2\\
&=-\eta_{2,k}\left\Vert B_k{\incxk}\right\Vert^2-2\eta_{2,k}\incxk^T B_k\nablaxl-\frac{\eta_{2,k}}{2}\normnablaxl^2\\
&\stackrel{\mathclap{\eqref{eq:kkt_sys_reduced}}}{=}-\eta_{2,k}\left\Vert B_k{\incxk}\right\Vert^2+2\eta_{2,k}\incxk^T B_k\left(B_k\incxk+G_k^T\inclk\right)-\frac{\eta_{2,k}}{2}\left\Vert B_k\incxk+G_k^T\inclk\right\Vert^2\\
&=\eta_{2,k}\incxk^T B_k^2{\incxk}+\eta_{2,k}\incxk^TB_k G_k^T\inclk-\frac{\eta_{2,k}}{2}\left\Vert B_k{\incxk}\right\Vert^2-\frac{\eta_{2,k}}{2}\left\Vert G_k^T\inclk\right\Vert^2\\
&\le\eta_{2,k}\incxk^TB_k^2{\incxk}+\eta_{2,k}\incxk^T B_k G_k^T\inclk-\frac{\eta_{2,k}}{2}\left\Vert G_k^T\inclk\right\Vert^2.
\end{align*}
Combining the above two displays, we get
\begin{align*}
\begin{pmatrix}\nablaxL\\\nablalL\end{pmatrix}^T&\incxlk\\ &\le-\incxk^T\left((I+\eta_{2,k} H_k)B_k+\frac{\eta_{1,k}}{2} G_k^TG_k\right)\incxk-\frac{\eta_{1,k}}{2}\normc^2-\inclk^T G_k\left(2I+\eta_{2,k}(B_k+H_k)\right)\incxk\\
&\quad -\frac{\eta_{2,k}}{2}\normnablaxl^2 +\eta_{2,k}\incxk^T  B_k^2\incxk+\eta_{2,k}\inclk^TG_kB_k\incxk-\frac{\eta_{2,k}}{2}\Vert G_k^T\inclk\Vert_2^2.
\end{align*}
Assuming $\eta_{1,k}\ge\eta_{2,k}$ at the moment and using Cauchy-Schwarz inequality, we get
\begin{align*}
\begin{pmatrix}\nablaxL\\\nablalL\end{pmatrix}^T\incxlk
&\le-\frac{\eta_{2,k}}{2}\normnablal^2-\incxk^T\left((I+\eta_{2,k}(H_k-B_k))B_k+\frac{\eta_{1,k}}{2} G_k^TG_k\right)\incxk-\frac{\eta_{2,k}}{2}\left\Vert G_k^T\inclk\right\Vert^2\\
&\quad -\inclk^T G_k\left(2I+\eta_{2,k}H_k\right)\incxk\\
&=-\frac{\eta_{2,k}}{2}\normnablal^2-\eta_{2,k}\incxk^T(H_k-B_k)B_k\incxk-\incxk^T B_k\incxk-\frac{\eta_{1,k}}{2}\incxk^T G_k^TG_k\incxk\\
&\quad -\frac{\eta_{2,k}}{2}\left\Vert G_k^T\inclk\right\Vert^2 -2\inclk^T G_k\incxk -\eta_{2,k}\inclk^TG_kH_k\incxk\\
&\le-\frac{\eta_{2,k}}{2}\normnablal^2+\eta_{2,k}\left\Vert(H_k-B_k)\incxk\right\Vert\left\Vert B_k\incxk\right\Vert-\incxk^T B_k\incxk-\frac{\eta_{1,k}}{2}\incxk^T G_k^TG_k\incxk\\
&\quad -\frac{\eta_{2,k}}{2}\left\Vert G_k^T\inclk\right\Vert^2 +2\left\Vert{\inclk}\right\Vert\left\Vert G_k\incxk\right\Vert+\eta_{2,k}\Upsilon \left\Vert { G_k^T\inclk}\right\Vert\left\Vert\incxk\right\Vert\\
&\le-\frac{\eta_{2,k}}{2}\normnablal^2+2\eta_{2,k}\Upsilon^2\left\Vert \incxk\right\Vert^2-\incxk^T\left(B_k+\frac{\eta_{1,k}}{2}G_k^TG_k\right)\incxk-\frac{\eta_{2,k}}{2}\left\Vert G_k^T\inclk\right\Vert^2\\
&\quad+2\left\Vert{\inclk}\right\Vert\left\Vert G_k\incxk\right\Vert+\eta_{2,k}\Upsilon \left\Vert { G_k^T\inclk}\right\Vert\left\Vert\incxk\right\Vert.
\end{align*}
Applying Young's inequalities for the last two terms, we obtain 
\begin{align*}
&2\left\Vert\inclk\right\Vert\left\Vert G_k\incxk\right\Vert\le\frac{\eta_{2,k}\xi_G}{8}\left\Vert\inclk\right\Vert^2+\frac{8}{\eta_{2,k}\xi_G}\left\Vert G_k\incxk\right\Vert^2,\\
&\eta_{2,k}\Upsilon\left\Vert G_k^T{\inclk}\right\Vert\left\Vert\incxk\right\Vert\le\frac{\eta_{2,k}}{4}\left\Vert G_k^T\inclk\right\Vert^2+\eta_{2,k}\Upsilon^2\left\Vert\incxk\right\Vert^2.
\end{align*}
Combining the above two displays and using Assumption \ref{assum2}, we get
\begin{align}\label{eq:descent2}
&\begin{pmatrix}\nablaxL\\\nablalL\end{pmatrix}^T\incxlk\nonumber\\
&\quad\le-\frac{\eta_{2,k}}{2}\normnablal^2+3\eta_{2,k}\Upsilon^2\left\Vert \incxk\right\Vert^2-\frac{\eta_{2,k}}{4}\left\Vert G_k^T\inclk\right\Vert^2+\frac{\eta_{2,k}\xi_G}{8}\left\Vert\inclk\right\Vert^2+\frac{8}{\eta_{2,k}\xi_G}\left\Vert G_k\incxk\right\Vert^2\nonumber\\
&\qquad\quad -\incxk^T\left( B_k+\frac{\eta_{1,k}}{2} G_k^T G_k\right)\incxk\nonumber\\
&\quad\le-\frac{\eta_{2,k}}{2}\normnablal^2-\frac{\eta_{2,k}\xi_G}{8}\left\Vert \inclk\right\Vert^2+3\eta_{2,k}\Upsilon^2\left\Vert\incxk\right\Vert^2-\incxk^T\left(B_k+\left(\frac{\eta_{1,k}}{2}-\frac{8}{\eta_{2,k}\xi_G}\right)G_k^TG_k\right)\incxk.
\end{align}
In order to bound the last two terms of the above inequality, we decompose $\incxk$ as $\incxk=\Delta\vv_k+\Delta\vu_k$ where $\Delta\vv_k=G_k^T\Delta\bar{\vv}_k$ for some $\Delta\bar{\vv}_k$ and $\Delta\vu_k$ satisfies $G_k\Delta\vu_k=\vzero$. Then, we have $\Vert\incxk\Vert^2=\Vert\Delta\vv_k\Vert^2+\Vert\Delta\vu_k\Vert^2$ and $\Vert\Delta\vv_k\Vert^2= \Vert G_k^T\Delta\bar{\vv}_k\Vert^2\le\Upsilon^2\Vert \Delta\bar{\vv}_k\Vert^2$. Thus, $\left\Vert G_k\incxk\right\Vert^2=\left\Vert G_k\Delta\vv_k\right\Vert^2=\left\Vert G_kG_k^T\Delta\bar{\vv}_k\right\Vert^2\ge \xi_G^2\left\Vert\Delta\bar{\vv}_k\right\Vert^2\ge(\xi_G^2/\Upsilon^2)\left\Vert\Delta\vv_k\right\Vert^2$. Assuming $\eta_{1,k} \ge 16/(\eta_{2,k}\xi_G)$ at the moment and using Assumption~\ref{assum2} and Cauchy-Schwarz inequality, we get
\begin{align*}
3\eta_{2,k}\Upsilon^2\left\Vert \incxk\right\Vert^2 & - \incxk^T\left(B_k+\left(\frac{\eta_{1,k}}{2}-\frac{8}{\eta_{2,k}\xi_G}\right)G_k^TG_k\right)\incxk\\
&=3\eta_{2,k}\Upsilon^2\left\Vert \incxk\right\Vert^2-\Delta\vu_k^T B_k{\Delta\vu_k}-2\Delta\vu_k^T B_k{\Delta\vv_k}-\Delta\vv_k^T B_k{\Delta\vv_k}-\left(\frac{\eta_{1,k}}{2}-\frac{8}{\eta_{2,k}\xi_G}\right)\left\Vert G_k\incxk\right\Vert^2\\
&\le3\eta_{2,k}\Upsilon^2\left\Vert\incxk\right\Vert^2-\xi_B\left\Vert \Delta\vu_k\right\Vert^2+2\Upsilon\left\Vert\Delta\vu_k\right\Vert\left\Vert \Delta\vv_k\right\Vert+\Upsilon\left\Vert \Delta\vv_k\right\Vert^2-\left(\frac{\eta_{1,k}}{2}-\frac{8}{\eta_{2,k}\xi_G}\right)\dfrac{\xi_G^2}{\Upsilon^2}\left\Vert\Delta\vv_k\right\Vert^2\\
&\le\left(3\eta_{2,k}\Upsilon^2-\xi_B\right)\left\Vert\incxk\right\Vert^2+2\Upsilon\left\Vert \Delta\vu_k\right\Vert \left\Vert\Delta\vv_k\right\Vert+\left(\xi_B+\Upsilon\right)\left\Vert \Delta\vv_k\right\Vert^2-\left(\frac{\eta_{1,k}\xi_G^2}{2\Upsilon^2}-\frac{8\xi_G}{\eta_{2,k}\Upsilon^2}\right)\left\Vert\Delta\vv_k\right\Vert^2.
\end{align*}
Applying Young's inequality for the second term, we obtain
\begin{equation*}
2\Upsilon\left\Vert\Delta\vu_k\right\Vert\left\Vert\Delta\vv_k\right\Vert\le\frac{\xi_B}{2}\left\Vert\Delta\vu_k\right\Vert^2+\frac{2\Upsilon^2}{\xi_B}\left\Vert\Delta\vv_k\right\Vert^2\le\frac{\xi_B}{2}\left\Vert\Delta\vx_k\right\Vert^2+\frac{2\Upsilon^2}{\xi_B}\left\Vert\Delta\vv_k\right\Vert^2.
\end{equation*}
Combining the above two displays, we get
\begin{align*}
3\eta_{2,k}\Upsilon^2\left\Vert \incxk\right\Vert^2-&\incxk^T\left(B_k+\left(\frac{\eta_{1,k}}{2}-\frac{8}{\eta_{2,k}\xi_G}\right) G_k^TG_k\right)\incxk\\
&\le\left(3\eta_{2,k}\Upsilon^2-\frac{\xi_B}{2}\right)\left\Vert \incxk\right\Vert^2+\frac{2\Upsilon^2}{\xi_B}\left\Vert \Delta\vv_k\right\Vert^2+\left(\xi_B+\Upsilon\right)\left\Vert \Delta\vv_k\right\Vert^2-\left(\frac{\eta_{1,k}\xi_G^2}{2\Upsilon^2}-\frac{8\xi_G}{\eta_{2,k}\Upsilon^2}\right)\left\Vert\Delta\vv_k\right\Vert^2\\
& = \left(3\eta_{2,k}\Upsilon^2-\frac{\xi_B}{2}\right)\left\Vert \incxk\right\Vert^2+\left(\frac{2\Upsilon^2}{\xi_B}+\xi_B+\Upsilon+\frac{8\xi_G}{\eta_{2,k}\Upsilon^2}-\frac{\eta_{1,k}\xi_G^2}{2\Upsilon^2}\right)\left\Vert\Delta\vv_k\right\Vert^2.
\end{align*}
Combining the above inequality with \eqref{eq:descent2}, we get
\begin{multline*}
\begin{pmatrix}\nablaxL\\\nablalL\end{pmatrix}^T\incxlk
\le -\frac{\eta_{2,k}}{2}\normnablal^2+\left(3\eta_{2,k}\Upsilon^2-\frac{\xi_B}{2}\right)\left\Vert \incxk\right\Vert^2 \\ +\left(\frac{2\Upsilon^2}{\xi_B}+\xi_B+\Upsilon+\frac{8\xi_G}{\eta_{2,k}\Upsilon^2}-\frac{\eta_{1,k}\xi_G^2}{2\Upsilon^2}\right)\left\Vert\Delta\vv_k\right\Vert^2-\frac{\eta_{2,k}\xi_G}{8}\left\Vert \inclk\right\Vert^2.
\end{multline*}
In order to make the second term on the right-hand side negative, we let
\begin{equation}\label{eq:eta2cond}
\eta_{2,k}\le\frac{\xi_B}{12\Upsilon^2}.
\end{equation}
Without loss of generality, we assume $\Upsilon/2\ge 1\ge(\xi_B\vee\xi_G)$. Otherwise, we replace $\Upsilon$ by $\Upsilon\vee 2$, $\xi_{B}$ by $\xi_{B}\wedge1$, and $\xi_{G}$ by $\xi_{G}\wedge 1$. Then, we obtain
\begin{equation*}
\frac{2\Upsilon^2}{\xi_B}+\xi_B+\Upsilon+\frac{8\xi_G}{\eta_{2,k}\Upsilon^2}\le \frac{2\Upsilon^2}{\xi_B}+\frac{3\Upsilon}{2}+\frac{8\xi_G}{\eta_{2,k}\Upsilon^2}\le \frac{3\Upsilon^2}{\xi_B}+\frac{8\xi_G}{\eta_{2,k}\Upsilon^2}\stackrel{\mathclap{\eqref{eq:eta2cond}}}{\le}\frac{1}{4\eta_{2,k}}+\frac{8\xi_G}{\eta_{2,k}\Upsilon^2}\le\frac{1}{4\eta_{2,k}}+\frac{2}{\eta_{2,k}}\le\frac{2.5}{\eta_{2,k}}.
\end{equation*}
Using the above inequality and~\eqref{eq:eta2cond}, we have
\begin{equation}\label{eq:descent3}
\begin{pmatrix}\nablaxL\\\nablalL\end{pmatrix}^T\incxlk\le -\frac{\eta_{2,k}}{2}\normnablal^2-\frac{\xi_B}{4}\left\Vert \incxk\right\Vert^2+\left(\frac{2.5}{\eta_{2,k}}-\frac{\eta_{1,k}\xi_G^2}{2\Upsilon^2}\right)\left\Vert\Delta\vv_k\right\Vert^2-\frac{\eta_{2,k}\xi_G}{8}\left\Vert \inclk\right\Vert^2.
\end{equation}
In order to make the right-hand side negative, we let
\begin{equation}\label{eq:eta1cond}
\eta_{1,k}\ge \frac{5\Upsilon^2}{\eta_{2,k}\xi_G^2}.
\end{equation}
Notice that~\eqref{eq:eta2cond} and \eqref{eq:eta1cond} imply $\eta_{1,k}\ge\eta_{2,k}$ and $\eta_{1,k} \ge 16/(\eta_{2,k}\xi_G)$, hence, justify our presumptions. From \eqref{eq:descent3}, we get
\begin{align}\label{eq:descent4}
\begin{pmatrix}\nablaxL\\\nablalL\end{pmatrix}^T\incxlk &\stackrel{\mathclap{\eqref{eq:eta1cond}}}{\le}-\frac{\eta_{2,k}}{2}\normnablal^2-\frac{\eta_{2,k}\xi_G}{8}\normincl^2-\frac{\xi_B}{4}\normincx^2\nonumber\\
&\stackrel{\mathclap{\eqref{eq:eta2cond}}}{\le} -\frac{\eta_{2,k}}{2}\normnablal^2-\frac{\eta_{2,k}\xi_G}{8}\normincl^2-\frac{\eta_{2,k}\xi_G}{8}\normincx^2 \nonumber\\
& \le -\frac{\eta_{2,k}}{2}\normnablal^2-\frac{\eta_{2,k}\xi_G}{8}\|\Delta\vz_k\|^2.
\end{align}
Now we develop the second term in~\eqref{eq:descent1}. By Cauchy-Schwarz inequality, we get
\begin{align*}
\begin{pmatrix}\nablaxL\\\nablalL\end{pmatrix}^T\incedj
&\stackrel{\mathclap{\eqref{eq:opt_merit}}}{=}\incedj^T\begin{pmatrix} I+\eta_{2,k} H_k&\eta_{1,k} G_k^T \\ \eta_{2,k} G_k & I\end{pmatrix}\begin{pmatrix}\nablaxl\\c_k\end{pmatrix}\\
&\stackrel{\mathclap{\eqref{eq:kkt_sys_reduced}}}{=}-\incedj^T\begin{pmatrix}(I+\eta_{2,k}H_k)B_k+\eta_{1,k}G_k^TG_k&(I+\eta_{2,k}H_k)G_k^T \\ G_k(I+\eta_{2,k}B_k) & \eta_{2,k}G_kG_k^T\end{pmatrix}\incxlk\\
&\le\normedj\normexact\left((1+\eta_{2,k}\Upsilon)\Upsilon+(\eta_{1,k}+\eta_{2,k})\Upsilon^2+2(1+\eta_{2,k}\Upsilon)\Upsilon\right)\\
&=\normedj\normexact\left(3\Upsilon+4\eta_{2,k}\Upsilon^2+\eta_{1,k}\Upsilon^2\right)\\
&\stackrel{\mathclap{\eqref{eq:eta2cond}}}{\le}\normedj\normexact\left(3\Upsilon+\xi_B/3+\eta_{1,k}\Upsilon^2\right).
\end{align*}
Using $\Upsilon/2\ge 1\ge (\xi_B\vee \xi_G)$ and \eqref{eq:eta2cond}, we get $\eta_{1,k}\stackrel{\mathclap{\eqref{eq:eta1cond}}}{\ge}(5\Upsilon^2)/(\eta_{2,k}\xi_G^2)\ge 20/\eta_{2,k}\ge 40\cdot24\Upsilon$. Then, we have $19/6\le \eta_{1,k}\Upsilon$ and further obtain
\begin{equation*}
3\Upsilon+\frac{\xi_B}{3}+\eta_{1,k}\Upsilon^2\le\frac{19\Upsilon}{6}+\eta_{1,k}\Upsilon^2\le 2\eta_{1,k}\Upsilon^2.
\end{equation*}
Using the above inequality, we get
\begin{equation}\label{eq:descent5}
\begin{pmatrix}\nablaxL\\\nablalL\end{pmatrix}^T\incedj\le2\eta_{1,k}\Upsilon^2\normedj\normexact.
\end{equation}
Finally, we obtain
\begin{align*}
\begin{pmatrix}\nablaxL\\\nablalL\end{pmatrix}^T\incxltkj
&\stackrel{\mathclap{\eqref{eq:descent1}}}{=}\begin{pmatrix}\nablaxL\\\nablalL\end{pmatrix}^T\incxlk + \begin{pmatrix}\nablaxL\\\nablalL\end{pmatrix}^T\incedj\\
&\stackrel{\mathclap{\substack{\eqref{eq:descent4}\\\eqref{eq:descent5}}}}{\le}-\frac{\eta_{2,k}}{2}\normnablal^2-\frac{\eta_{2,k}\xi_G}{8}\normexact^2+2\eta_{1,k}\Upsilon^2\normedj\normexact\\
&\stackrel{\mathclap{\eqref{eq:errorstep}}}{\le}-\frac{\eta_{2,k}}{2}\normnablal^2-\frac{\eta_{2,k}\xi_G}{8}\normexact^2+2\delta_k\eta_{1,k}\Upsilon^2\normexact^2\\
&\le-\frac{\eta_{2,k}}{2}\normnablal^2-\left(\frac{\eta_{2,k}\xi_G}{8}-2\delta_k\eta_{1,k}\Upsilon^2\right)\|\Delta\vz_k\|^2.
\end{align*}
In order to make the right-hand side negative, we let 
\begin{equation*}
\delta_k\le\frac{\eta_{2,k}\xi_G}{16\eta_{1,k}\Upsilon^2}.
\end{equation*}
Thus, the descent direction condition is satisfied as long as $\eta_{2,k}\le \xi_B/12\Upsilon^2$, $\eta_{1,k}\eta_{2,k}\ge 5\Upsilon^2/\xi_G^2$, and $\delta_k\eta_{1,k}/\eta_{2,k}\le\xi_G/(16\Upsilon^2)$. Finally, we let $\Upsilon\leftarrow 12\Upsilon^2/\xi_B\vee 5\Upsilon^2/\xi_G^2\vee 16\Upsilon^2/\xi_G$, and complete the proof.

\subsection{Proof of Lemma \ref{lemma12}}

We denote the event that the algorithm reaches $\vz_k$ as $\mZ_k$. Then the event $\cap_{k=0}^\infty \mathcal{Z}_k$ implies that the algorithm generates the iterates infinitely. From Lemmas \ref{lemma10} and \ref{lemma11}, and noting that the Armijo condition \eqref{eq:armijo} can be satisfied for small enough $\alpha_k$ as long as $\tilde{\Delta}\vz_k$ is a descent direction (as implied by Lemma \ref{lemma11}), we have $P(\mZ_{k+1}|\mZ_k)=1$, $\forall k\geq 0$. Thus,
\begin{equation}\label{eq:event}
P\left(\cap_{k=0}^\infty \mathcal{Z}_k\right)\ge 1-\sum_{k=0}^\infty P(\mZ_k^c)= 1 - \sum_{k=0}^\infty\int P(\mZ_{k+1}^c|\mZ_{k})P(\mZ_{k})d\vz_k=1.
\end{equation}
We now start from finding the lower bound of $\delta_k^{\text{trial}}\eta_{1,k}/\eta_{2,k}$. Since the updating rule of the parameters \eqref{eq:parameter_update} increases $\eta_{1,k}$ by a factor of $\nu^2$ and decreases $\eta_{2,k}$ by a factor of $\nu$, we know that $\eta_{1,0}\le\eta_{1,k}$ and $\eta_{2,0}\ge\eta_{2,k}$ for all $k\ge0$. By Assumptions \ref{assum1} and \ref{assum2}, we know there exist constants $\Upsilon_{\upsilon},\Psi>0$ such that $\Upsilon_k\le\Upsilon_{\upsilon}$ and $\Psi_k\le\Psi$, $\forall k\geq 0$. Thus, for $k\ge 0$, we have 
\begin{equation*}
(1+\eta_{1,k}+\eta_{2,k})\Upsilon_k^2\Psi_k^2\le \frac{\eta_{1,k}}{\eta_{1,0}}(\Upsilon_{\upsilon}^2\Psi^2+\eta_{2,0}\Upsilon_{\upsilon}^2\Psi^2)+\eta_{1,k}\Upsilon_{\upsilon}^2\Psi^2=\frac{\eta_{1,k}}{\eta_{1,0}}(\Upsilon_{\upsilon}^2\Psi^2+\eta_{1,0}\Upsilon_{\upsilon}^2\Psi^2+\eta_{2,0}\Upsilon_{\upsilon}^2\Psi^2).
\end{equation*}
Using the above inequality, we have
\begin{equation}\label{eq:delta_trial_bound} 
\dfrac{\delta_k^{\text{trial}}\eta_{1,k}}{\eta_{2,k}} \stackrel{\eqref{eq:delta_trial}}{=}\dfrac{(0.5-\beta)\eta_{1,k}}{(1+\eta_{1,k}+\eta_{2,k})\Upsilon_k^2\Psi_k^2}\ge\dfrac{(0.5-\beta)\eta_{1,0}}{\Upsilon_{\upsilon}^2\Psi^2+\eta_{1,0}\Upsilon_{\upsilon}^2\Psi^2+\eta_{2,0}\Upsilon_{\upsilon}^2\Psi^2}.
\end{equation}
By Lemma~\ref{lemma11} and~\eqref{eq:delta_trial_bound}, we have
\begin{equation*}
\eta_{1,k}\eta_{2,k}\ge \Upsilon,\quad\eta_{2,k}\le\frac{1}{\Upsilon},\quad\frac{\delta_k\eta_{1,k}}{\eta_{2,k}}\le\frac{1}{\Upsilon}\wedge\dfrac{(0.5-\beta)\eta_{1,0}}{\Upsilon_{\upsilon}^2\Psi^2+\eta_{2,0}\Upsilon_{\upsilon}^2\Psi^2+\eta_{1,0}\Upsilon_{\upsilon}^2\Psi^2}.
\end{equation*}
Notice that the lower bound of $\eta_{1,k}\eta_{2,k}$ and the upper bounds of $\eta_{2,k}$, $\delta_k\eta_{1,k}/\eta_{2,k}$ do not depend on $k$. Since the updating rule of the parameters \eqref{eq:parameter_update} implies that $\eta_{1,k}\eta_{2,k}$ increases by a factor of $\nu$, $\eta_{2,k}$ decreases by a factor of $\nu$, and $\delta_k\eta_{1,k}/\eta_{2,k}$ decreases by at least a factor of $\nu$, there exists an iteration threshold $K$ such that $(\eta_{1,k},\eta_{2,k},\delta_k)=(\eta_{1,K},\eta_{2,K},\delta_K)$ for all $k$. Using \eqref{eq:event}, we complete the proof.

\subsection{Proof of Lemma~\ref{lemma13}}

By Assumption \ref{assum1}, we know that $\nabla\mathcal{L}_{\veta}$ is Lipschitz continuous. Let us denote the Lipschitz constant of $\nabla\mathcal{L}_{\veta_k}$ by $\Upsilon_{\eta_k}$. We have for any $k\ge 0$,
\begin{align*}
\Lagtrial
&= \Lag+\alpha_k(\nablaL)^T\incztk+\alpha_k\int_0^1\{\nabla\mathcal{L}^T_{\veta_k}(\vz_k + t\alpha_k\incztk)-\nabla\mathcal{L}^k_{\veta_k}\}^T\incztk dt\\
&\le \Lag+\alpha_k(\nablaL)^T\incztk+\Upsilon_{\eta_k}\alpha_k^2\|\incztk\|^2 \int_0^1 t dt\\
&= \Lag+\alpha_k(\nablaL)^T\incztk+\Upsilon_{\eta_k}\alpha_k^2 \|\incztk\|^2/2.
\end{align*}
By Assumption \ref{assum2}, we know there exists a constant $\Psi>0$ such that $\Psi_k\le\Psi$ for all $k\ge 0$. Using \eqref{eq:errorstep}, we have
\begin{equation}\label{eq:errorstep2}
\|\incztk\|\le\delta_k\|\inczk\|+\|\inczk\|=\left(\delta_k+1\right)\|\inczk\|\leq2\|\inczk\|.
\end{equation}
By Lemma \ref{lemma12}, we know that $\eta_{2,k}\ge\eta_{2,K}$ for all $k\ge0$. We let $\Upsilon_{\eta}\coloneqq\Upsilon_{\eta_0}\vee\cdots\vee\Upsilon_{\eta_K}$. Then, we have
\begin{align*}
\Lagtrial
&\le \Lag+\alpha_k(\nablaL)^T\incztk+\alpha_k^2\Upsilon_\eta\|\incztk\|^2/2\\
& \stackrel{\mathclap{\eqref{eq:errorstep2}}}{\leq} \Lag+\alpha_k(\nablaL)^T\incztk+2\alpha_k^2\Upsilon_\eta \|\inczk\|^2\\
&\stackrel{\mathclap{\substack{\eqref{eq:kkt_sys_reduced}\\\eqref{eq:Psi}}}}{\le} \Lag+\alpha_k(\nablaL)^T\incztk+2\alpha_k^2\Upsilon_\eta\Psi^2\normnablal^2\\
&\stackrel{\mathclap{\eqref{eq:cond2}}}{\leq} \Lag+\alpha_k(\nablaL)^T\incztk-(4\alpha_k^2\Upsilon_\eta\Psi^2/\eta_{2,k})(\nablaL)^T\incztk\\
&\le \Lag+\alpha_k(1-4\alpha_k\Upsilon_\eta\Psi^2/\eta_{2,K})(\nablaL)^T\incztk.
\end{align*}
Now, we let
\begin{equation*}
1-\frac{4\alpha_k\Upsilon_\eta\Psi^2}{\eta_{2,K}}\ge\beta\Leftrightarrow \alpha_k\le\frac{(1-\beta)\eta_{2,K}}{4\Upsilon_\eta\Psi^2}.
\end{equation*}
Since the upper bound of $\alpha_k$ does not depend on $k$, there exists $\alpha_{\min}>0$ such that for any $k$, $\alpha_{\min}\le\alpha_k$, when we do, for example, backtracking. Using \eqref{eq:event}, we complete the proof.

\subsection{Proof of Theorem~\ref{theorem1}}

By Lemmas \ref{lemma11}, \ref{lemma12}, and \ref{lemma13}, we have for any $k\ge K$,
\begin{equation*}
\mathcal{L}_{\veta_K}^{k+1}-\mathcal{L}_{\veta_K}^{k} \leq \alpha_k\beta(\nabla\mL_{\veta_K}^k)^T\tilde{\Delta}\vz_k \stackrel{\eqref{eq:cond2}}{\leq} -\eta_{2,k}\alpha_k\beta\|\nabla\mL_k\|^2/2 \le-\eta_{2,K}\alpha_{\min}\beta\normnablal^2/2.
\end{equation*}
Summing over $k\ge K$, we have
\begin{equation*}
\sum_{k=K}^{\infty}\normnablal^2\le\frac{2}{\eta_{2,K}\alpha_{\min}\beta}\left(\mathcal{L}_{\veta_K}^{K}-\underset{\mX\times \Lambda}{\min}\mathcal{L}(\vx,\vlambda)\right)<\infty.
\end{equation*}
Therefore, $\normnablal \rightarrow 0\text{ as }k\rightarrow\infty$. Using~\eqref{eq:event}, we have $P\left(\normnablal \rightarrow 0\text{ as }k\rightarrow\infty\right)=1$. This completes the proof.

\section{Proofs of Section~\ref{sec:4}}

\subsection{Proof of Theorem \ref{theorem2}}\label{subsec:theorem2}

We first show that for all sufficiently large $k$, a unit stepsize is admissible. It suffices to show that for all sufficiently large $k$,
\begin{equation*}
\mathcal{L}_{\veta_k}(\vz_k+\incztk)\le \Lag+\beta(\nablaL)^T\incztk.
\end{equation*}	
Let $\vu\in\mR^l$, $M:\mR^l\rightarrow \mR^{p\times q}$ and $w:\mR^l\rightarrow\mR^{q}$. Let us denote the columns of $M(\vu)$ as $m_i(\vu)\in\mR^p$ for $i=1,\dots,q$ and the entries of $w(\vu)$ as $w_i(\vu)$ for $i=1,\dots,q$. We can write $\nabla_{\vu} (M(\vu)w(\vu))\in \mR^{l\times p}$ as
\begin{equation*}
\nabla_{\vu} (M(\vu)w(\vu)) = \sum_{i=1}^qw_i(\vu)\nabla_{\vu}m_i(\vu) + \nabla_{\vu} w(\vu)M(\vu)^T \eqqcolon \langle w(\vu), \nabla_{\vu}M(\vu)\rangle + \nabla_{\vu} w(\vu)M(\vu)^T.
\end{equation*}
With the above definition and using \eqref{eq:opt_merit}, we can compute
\begin{align*}
\nabla_{\vx}^2\mL_{\veta_k}^k&=\nabla_\vx\left(\nablaxl+\eta_{2,k}H_k\nablaxl+\eta_{1,k}G_k^Tc_k \right)=H_k+\eta_{2,k}\left(\langle\nablaxl, \nabla_\vx H_k\rangle+H_k^2\right)+\eta_{1,k}\left(\langle c_k,\nabla G_k\rangle+G_k^TG_k\right),\\
\nabla_{\vlambda}^2\mL_{\veta_k}^k&=\nabla_{\vlambda}\left(c_k+\eta_{2,k}G_k\nablaxl\right)=\eta_{2,k}G_kG_k^T,\\
\nabla_{\vlambda\vx}^2\mL_{\veta_k}^k&=\nabla_{\vlambda}\left(\nablaxl+\eta_{2,k}H_k\nablaxl+\eta_{1,k}G_k^Tc^k\right)=G_k+\eta_{2,k}\left(\langle\nablaxl,\nabla_{\vlambda}H_k\rangle+G_kH_k\right).
\end{align*}
By Assumption~\ref{assum4}, we know $\nabla^2\mathcal{L}_{\veta}$ is continuous over $\mathcal{X}$. We define
\begin{equation*}
\mathcal{H}_k=\begin{pmatrix}H_k+\eta_{2,k}H_k^2+\eta_{1,k}G_k^TG_k & G_k^T+\eta_{2,k}H_kG_k^T\\ G_k+\eta_{2,k}G_kH_k & \eta_{2,k}G_kG_k^T\end{pmatrix}.
\end{equation*}
Using $\normnablal=\Vert(\nablaxl,c_k)\Vert=o(1)$, we have $\nablaLL=\mathcal{H}_k+o(1)$. We now apply Taylor's theorem and obtain
\begin{align*}
&\mathcal{L}_{\veta_k}(\vx_k+\incxtk,\vlambda_k+\incltk)\\
&\quad\le \Lag+\begin{pmatrix}\nablaxL \\ \nablalL\end{pmatrix}^T\incxltk+\frac{1}{2}\incxltk^T\nablaLL\incxltk+o\left(\left\Vert \incxltk\right\Vert^2\right)\\
&\quad= \Lag+\begin{pmatrix}\nablaxL \\ \nablalL\end{pmatrix}^T\incxltk+\frac{1}{2}\incxltk^T\mathcal{H}_k\incxltk+o\left(\left\Vert \incxltk\right\Vert^2\right)\\
&\quad\stackrel{\mathclap{\substack{\eqref{eq:opt_merit}\\\eqref{eq:kkt_sys_reduced}}}}{=} \Lag+\frac{1}{2}\begin{pmatrix}\nablaxL \\ \nablalL\end{pmatrix}^T\incxltk+\frac{1}{2}\incxltk^T\mathcal{H}_k\incxltk\\
&\quad\quad-\frac{1}{2}\incxltk^T\begin{pmatrix}(I+\eta_{2,k}H_k)B_k+\eta_{1,k}G_k^TG_k&(I+\eta_{2,k}H_k)G_k^T\\G_k(I+\eta_{2,k}B_k)&\eta_{2,k}G_kG_k^T\end{pmatrix}\incxlk+o\left(\left\Vert \incxltk\right\Vert^2\right)\\
&\quad= \Lag+\frac{1}{2}\begin{pmatrix}\nablaxL \\ \nablalL\end{pmatrix}^T\incxltk+\frac{1}{2}\incxltk^T\mathcal{H}_k\inced\\
&\quad\quad+\frac{1}{2}\incxltk^T\left(\mathcal{H}_k-\begin{pmatrix}(I+\eta_{2,k}H_k)B_k+\eta_{1,k}G_k^TG_k&(I+\eta_{2,k}H_k)G_k^T\\G_k(I+\eta_{2,k}B_k)&\eta_{2,k}G^kG_k^T\end{pmatrix}\right)\incxlk+o\left(\left\Vert \incxltk\right\Vert^2\right)\\
&\quad= \Lag+\frac{1}{2}\begin{pmatrix}\nablaxL \\ \nablalL\end{pmatrix}^T\incxltk+\frac{1}{2}\incxltk^T\mathcal{H}_k\inced\\
&\quad\quad+\frac{1}{2}\incxltk^T\begin{pmatrix}(I+\eta_{2,k}H_k)(H_k-B_k)&\vzero\\\eta_{2,k}G_k(H_k-B_k)&\vzero\end{pmatrix}\incxlk+o\left(\left\Vert \incxltk\right\Vert^2\right).
\end{align*}
By Assumption \ref{assum5}, we have $\left\Vert(H_k-B_k)\incxk\right\Vert\le \left\Vert(H_k-B_k)\right\Vert\left\Vert \incxk\right\Vert=o\left(\left\Vert \incxk\right\Vert\right)$. Thus, we have for any $k\ge0$, 
\begin{align}\label{eq:unit_stepsize1}
&\mathcal{L}_{\veta_k}(\vx_k+\incxtk,\vlambda_k+\incltk)\nonumber\\
&\quad\le \Lag+\frac{1}{2}\begin{pmatrix}\nablaxL\\\nablalL\end{pmatrix}^T\incxltk+\frac{1}{2}\left\Vert\incxltk\right\Vert\left\Vert\mathcal{H}_k\right\Vert\left\Vert\inced\right\Vert\nonumber\\
&\quad\quad+\frac{1}{2}\norminexact\left(\left\Vert I+\eta_{2,k}H_k\right\Vert\left\Vert(H_k-B_k)\incxk\right\Vert+\left\Vert\eta_{2,k} G_k\right\Vert\left\Vert(H_k-B_k)\incxk\right\Vert\right)+o\left(\left\Vert \incxltk\right\Vert^2\right)\nonumber\\
&\quad\stackrel{\mathclap{\substack{\eqref{eq:errorstep}\\ \eqref{eq:errorstep2}}}}{\le} \Lag+\frac{1}{2}\begin{pmatrix}\nablaxL \\ \nablalL\end{pmatrix}^T\incxltk+\frac{1}{2}\left\Vert\mathcal{H}_k\right\Vert\cdot 2\normexact\cdot\delta_k\normexact+o\left(\norminexact^2\right)\nonumber\\
&\quad\stackrel{\eqref{eq:aux}}{\le} \Lag+\frac{1}{2}\begin{pmatrix}\nablaxL \\ \nablalL\end{pmatrix}^T\incxltk+\delta_k(3\Upsilon_k+\eta_{1,k}\Upsilon_k^2+4\eta_{2,k}\Upsilon_k^2)\normexact^2+o\left(\normexact^2\right)\nonumber\\
&\quad\stackrel{\eqref{eq:aux}}{\le} \Lag+\frac{1}{2}\begin{pmatrix}\nablaxL \\ \nablalL\end{pmatrix}^T\incxltk+4\delta_k\Upsilon_k^2(1+\eta_{1,k}+\eta_{2,k})\normexact^2+o\left(\normexact^2\right).
\end{align}
Using the fact that $\delta_k\le\delta_k^{\text{trial}}=\dfrac{(0.5-\beta)\eta_{2,k}}{(1+\eta_{1,k}+\eta_{2,k})\Upsilon_k^2\Psi_k^2}$, we have for any $k\ge0$,
\begin{align*}
& 4.05\delta_k\Upsilon_k^2(1+\eta_{1,k}+\eta_{2,k})\|\Delta\vz_k\|^2 \le \left(0.5-\beta\right)\frac{4.05\eta_{2,k}}{\Psi_k^2}\|\Delta\vz_k\|^2 \stackrel{\eqref{eq:kkt_sys_reduced}}{\leq} \left(0.5-\beta\right)\frac{\eta_{2,k}}{2}\frac{8.1}{\Psi_k^2}\|\KKT^{-1}\|^2\normnablal^2\\
&\stackrel{\mathclap{\substack{\eqref{eq:Psi}}}}{\le} \left(0.5-\beta\right)\eta_{2,k}\normnablal^2/2 \stackrel{\eqref{eq:cond2}}{\leq} -\left(0.5-\beta\right) \begin{pmatrix}\nablaxL \\ \nablalL\end{pmatrix}^T\incxltk.
\end{align*}
Thus, for any $k\ge0$, we have
\begin{equation*}\frac{1}{2}\begin{pmatrix}\nablaxL\\\nablalL\end{pmatrix}^T\incxltk+4.05\delta_k\Upsilon_k^2(1+\eta_{1,k}+\eta_{2,k})\normexact^2\le\beta\begin{pmatrix}\nablaxL \\ \nablalL\end{pmatrix}^T\incxltk.    
\end{equation*}
Thus, there exists an iteration threshold $K_1$ such that for any $k\ge K_1$, we have
\begin{equation*}
\frac{1}{2}\begin{pmatrix}\nablaxL\\\nablalL\end{pmatrix}^T\incxltk+4\delta_k\Upsilon_k^2(1+\eta_{1,k}+\eta_{2,k})\normexact^2+o\left(\normexact^2\right)\le\beta\begin{pmatrix}\nablaxL \\ \nablalL\end{pmatrix}^T\incxltk.    
\end{equation*}
Plugging the above inequality back into \eqref{eq:unit_stepsize1}, we know for any $k\ge K_1$, 
\begin{equation*}
\mathcal{L}_{\veta_k}(\vx_k+\incxtk,\vlambda_k+\incltk)\le \Lag+\beta\begin{pmatrix}\nablaxL \\ \nablalL\end{pmatrix}^T\incxltk\Leftrightarrow\mathcal{L}_{\veta_k}(\vz_k+\incztk)\le \Lag+\beta(\nablaL)^T\incztk.
\end{equation*}
This completes the first part of the proof. Next, we show for all sufficiently large $k$,
\begin{equation*}
\|\vz_k+\incztk-\vz^\star\|\le(1+\varphi)\theta\delta_K\|\vz_k-\vz^\star\|,\quad \text{ for any } \varphi>0.
\end{equation*}
We start from dividing $\vz_k+\incztk-\vz^\star$ into two terms as follows:
\begin{equation}\label{eq:local2-1}
\vz_k+\incztk-\vz^\star=\begin{pmatrix} \vx_k+\incxtk -\vx^\star\\ \vlambda_k+\incltk -\vlambda^\star\end{pmatrix} =\begin{pmatrix}\vx_k+\incxk-\vx^\star\\\vlambda_k+\inclk-\vlambda^\star\end{pmatrix}+\begin{pmatrix}\incxtk-\incxk\\\incltk-\inclk\end{pmatrix}.
\end{equation}
For the first term, we apply Assumption \ref{assum2} ($\Gamma_k$ is invertible) and have for any $k\ge0$,
\begin{align}\label{eq:local2-2}
\begin{pmatrix}\vx_k+\incxk-\vx^\star\\\vlambda_k+\inclk-\vlambda^\star\end{pmatrix}
&=\kktinv\kkt\begin{pmatrix} \vx_k-\vx^\star\nonumber\\ \vlambda_k-\vlambda^\star\end{pmatrix}+\incxlk\nonumber\\
&\stackrel{\mathclap{\eqref{eq:kkt_sys_reduced}}}{=}\kktinv\kkt\begin{pmatrix} \vx_k-\vx^\star\\ \vlambda_k-\vlambda^\star\end{pmatrix}-\kktinv\nablal\nonumber\\
&=\kktinv\left(\kkt\begin{pmatrix} \vx_k-\vx^\star\\ \vlambda_k-\vlambda^\star\end{pmatrix}-\nablal\right)\nonumber\\
&\stackrel{\eqref{eq:kkt}}{=}\kktinv\left(\kkt\begin{pmatrix} \vx_k-\vx^\star\\ \vlambda_k-\vlambda^\star\end{pmatrix}-\left(\nablal-\nabla\mathcal{L}_{\star}\right)\right).
\end{align}
By Assumption~\ref{assum1}, we know $\nabla^2\mathcal{L}$ is continuous over $\mathcal{X}$. We apply Taylor's theorem and obtain
\begin{align*}
\nablal-\nabla\mathcal{L}_{\star}
&=\int_0^1 \nabla^2\mathcal{L}\left(\vx_k+t(\vx^\star-\vx_k),\vlambda_k+t(\vlambda^\star-\vlambda_k)\right)\begin{pmatrix} \vx_k-\vx^\star\\ \vlambda_k-\vlambda^\star\end{pmatrix}dt\\
&=\int_0^1 \begin{pmatrix}H\left(\vx_k+t(\vx^\star-\vx_k),\vlambda_k+t(\vlambda^\star-\vlambda_k)\right)& G^T(\vx_k+t(\vx^\star-\vx_k))\\ G(\vx_k+t(\vx^\star-\vx_k))&\vzero\end{pmatrix}\begin{pmatrix}\vx_k-\vx^\star\\\vlambda_k-\vlambda^\star\end{pmatrix}dt\\
& \coloneqq \int_0^1 \begin{pmatrix}H_k(t)& G_k^T(t)\\ G_k(t)&\vzero\end{pmatrix}\begin{pmatrix}\vx_k-\vx^\star\\\vlambda_k-\vlambda^\star\end{pmatrix}dt.
\end{align*}
Plugging the above equation back into~\eqref{eq:local2-2}, we have
\begin{equation*}
\begin{pmatrix}\vx_k+\incxk-\vx^\star\\\vlambda_k+\inclk-\vlambda^\star\end{pmatrix}
=\kktinv\left(\int_0^1 \begin{pmatrix}B_k-H_k(t)&G_k^T-G_k(t)^T\\G_k-G_k(t)&\vzero\end{pmatrix}\begin{pmatrix}\vx_k-\vx^\star\\\vlambda_k-\vlambda^\star\end{pmatrix}dt\right).
\end{equation*}
By Assumption~\ref{assum1}, we know $H$ and $G$ are Lipschitz continuous over $\mathcal{X}$. Thus, we apply Assumption~\ref{assum5} and have
\begin{align}
\left\Vert\begin{pmatrix}\vx_k+\incxk-\vx^\star\\\vlambda_k+\inclk-\vlambda^\star\end{pmatrix}\right\Vert
&\le\left\Vert\kktinv\right\Vert\left\Vert\int_0^1 \begin{pmatrix}B_k-H_k(t)& G_k^T-G_k(t)^T\\ G_k-G_k(t)&\vzero\end{pmatrix}\begin{pmatrix}\vx_k-\vx^\star\\\vlambda_k-\vlambda^\star\end{pmatrix}dt\right\Vert\nonumber\\
&\le\left\Vert\kktinv\right\Vert\int_0^1\left\Vert\begin{pmatrix} B_k-H_k(t)& G_k^T-G_k(t)^T\\G_k-G_k(t)&\vzero\end{pmatrix}\right\Vert\left\Vert\begin{pmatrix}\vx_k-\vx^\star\\\vlambda_k-\vlambda^\star\end{pmatrix}\right\Vert dt\nonumber\\
&\le O\left(\tau_k\left\Vert\begin{pmatrix}\vx_k-\vx^\star\\\vlambda_k-\vlambda^\star\end{pmatrix}\right\Vert\right)+O\left(\left\Vert\begin{pmatrix}\vx_k-\vx^\star\\\vlambda_k-\vlambda^\star\end{pmatrix}\right\Vert^2\right).\label{eq:local2-4}
\end{align}
Taking $\ell_2$ norm on both sides of~\eqref{eq:local2-1}, we have
\begin{align*}
\left\Vert\begin{pmatrix} \vx_k+\incxtk -\vx^\star\\ \vlambda_k+\incltk -\vlambda^\star\end{pmatrix} \right\Vert
&\le\normed+\left\Vert\begin{pmatrix} \vx_k+\incxk -\vx^\star\\\vlambda_k+\inclk -\vlambda^\star\end{pmatrix} \right\Vert \stackrel{\mathclap{\eqref{eq:errorsteptheta}}}{\le}\theta\delta_k\normexact+\left\Vert\begin{pmatrix} \vx_k+\incxk -\vx^\star\\\vlambda_k+\inclk -\vlambda^\star\end{pmatrix} \right\Vert\\
&\le\theta\delta_k\left(\left\Vert\begin{pmatrix}\vx_k-\vx^\star\\\vlambda_k-\vlambda^\star\end{pmatrix}\right\Vert+\left\Vert\begin{pmatrix} \vx_k+\incxk -\vx^\star\\\vlambda_k+\inclk -\vlambda^\star\end{pmatrix} \right\Vert\right)+\left\Vert\begin{pmatrix} \vx_k+\incxk -\vx^\star\\\vlambda_k+\inclk -\vlambda^\star\end{pmatrix} \right\Vert\\
&\stackrel{\mathclap{\eqref{eq:local2-4}}}{\le}\theta\delta_k\left\Vert\begin{pmatrix}\vx_k-\vx^\star\\\vlambda_k-\vlambda^\star\end{pmatrix}\right\Vert+o\left(\left\Vert\begin{pmatrix}\vx_k-\vx^\star\\\vlambda_k-\vlambda^\star\end{pmatrix}\right\Vert\right).
\end{align*}
This completes the proof.

\subsection{Proof of Corollary~\ref{collorary1}}

By the proof of Theorem \ref{theorem2}, we know that $\alpha_k=1$ for all $k\ge K_1$ regardless of the value of $\theta_k$. Taking $\ell_2$ norm on both sides of~\eqref{eq:local2-1}, we have for any $k\ge0$,
\begin{align*}
\left\Vert\begin{pmatrix} \vx_k+\incxtk -\vx^\star\\ \vlambda_k+\incltk -\vlambda^\star\end{pmatrix} \right\Vert
&\le\normed+\left\Vert\begin{pmatrix} \vx_k+\incxk -\vx^\star\\\vlambda_k+\inclk -\vlambda^\star\end{pmatrix} \right\Vert\\
&\stackrel{\eqref{eq:errorsteptheta}}{\le}\theta_k\delta_k\normexact+\left\Vert\begin{pmatrix} \vx_k+\incxk -\vx^\star\\\vlambda_k+\inclk -\vlambda^\star\end{pmatrix} \right\Vert\\
&\le\theta_k\delta_k\left(\left\Vert\begin{pmatrix}\vx_k-\vx^\star\\\vlambda_k-\vlambda^\star\end{pmatrix}\right\Vert+\left\Vert\begin{pmatrix} \vx_k+\incxk -\vx^\star\\\vlambda_k+\inclk -\vlambda^\star\end{pmatrix} \right\Vert\right)+\left\Vert\begin{pmatrix} \vx_k+\incxk -\vx^\star\\\vlambda_k+\inclk -\vlambda^\star\end{pmatrix} \right\Vert\\
&\stackrel{\eqref{eq:local2-4}}{\le} \theta_k\delta_k\left\Vert\begin{pmatrix}\vx_k-\vx^\star\\\vlambda_k-\vlambda^\star\end{pmatrix}\right\Vert+O\left(\tau_k\left\Vert\begin{pmatrix}\vx_k-\vx^\star\\\vlambda_k-\vlambda^\star\end{pmatrix}\right\Vert\right)+O\left(\left\Vert\begin{pmatrix}\vx_k-\vx^\star\\\vlambda_k-\vlambda^\star\end{pmatrix}\right\Vert^2\right)\\
&= O(\theta_k\delta_k+\tau_k)\left\Vert\begin{pmatrix}\vx_k-\vx^\star\\\vlambda_k-\vlambda^\star\end{pmatrix}\right\Vert+O\left(\left\Vert\begin{pmatrix}\vx_k-\vx^\star\\\vlambda_k-\vlambda^\star\end{pmatrix}\right\Vert^2\right).
\end{align*}
This completes the proof.

\section{Additional Algorithms, Tables, and Figures}\label{sec:appenB}

Algorithm \ref{alg:alg2} and Algorithm \ref{alg:alg3} use the $\ell_1$ penalized merit function of the form $\phi_{\pi}(\vx)=f(\vx)+\pi\Vert c(\vx)\Vert_1$. \textit{Termination Test 1}, \textit{Termination Test 2}, \textit{Model Reduction Condition}, and $\pi_k^{\text{trial}}$ are referred to in \citet{Byrd2008Inexact}. Furthermore, we set $\sigma=\tau(1-\epsilon)$, $\beta=\kappa_2=\|\nabla\mathcal{L}_0\|_1/(\|c_0\|_1+1)\vee 1$ as in \citet{Byrd2008Inexact}. Algorithm~\ref{alg:alg4} uses the augmented Lagrangian function of the form $\mathcal{L}_{\mu}(\vx,\vlambda)=\mathcal{L}(\vx,\vlambda)+(\mu/2)\|c(\vx)\|^2$.
 
\begin{figure*}[h]
\begin{center}
\begin{minipage}[ht]{.49\linewidth}
\centerline{\includegraphics[scale=0.35]{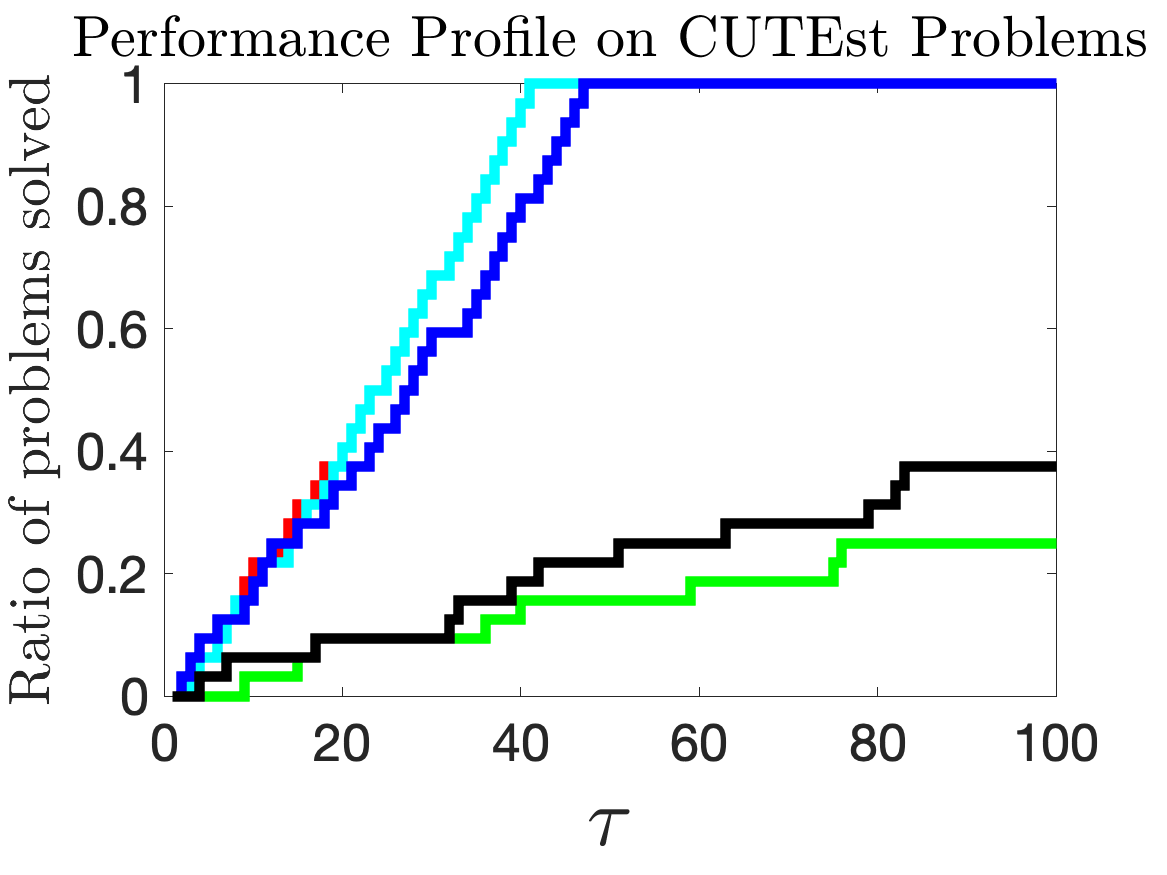}}
\end{minipage}
\begin{minipage}[ht]{.49\linewidth}
\centerline{\includegraphics[scale=0.35]{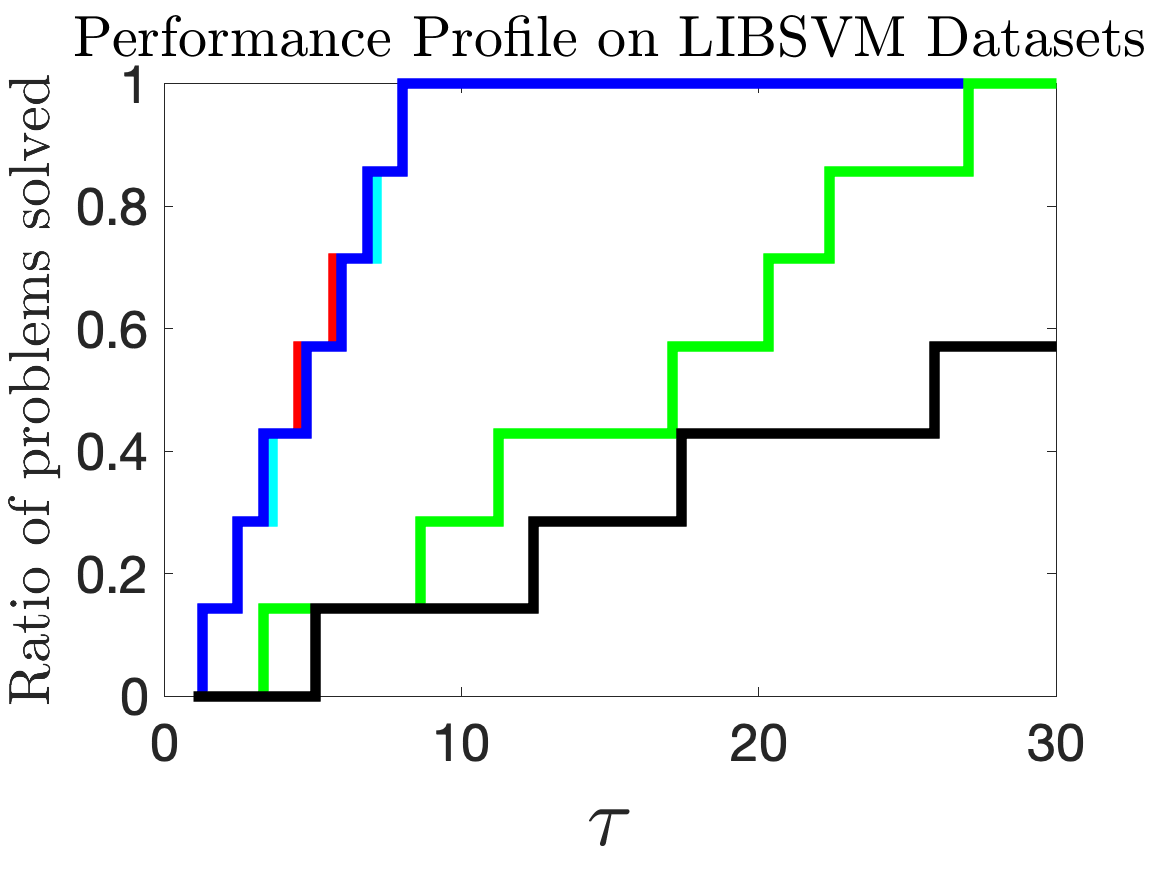}}
\end{minipage}\\
\begin{minipage}[ht]{.24\linewidth}
\centerline{\includegraphics[scale=0.60]{Figure/Legend.png}}
\end{minipage}
\caption{The performance profiles of the total number of flops for $\myalg$, Algorithms \ref{alg:alg2}, Algorithm \ref{alg:alg3}, and Augmented Lagrangian on CUTEst problems (Left) and on 7 LIBSVM datasets (right); the ratio of the problems  solved on the $y$ axis, while the proportion of the total number of flops (called the performance ratio and denoted as $\tau$) on the $x$ axis. See \citet{Dolan2002Benchmarking} for more details.}\label{fig:fig3}
	\end{center}
\end{figure*}

\begin{figure*}[t]
	\begin{center}
		\begin{minipage}[ht]{.33\linewidth}
			\centerline{\includegraphics[scale=0.3]{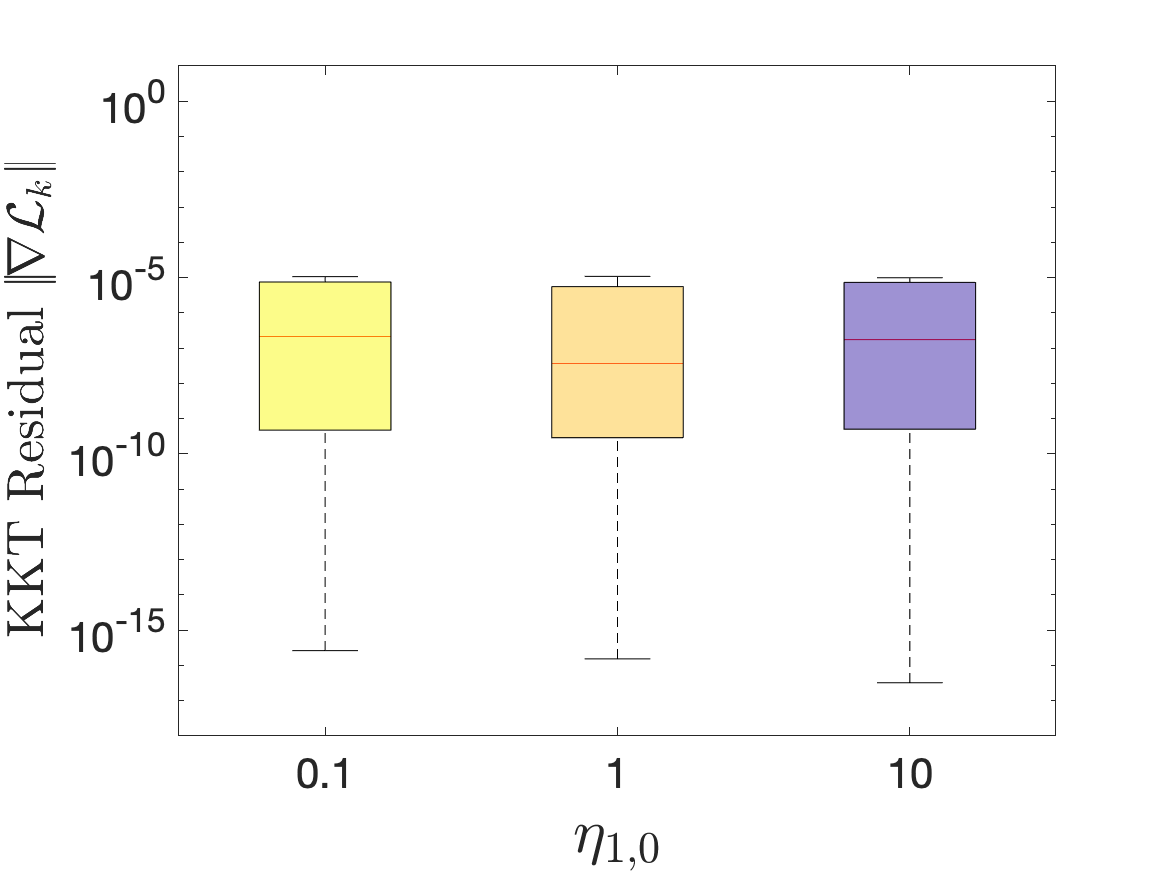}}
		\end{minipage}%
		\begin{minipage}[ht]{.33\linewidth}
			\centerline{\includegraphics[scale=0.3]{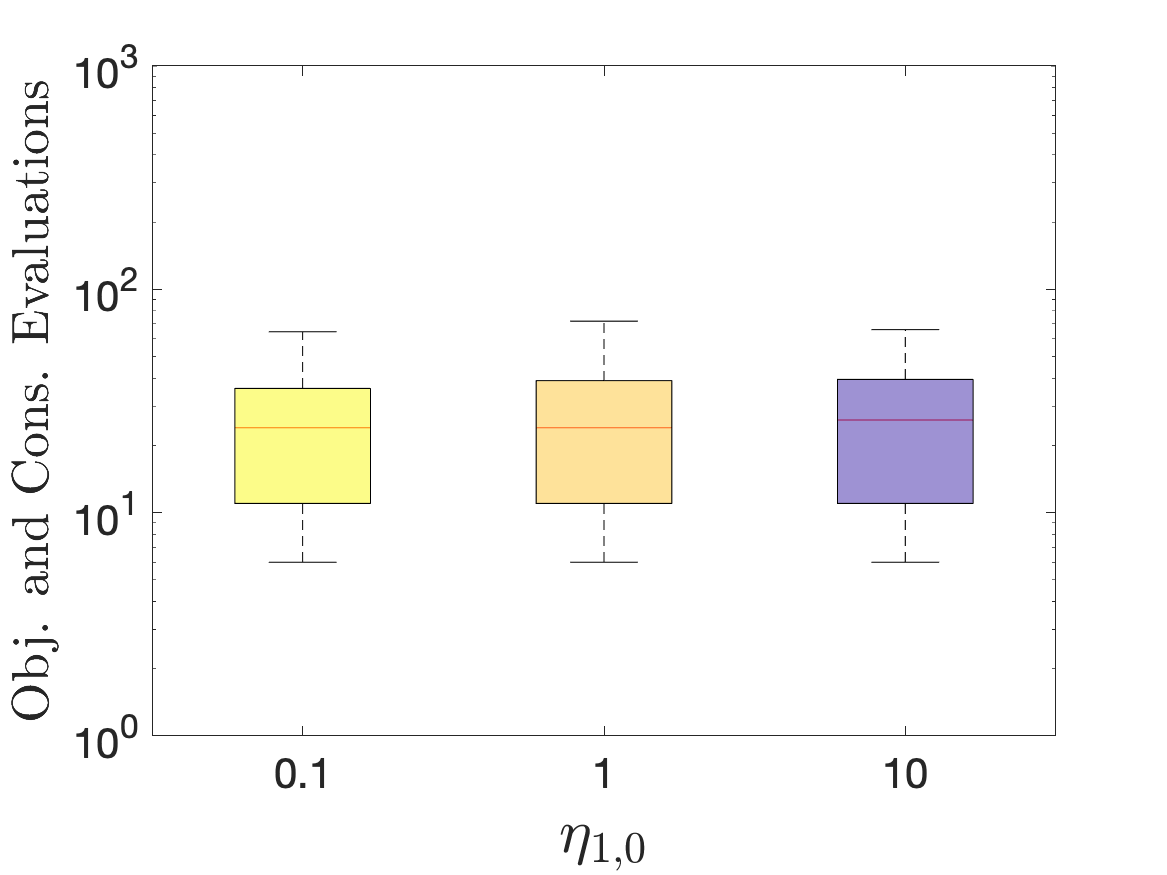}}
		\end{minipage}
		\begin{minipage}[ht]{.33\linewidth}
			\centerline{\includegraphics[scale=0.3]{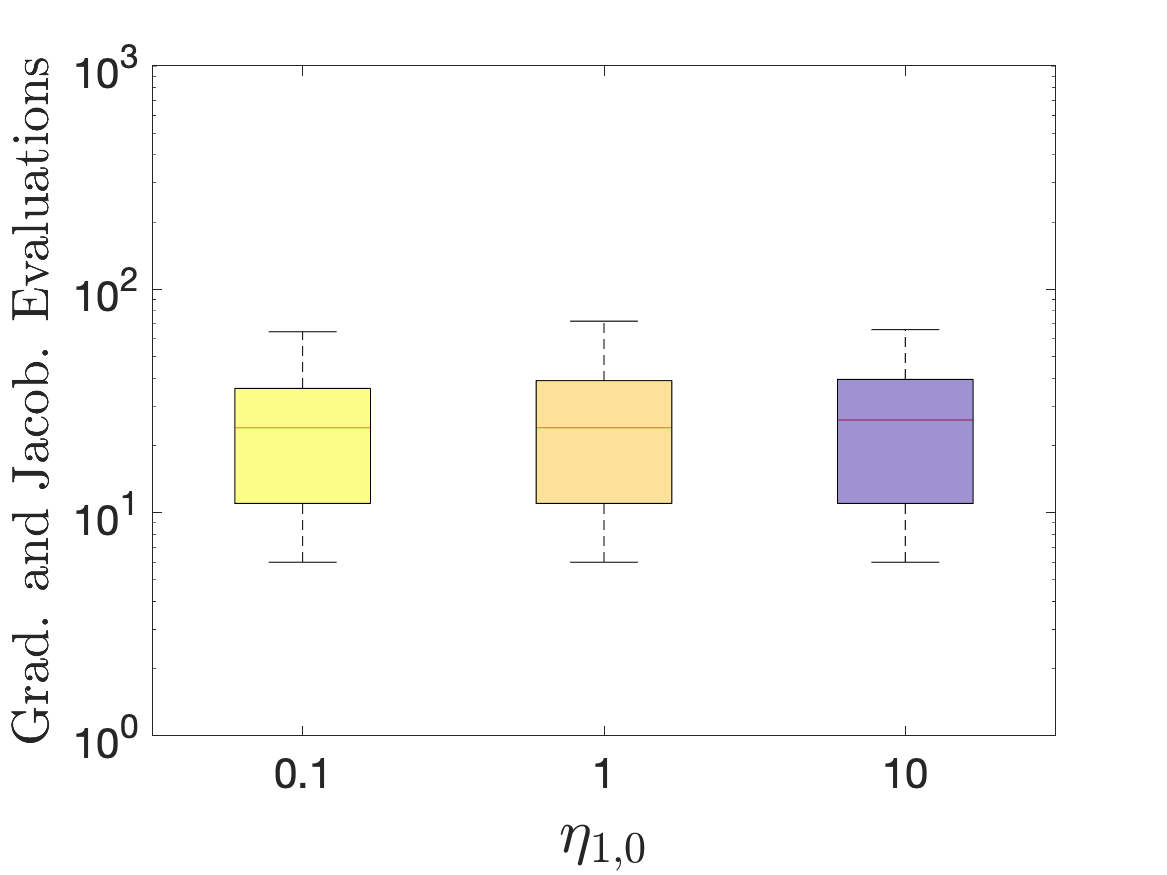}}
		\end{minipage}\\
		\begin{minipage}[ht]{.33\linewidth}
			\centerline{\includegraphics[scale=0.3]{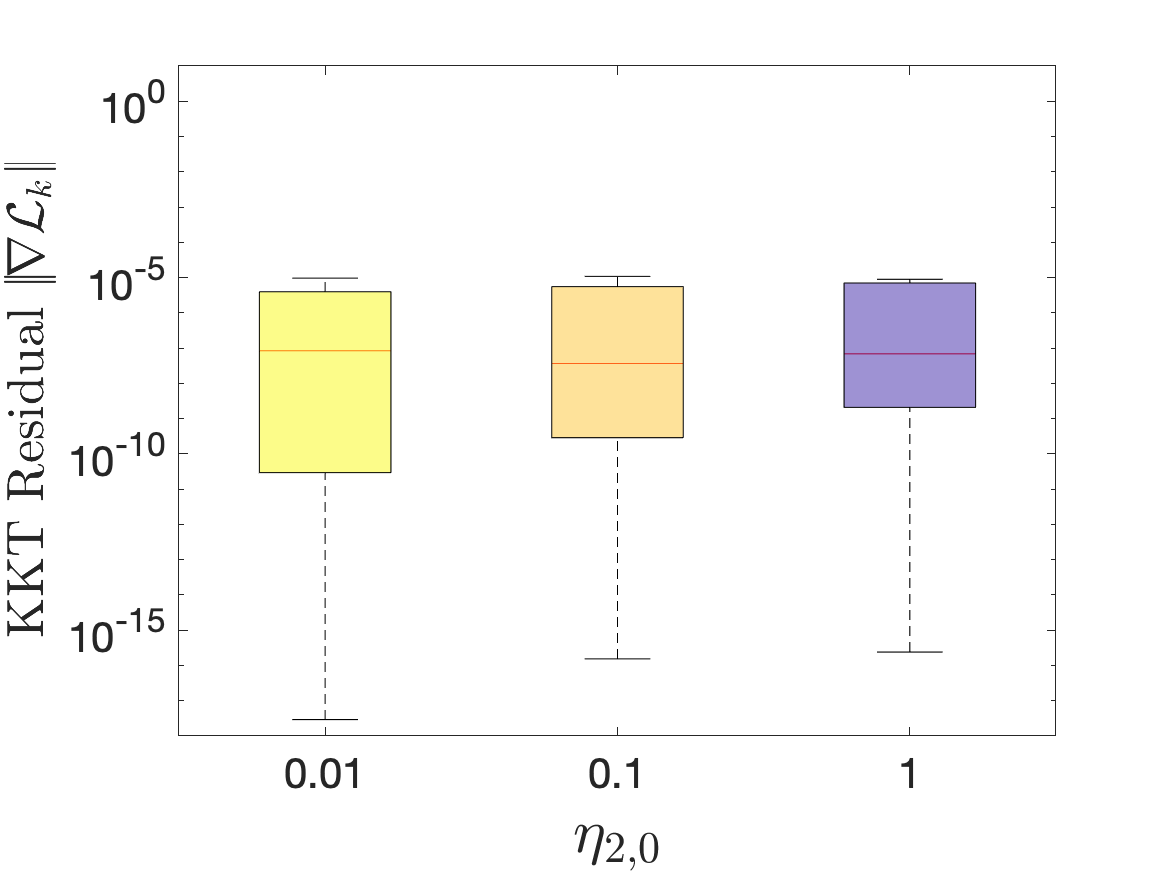}}
		\end{minipage}%
		\begin{minipage}[ht]{.33\linewidth}
			\centerline{\includegraphics[scale=0.3]{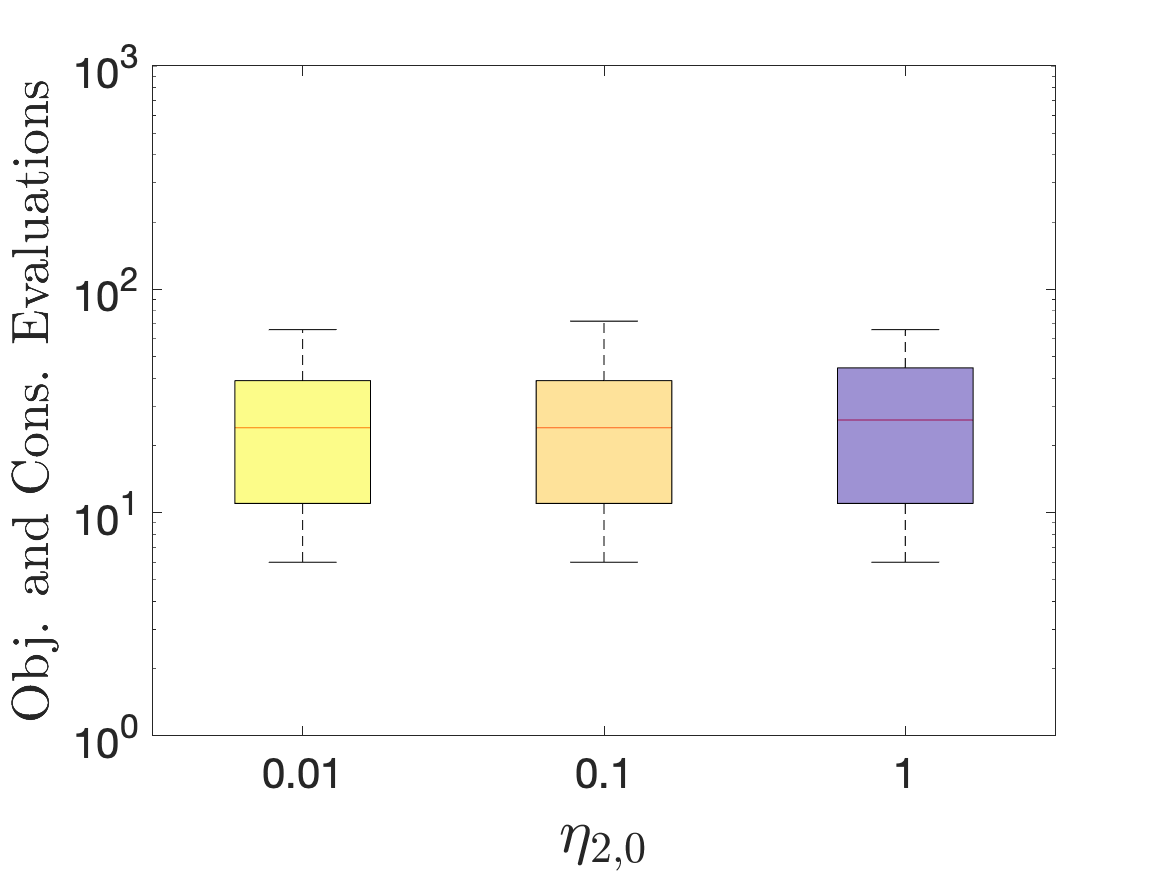}}
		\end{minipage}
		\begin{minipage}[ht]{.33\linewidth}
			\centerline{\includegraphics[scale=0.3]{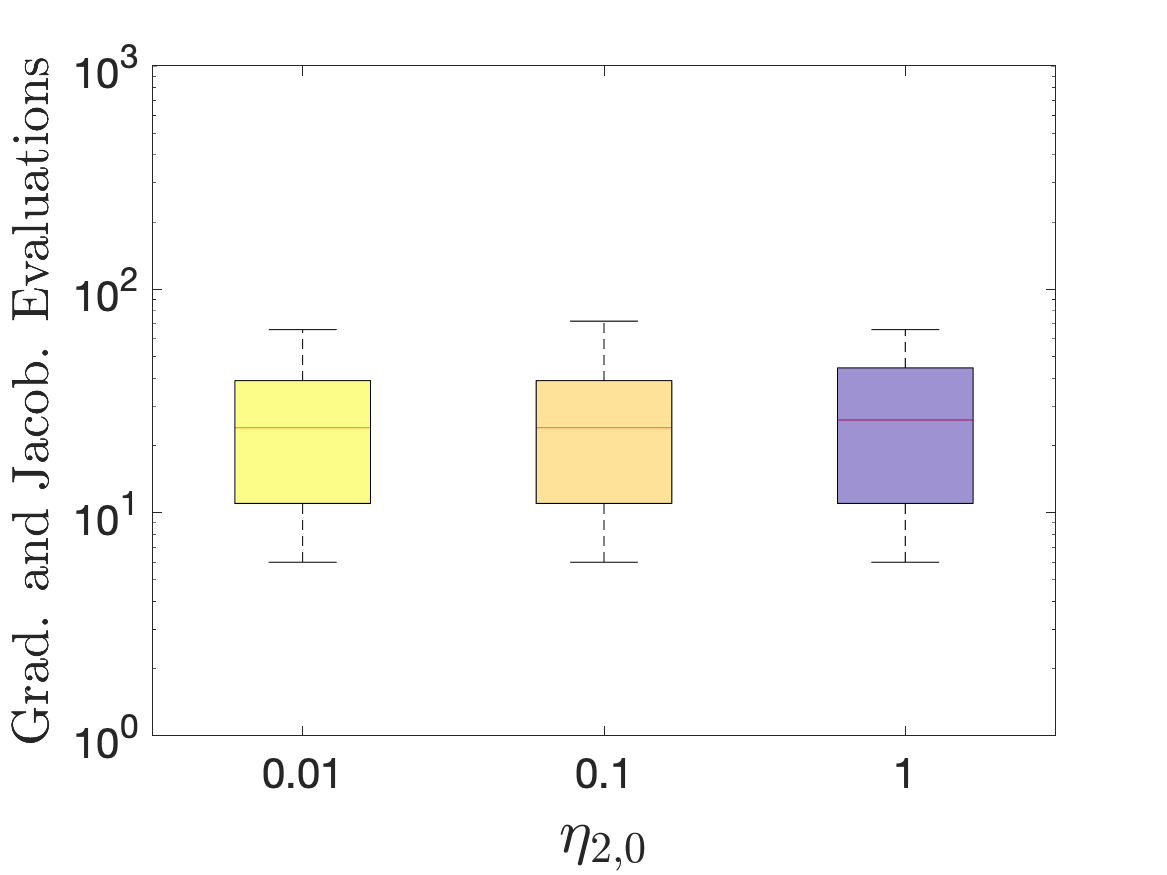}}
		\end{minipage}\\
		\begin{minipage}[ht]{.33\linewidth}
			\centerline{\includegraphics[scale=0.3]{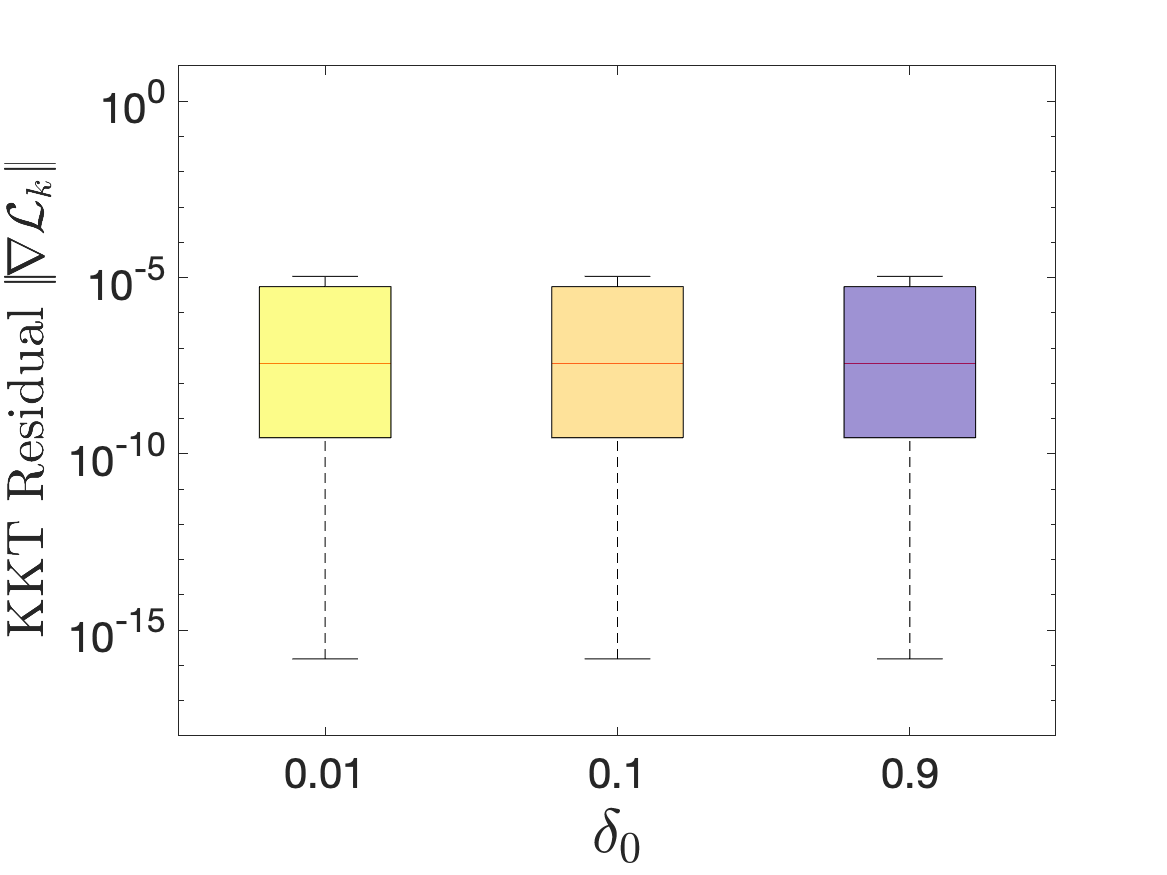}}
		\end{minipage}%
		\begin{minipage}[ht]{.33\linewidth}
			\centerline{\includegraphics[scale=0.3]{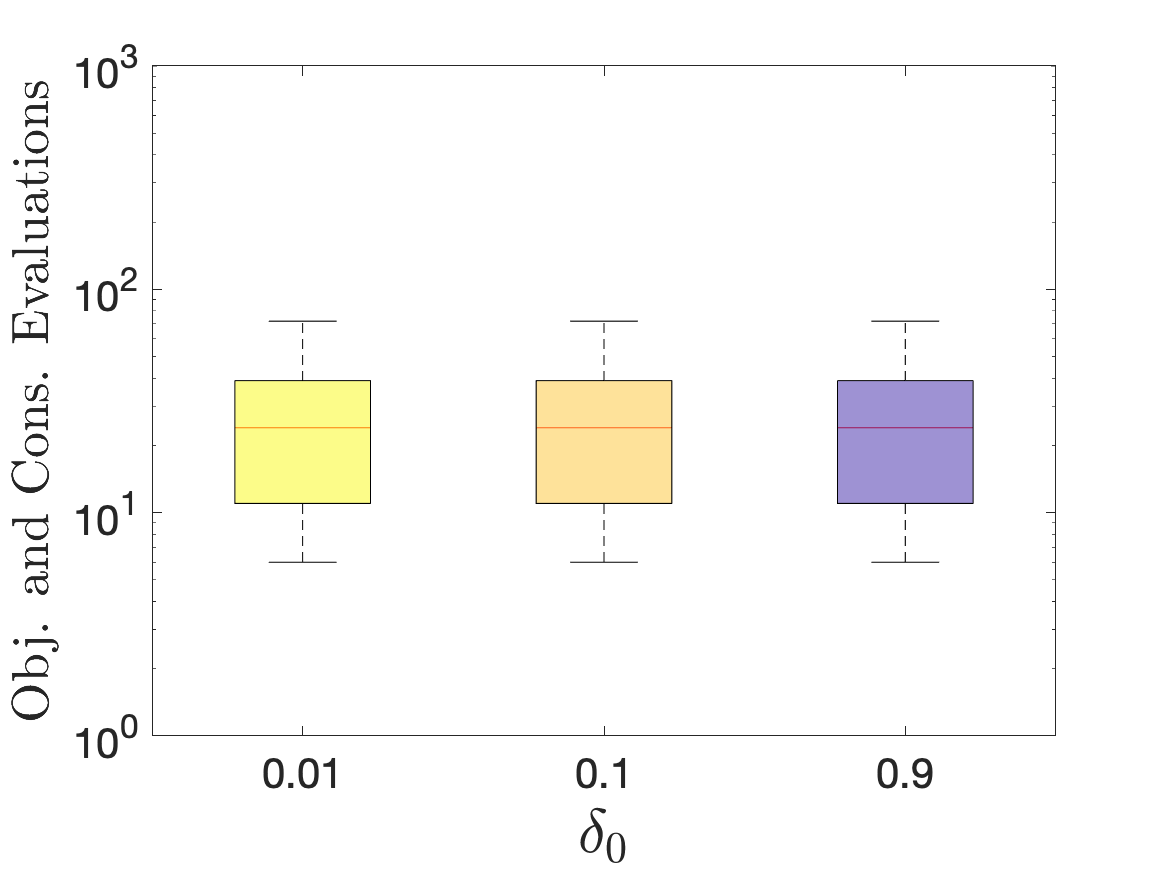}}
		\end{minipage}
		\begin{minipage}[ht]{.33\linewidth}
			\centerline{\includegraphics[scale=0.3]{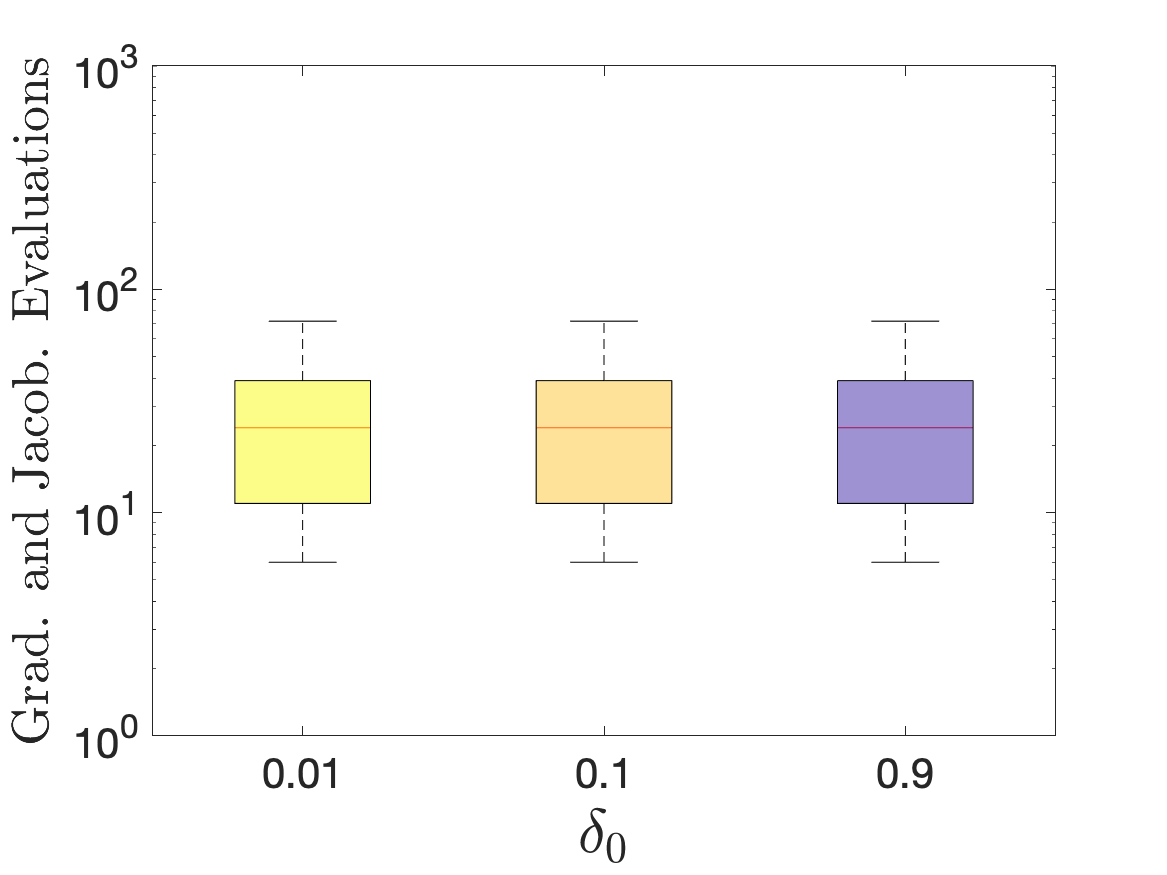}}
		\end{minipage}
		\begin{minipage}[ht]{.33\linewidth}
			\centerline{\includegraphics[scale=0.3]{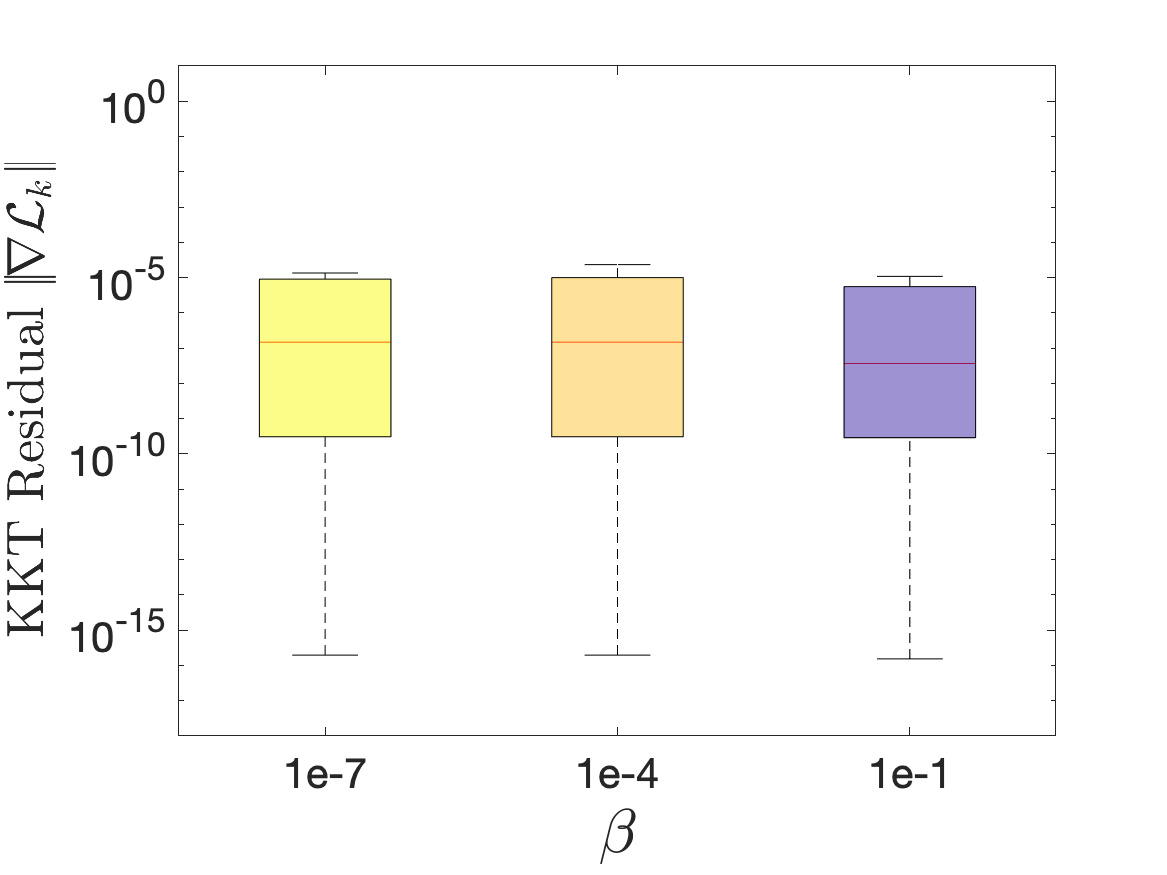}}
		\end{minipage}%
		\begin{minipage}[ht]{.33\linewidth}
			\centerline{\includegraphics[scale=0.3]{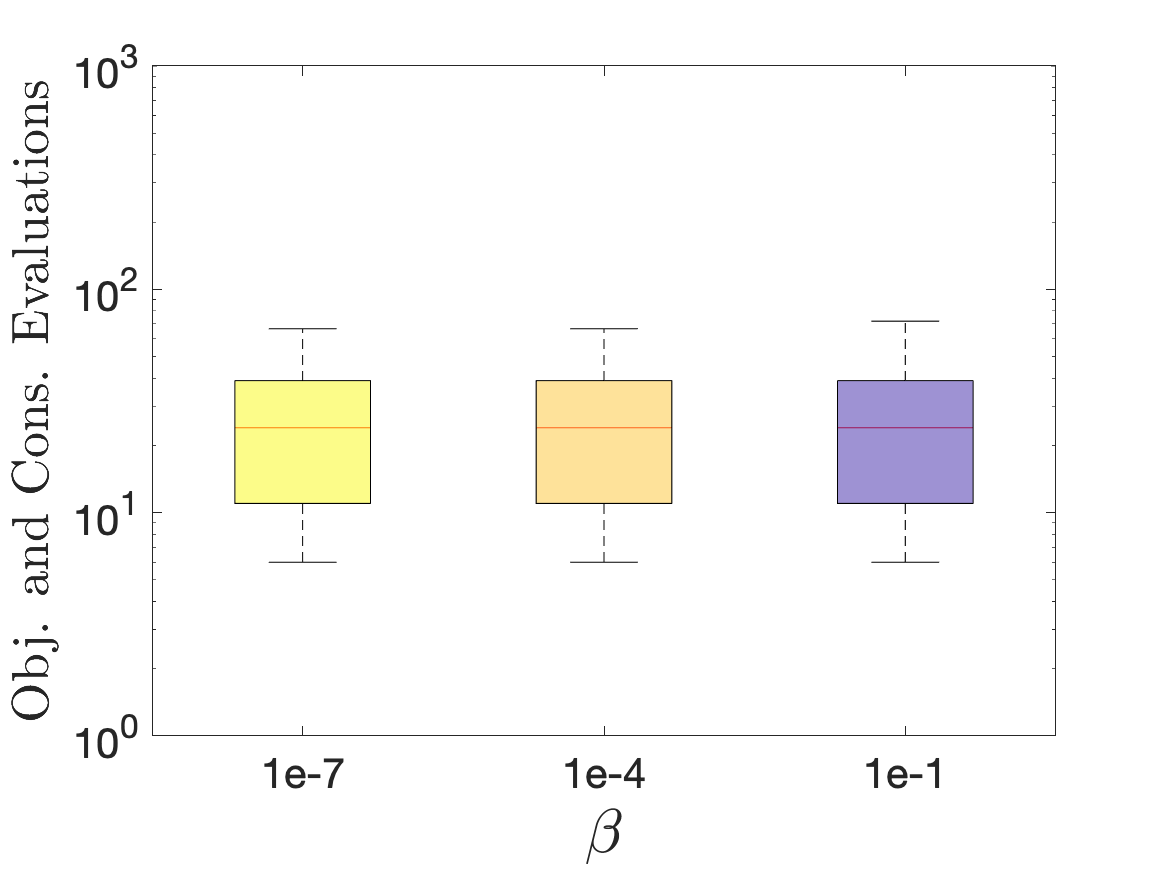}}
		\end{minipage}
		\begin{minipage}[ht]{.33\linewidth}
			\centerline{\includegraphics[scale=0.3]{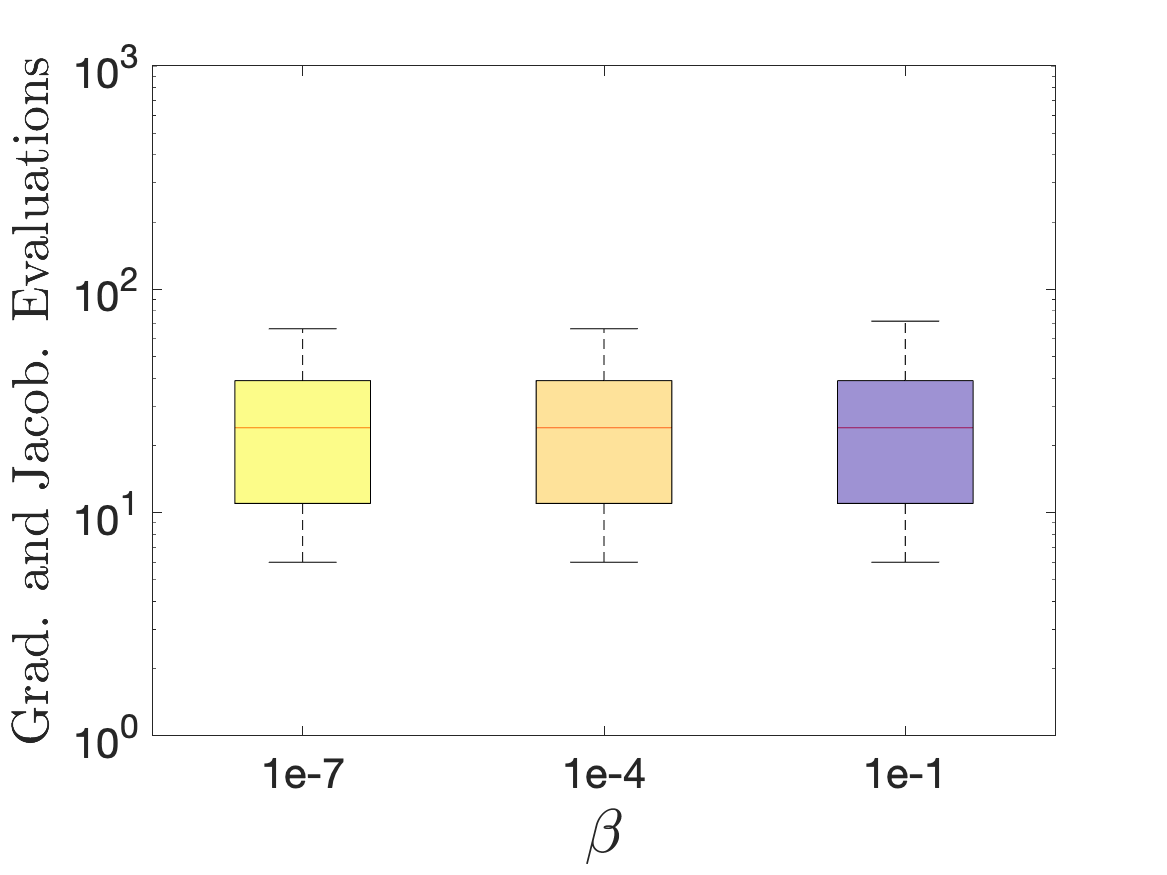}}
		\end{minipage}
		\caption{The boxplots of the KKT residual, the number of objective and constraints evaluations, and the number of gradient and Jacobian evaluations for $\myalg$-GV with different settings of the tuning parameters  $(\eta_{1,0},\eta_{2,0},\delta_0,\beta)$ on CUTEst problems.}\label{fig:fig4}
	\end{center}
\end{figure*}

 Figure~\ref{fig:fig3},
 Table~\ref{tab:table2}, and
 Figure~\ref{fig:fig4}
 present additional results.

\begin{algorithm}[H]
\caption{\citet{Byrd2008Inexact} with the $\ell_1$ penalized merit function}\label{alg:alg2}
\begin{algorithmic}[1]
\STATE {\bfseries Input:} initial iterate $\vz_0$; scalars $\xi_B,\pi_0,\beta, \kappa>0$, $\kappa_1,\epsilon,\tau,\eta\in(0,1)$;
\FOR{$k=0,1,2,\dots$}
\STATE Compute $f_k$, $\nablaxf$, $c_k$, $G_k$, $H_k$, and generate $B_k$;
\STATE Set $\incztk\gets\vzero$ and compute $\vr_k$ by \eqref{eq:residual};
\WHILE{\textit{Termination Test 1} AND \textit{Termination Test 2} does not hold}
\STATE Update $\incztk$ and $\vr_k$ by GMRES;
\ENDWHILE
\IF {\textit{Termination Test 2}  is satisfied and \textit{Model Reduction Condition} does not hold}
\STATE Set $\pi_k\gets\pi_k^{\text{trial}}+10^{-4}$;
\ENDIF
\STATE Select $\alpha_k$ to satisfy the Armijo condition and update the iterate by~\eqref{eq:iterate_update};
\STATE	Set $\pi_{k+1}\gets \pi_{k}$;
\ENDFOR
\end{algorithmic}
\end{algorithm}

\begin{algorithm}[H]
\caption{A modified \citet{Byrd2008Inexact} scheme with adaptive design}\label{alg:alg3}
\begin{algorithmic}[1]
\STATE {\bfseries Input:} initial iterate $\vz_0$; scalars $\xi_B,\pi_0,\kappa_0>0$, $\nu>1$, $\eta\in(0,1)$;
\FOR{$k=0,1,2,\dots$}
\STATE Compute $f_k$, $\nablaxf$, $c_k$, $G_k$, $H_k$, and generate $B_k$;
\STATE Set $\incztk\gets\vzero$ and compute $\vr_k$ by~\eqref{eq:residual}
\WHILE{\textit{Termination Test 1} does not hold}
\WHILE{$\normr_1>\kappa_k\normnablal_1$}
\STATE Update $\incztk$ and $\vr_k$ by GMRES;
\ENDWHILE
\IF {\textit{Model Reduction Condition} does not hold}
\STATE Set $\pi_k\gets\pi_k\nu$ and $\kappa_k\gets\kappa_k/\nu^2$;
\ENDIF
\ENDWHILE
\STATE Select $\alpha_k$ to satisfy the Armijo condition and update the iterate by~\eqref{eq:iterate_update};
\STATE	Set $\pi_{k+1}\gets \pi_{k}$ and $\kappa_{k+1}\gets \kappa_{k}$;
\ENDFOR
\end{algorithmic}
\end{algorithm}

\begin{algorithm}[H]
\caption{Augmented Lagrangian method \citep[Algorithm 17.3]{Nocedal2006Numerical}}\label{alg:alg4}
\begin{algorithmic}[1]
\STATE {\bfseries Input:} initial iterate $\vz_0=(\vx_0,\vlambda_0)$; scalars $\mu_0,\tau_0,\kappa>0$, $\nu_{\mu}>1$, $\nu_{\tau}\in(0,1)$, $\eta\in(0,1)$;
\FOR{$k=0,1,2,\dots$}
\STATE Set $\vx^s_k\gets\vx_k$;
\WHILE{$\|\nabla_{\vx}\mathcal{L}_{\mu_k}(\vx^s_k,\vlambda_k)\|>\tau_k$}
\STATE Compute the gradient and modified Hessian of the augmented Lagrangian at $(\vx^s_k,\vlambda_k)$;
\STATE Find a search direction $\Delta\vx^s_k$ via an inexact Newton subproblem solver with GMRES and a forcing term $\kappa$;
\STATE Select $\alpha_k$ to satisfy the Armijo condition and set $\vx_k^s\gets\vx_k^s+\alpha_k\Delta\vx^s_k$ ;
\STATE Update $\kappa$;
\ENDWHILE
\STATE Set $\vx_{k+1}\gets\vx^s_k$ and $\vlambda_{k+1}\gets\vlambda_k+\mu_kc_k$;
\STATE	Set $\mu_{k+1}\gets\nu_{\mu}\mu_{k}$ and $\tau_{k+1}\gets\nu_{\tau}\tau_{k}$;
\ENDFOR
\end{algorithmic}
\end{algorithm}

\begin{table*}[t]
\caption{Numerical Results for PDE-constrained Problem}\label{tab:table2}
\begin{center}
\begin{small}
\begin{sc}
\begin{tabular}{lccc}
\toprule
Method & KKT residual & obj. and cons. eval & grad. and jacob. eval \\ 
\midrule
AdaSketch-Newton-GV & 1.79e-5  & 18 & 10 \\
AdaSketch-Newton-RK & 4.22e-6  & 14 & 8 \\ 
Algorithm 3-GMRES & 9.10e-5 & 42 & 12  \\ 
Algorithm 2-GMRES & 9.95e-5 & 130 & 34  \\ 
Augmented Lagrangian & 9.99e-5 & 278 & 73 \\ 
\bottomrule
\end{tabular}
\end{sc}
\end{small}
\end{center}
\end{table*}

\end{document}